\documentclass[reqno]{amsart}

\usepackage{pstricks}
\usepackage{verbatim}
\usepackage{amssymb}
\usepackage{mathrsfs}
\usepackage{dsfont}

\AtBeginDocument{%
}


\title{Buildings and Hecke Algebras}

\numberwithin{equation}{section}

\newtheorem{lem}{Lemma}[section]
\newtheorem{thm}[lem]{Theorem}
\newtheorem{cor}[lem]{Corollary}
\newtheorem{prop}[lem]{Proposition}

\theoremstyle{definition}
\newtheorem{defn}[lem]{Definition}

\newtheorem{rem}[lem]{Remark}

\newtheorem{exmp}[lem]{Example}
\newtheorem*{ack}{Acknowledgements}

\newcommand{\ooo}{\mathds{1}}
\newcommand{\oo}{\mathds{1}_0}

\renewcommand{\b}{\beta}

\renewcommand{\S}{\Sigma}

\newcommand{\la}{\lambda}

\newcommand{\s}{\sigma}
\newcommand{\ca}{\mathcal{A}}

\newcommand{\cc}{\mathcal{C}}

\newcommand{\ch}{\mathcal{H}}

\newcommand{\cs}{\mathcal{S}}

\newcommand{\bc}{\mathbb{C}}
\newcommand{\bn}{\mathbb{N}}
\newcommand{\br}{\mathbb{R}}

\newcommand{\bz}{\mathbb{Z}}
\newcommand{\bk}{\mathbb{C}}
\newcommand{\bq}{\mathbb{Q}}

\newcommand{\sca}{\mathscr{A}}
\newcommand{\scb}{\mathscr{B}}
\newcommand{\sch}{\mathscr{H}}

\newcommand{\lan}{\langle}
\newcommand{\ran}{\rangle}
\newcommand{\Res}{R}
\newcommand{\Aut}{\mathrm{Aut}}
\newcommand{\Auttr}{\mathrm{Aut_{tr}}}

\newcommand{\scx}{\mathscr{X}}

\newcommand{\st}{\mathrm{st}}

\newcommand{\ts}{\textsection}

\newcommand{\cha}{\alpha^{\vee}}
\newcommand{\chR}{R^{\vee}}
\newcommand{\chb}{\beta^{\vee}}
\newcommand{\bct}{B\widetilde{\vphantom{A}}C}

\DeclareMathAlphabet{\mathpzc}{OT1}{pzc}{m}{it}

\begin{document}

\author{James Parkinson}

\dedicatory{\upshape
  School of Mathematics and Statistics\\
  University of Sydney, NSW 2006\\
  Australia\\[3pt]
  \texttt{jamesp@maths.usyd.edu.au}\\[8pt]
  \today}

\begin{abstract} In this paper we establish a strong connection between
buildings and Hecke algebras by studying two algebras of
averaging operators on buildings. To each locally finite regular
building we associate a natural algebra $\scb$ of chamber set
averaging operators, and when the building is affine we also
define an algebra $\sca$ of vertex set averaging operators. We
show that for appropriately parametrised Hecke algebras $\sch$ and
$\tilde{\sch}$, the algebra $\scb$ is isomorphic to $\sch$ and
the algebra $\sca$ is isomorphic to the centre of $\tilde{\sch}$.
On the one hand these results give a thorough understanding of the
algebras $\sca$ and $\scb$. On the other hand they give a nice
geometric and combinatorial understanding of Hecke algebras, and
in particular of the Macdonald spherical functions and the centre
of affine Hecke algebras. Our results also produce interesting
examples of association schemes and polynomial hypergroups. In
later work we use the results here to study random walks on affine
buildings.
\end{abstract}

\maketitle

\section*{Introduction}

Let $G=PGL(n+1,F)$ where $F$ is a local field, and let
$K=PGL(n+1,\mathcal{O})$, where $\mathcal{O}$ is the valuation
ring of $F$. The space of bi-$K$-invariant compactly supported
functions on $G$ forms a commutative convolution algebra (see
\cite[Corollary~3.3.7]{macsph} for example). Associated to $G$
there is a building $\scx$ (of type $\widetilde{A}_{n}$), and the
above algebra is isomorphic to an algebra $\sca$ of averaging
operators defined on the space of all functions $G/K\to\bc$. In
\cite{C2} it was shown that these averaging operators may be
defined in a natural way using only the geometric and
combinatorial properties of $\scx$, hence removing the group $G$
entirely from the discussion. For example, in the case $n=1$,
$\scx$ is a homogeneous tree and $\sca$ is the algebra generated
by the operator $A_1$, where for each vertex, $(A_1f)(x)$ is
given by the average value of $f$ over the neighbours of~$x$.

In \cite{C2}, using this geometric approach, Cartwright showed
that $\sca$ is a commutative algebra, and that the algebra
homomorphisms $h:\sca\to\bc$ can be expressed in terms of the
classical Hall-Littlewood polynomials of \cite[III, \ts2]{m2}. It
was not assumed that $\scx$ was constructed from a group $G$
(although there always is such a group when $n\geq3$). Although
not entirely realised in \cite{C2}, as a consequence of our work
here we see that the commutativity of the algebra $\sca$ and the
description of the algebra homomorphisms $h:\sca\to\bc$ follow
from the fact that $\sca$ is isomorphic to the centre of an
appropriately parametrised affine Hecke algebra.

One objective of this paper is to put the above observations into
a more general setting. To do so we will demonstrate a close
connection between buildings and Hecke algebras through the
`combinatorial' study of two algebras of averaging operators
associated to buildings. Apart from establishing these important
connections, our results also have applications to the theory of
random walks on buildings, and provides interesting examples of
association schemes and polynomial hypergroups. We will elaborate
on the random walk applications in a later paper, where we
generalise the results in \cite{cw}. Let us briefly describe the
results we give here.

\subsection{Regularity and Parameter Systems} To begin with we consider buildings as certain \textit{chamber
systems}. Thus a \textit{building} $\scx$ is a set $\cc$ of
\textit{chambers} with an associated Coxeter system $(W,S)$ and a
$W$-\textit{distance function} $\delta:\cc\times\cc\to W$. For
each $c\in\cc$ and $w\in W$, define
$\cc_w(c)=\{d\in\cc\mid\delta(c,d)=w\}$. An important assumption
we make throughout is that $\scx$ is \textit{regular}, by which
we mean that for each $s\in S$, $|\cc_s(c)|=|\cc_s(d)|$ for all
$c,d\in\cc$. In a regular building we write $q_{s}=|\cc_s(c)|$,
and we call the set $\{q_s\}_{s\in S}$ the \textit{parameter
system} of the building. In Proposition~\ref{well} we show that
regularity implies the stronger result that
$|\cc_w(c)|=|\cc_w(d)|$ for all $c,d\in\cc$ and $w\in W$, and as
such we define $q_w=|\cc_w(c)|$. In Theorem~\ref{regular} we show that all
\textit{thick} buildings with no rank 2 residues of type
$\widetilde{A}_1$ are regular, generalising
\cite[Proposition~3.4.2]{scharlau}. This shows that regularity is
a very weak hypothesis.

\subsection{The Algebra $\scb$} Let $\scx$ be any (locally finite) regular building.
For each $w\in W$ we define an operator $B_w$, acting on the
space of functions $f:\cc\to \bc$,~by
\begin{align}
\label{short}(B_wf)(c)=\frac{1}{q_w}\sum_{d\in\cc_w(c)}f(d)\qquad\textrm{for
all $c\in\cc$}.
\end{align}
We call these operators \textit{chamber set averaging operators},
and write $\scb$ for the linear span of $\{B_w\}_{w\in W}$ over
$\bc$. Our main result here is Theorem~\ref{usebalg}, where we
show that $\scb$ is isomorphic to a suitably parametrised Hecke
algebra (the parametrisation depending on the parameter system of
the building). This result is a generalisation of results in
\cite[Chapter 6]{garrett} where an analogous algebra is studied
under the assumption that there is a group $G$ (of label
preserving simplicial complex automorphisms) acting
\textit{strongly transitively} on the building. We note that it
is simple to see that all buildings admitting such a group are
regular. However not all regular buildings admit such a group (see
\cite{ronan2} for the $\widetilde{A}_2$, $\widetilde{C}_2$ and
$\widetilde{G}_2$ buildings). Since we only assume regularity,
our results are more general. We require this additional
generality to prove the more difficult results concerning the
algebra $\sca$ of \textit{vertex set operators} in their full
generality. We note that some of our results in
Section~\ref{section3} are proved in \cite{zieschang} using the
quite different language of \textit{association schemes}.

\subsection{The Algebra $\sca$} The latter part of this paper is mainly devoted to the study
of an algebra $\sca$ of \textit{vertex set averaging operators}
associated to locally finite regular \textit{affine} buildings,
and the connections with \textit{affine Hecke algebras}. We
consider the study of $\sca$ to be the main contribution of this
paper. It is a considerably more complicated object than the
algebra $\scb$. Let us give a (simplified) description of this
algebra.

We now consider a building $\scx$ as a certain \textit{simplicial
complex} \cite[Chapter IV]{brown}, and we write $V$ for the
\textit{vertex set} of $\scx$. In Definition~\ref{goodbuilding}
we define a subset $V_P\subseteq V$ of \textit{good} vertices,
which, for the sake of this simplified description, can be
thought of as the \textit{special vertices} of $\scx$.

To each (locally finite regular) affine building we associate a
root system $R$. Let~$P$ be the \textit{coweight} lattice of
$R$ and write $P^+$ for a set of dominant coweights. For
each $x\in V_P$ and $\la\in P^+$ we define
(Definition~\ref{vertices in sigma}) sets $V_{\la}(x)$ in such a
way that $\{V_{\la}(x)\}_{\la\in P^+}$ forms a partition of
$V_P$. In Theorem~\ref{N} we show that regularity implies that
the cardinalities $|V_{\la}(x)|$, $\la\in P^+$, are independent
of the particular $x\in V_P$, and as such we write
$N_{\la}=|V_{\la}(x)|$. For each $\la\in P^+$ we define an
averaging operator $A_{\la}$, acting on the space of functions $f:V_P\to
\bc$,~by
\begin{align}
\label{short2}(A_{\la}f)(x)=\frac{1}{N_{\la}}\sum_{y\in
V_{\la}(x)}f(y)\qquad\textrm{for all $x\in V_P$}.
\end{align}
These operators specialise to the operators studied in \cite{C2}
when $\scx$ is an $\widetilde{A}_n$ building.

We write $\sca$ for the linear span of $\{A_{\la}\}_{\la\in P^+}$
over $\bc$. Our first main result concerning $\sca$ is
Theorem~\ref{maina}, where we show that $\sca$ is a commutative
algebra. We stress that we only assume regularity, and do not
require the existence of groups or $\mathrm{BN}$-pairs associated
with the building. This puts our results in a very general
setting.

To get a feel for the above definitions in a special case, let
$\scx$ be a homogeneous tree with degree $q+1$, which is a special
case of an $\widetilde{A}_1$ building. Let
$R=\{\alpha,-\alpha\}$, where $\alpha=e_1-e_2$, be the usual root
system of type $A_1$ in the vector space $E=\{x\in\br^2\mid\lan
x,e_1+e_2\ran=0\}$. Taking $\{\alpha\}$ as a base of $R$ we have
$P^+=\{\frac{k}{2}\alpha\}_{k\in\bn}$ where $\bn=\{0,1,\ldots\}$.
Here $V_P=V$, the set of all vertices, and, writing $V_{k}(x)$ in
place of $V_{\la}(x)$ when $\la=\frac{k}{2}\alpha$ with $k\geq0$,
we see that $V_{k}(x)$ is the set of vertices of distance $k$ from
$x$. Thus we compute $N_{k}=1$ if $k=0$ and $(q+1)q^{k-1}$ if
$k\geq1$. The algebra $\sca$ in this case is a well known
object (see \cite{FN} for example). It is generated by $A_1$,
where $(A_1f)(x)=\frac{1}{q+1}\sum_{y\sim x} f(y)$ and the sum is
over the neighbours of $x$.

Our results on the algebra $\sca$ give interesting examples of
\textit{association schemes} (see Remark~\ref{IDRG} and
Remark~\ref{AS}) which generalises the well known construction of
association schemes from \textit{infinite distance regular
graphs}.

\begin{rem} To increase the readability of this paper we have
restricted our attention to \textit{irreducible} affine
buildings. Everything we do here goes through perfectly well for
reducible affine buildings too, and the details will be given
elsewhere. Put briefly, when $\scx$ is a reducible building, it has
a natural description as a \textit{polysimplicial complex}, and by
associating a reducible root system to $\scx$ we can define the
algebra $\sca$ as in the irreducible case. It turns out that
$\scx$ decomposes (essentially uniquely) into the cartesian
product of certain \textit{irreducible components}
$\{\scx_j\}_{j=1}^{k}$, each of which is an irreducible building.
The results of this paper can be used on each irreducible
component $\scx_j$, thus obtaining a family
$\{\sca_j\}_{j=1}^{k}$ of algebras. It turns out that
$\sca\cong\sca_1\times\cdots\times\sca_k$, where $\times$ is
\textit{direct product}.
\end{rem}

\subsection{Connections with Affine Hecke Algebras} The main result of this paper
is Theorem~\ref{i}, where we considerably strengthen the
commutativity result of Theorem~\ref{maina} by showing that $\sca$
is isomorphic to the centre of an appropriately parametrised
\textit{affine Hecke algebra} (the parametrisation depending on
the parameters of the building). Let us briefly describe this
important isomorphism. Let $\tilde{\sch}$ be an \textit{affine
Hecke algebra}, and write $Z(\tilde{\sch})$ for the centre of
$\tilde{\sch}$. It is well known that $Z(\tilde{\sch})$ equals $
\bc[P]^{W_0}$, the algebra of $W_0$-invariant elements of the
group algebra of $P$ (here $P$ is considered as a multiplicative
group in exponential notation $\la\leftrightarrow x^{\la}$). For
$\la\in P^+$ let $P_{\la}(x)$ denote the \textit{Macdonald
spherical function}. This is a special element of $\bc[P]^{W_0}$
which arises naturally in connection with the \textit{Satake
isomorphism}. The isomorphism in Theorem~\ref{i} is then
$A_{\la}\mapsto P_{\la}$.

This isomorphism serves two purposes. Firstly it gives us an
essentially complete understanding the algebra $\sca$. For
example, in Theorem~\ref{generators} we use rather simple facts
about the Macdonald spherical functions to show that $\sca$ is
generated by $\{A_{\la_i}\}_{i\in I_0}$ where $\{\la_i\}_{i\in
I_0}$ is a set of fundamental coweights of $R$. On the other
hand, since $\sca$ is a purely combinatorial object, the above
isomorphism gives a nice combinatorial description of
$Z(\tilde{\sch})$ when a suitable building exists. In particular
the \textit{structure constants} $c_{\la,\mu;\nu}$ that appear in
$$P_{\la}(x)P_{\mu}(x)=\sum_{\nu\in
P^+}c_{\la,\mu;\nu}P_{\nu}(x)\qquad\textrm{are}\qquad
c_{\la,\mu;\nu}=\frac{N_{\nu}}{N_{\la}N_{\mu}}|V_{\la}(x)\cap
V_{\mu^*}(y)|,$$ for some $\mu^*\in P^+$ (depending only on $\mu$
in a simple way). This shows that (when a suitable building
exists) $c_{\la,\mu;\nu}\geq0$.

In Theorem~\ref{hyp} we extend this result by showing that the
$c_{\la,\mu;\nu}$'s are (up to positive normalisation factors)
polynomials in the variables $\{q_s-1\}_{s\in S}$ with nonnegative
integer coefficients (even when no building exists). This
generalises the main theorem in \cite{mm}, where the corresponding
result for the $A_n$ case (where the $c_{\la,\mu,\nu}$'s are
certain \textit{Hall polynomials}) is proved. Thus we see how to
construct a polynomial hypergroup from the structure constants
$c_{\la,\mu,\nu}$ as in \cite{bloom} (see also \cite{hyper}).

Since the submission of this paper we have learnt that
Theorem~\ref{hyp} has been proved independently by Schwer in
\cite{schwer}, where a formula for $c_{\la,\mu;\nu}$ is given.

In later papers we will use our results here to give a description
of the algebra homomorphisms $h:\sca\to\bc$ in terms of the
Macdonald spherical functions of \cite[Chapter 4]{macsph}. We
will also provide an integral formula for these algebra
homomorphisms (over the \textit{boundary of $\scx$}), and use
these results to study local limit theorems, central limit
theorems and rate of escape theorems for \textit{radial} random
walks on affine buildings.

\begin{ack} The author would like to thank Donald Cartwright for helpful
discussions and suggestions throughout the preparation of this
paper. Thank you also to Jon Kusilek for useful discussions
regarding affine Hecke algebras.
\end{ack}

\section{Coxeter Groups, Chamber Systems and Buildings}\label{s1} Let $I$ be an
index set, which we assume throughout is finite, and for $i,j\in I$
let $m_{i,j}$ be an integer or $\infty$ such that
$m_{i,j}=m_{j,i}\geq2$ for all $i\neq j$, and $m_{i,i}=1$ for all
$i\in I$. We call $M=(m_{i,j})_{i,j\in I}$ a \textit{Coxeter
matrix}. The \textit{Coxeter group} of type $M$ is the group
\begin{align}
W=\lan \{s_i\}_{i\in I}\mid(s_is_j)^{m_{i,j}}=1\textrm{ for all
$i,j\in I$}\ran,
\end{align}
where the relation $(s_is_j)^{m_{i,j}}=1$ is omitted if
$m_{i,j}=\infty$. Let $S=\{s_i\mid i\in I\}$. For subsets
$J\subset I$ we write $W_J$ for the subgroup of $W$ generated by
$\{s_i\}_{i\in J}$. Given $w\in W$, we define the \textit{length}
$\ell(w)$ of $w$ to be smallest $n\in\bn$ such that
$w=s_{i_1}\ldots s_{i_n}$, with $i_1,\ldots,i_n\in I$.

It will be useful on occasion to work with $I^*$, the free monoid
on $I$. Thus elements of $I^*$ are \textit{words} $f=i_1\cdots
i_n$ where $i_1,\ldots, i_n\in I$, and we write
$s_{f}=s_{i_1}\cdots s_{i_n}\in W$. Recall \cite[Chapter 2,
\ts1]{ronan} that an \textit{elementary homotopy} is an alteration
from a word of the form $f_1p(i,j)f_2$ to a word of the form
$f_1p(j,i)f_2$, where $p(i,j)=\cdots ijij$ ($m_{i,j}$ terms). We
say that the words $f$ and $f'$ are \textit{homotopic} if
$f$ can be transformed into $f'$ by a sequence of elementary
homotopies, in which case we write $f\sim f'$. A word $f$ is said
to be \textit{reduced} if it is not homotopic to a word of the
form $f_1iif_2$ for any $i\in I$. Thus $f=i_1\cdots i_n\in I^*$ is
reduced if and only if $s_f=s_{i_1}\cdots s_{i_n}$ is a reduced
expression in $W$ (that is, $\ell(s_f)=n$).

The \textit{Coxeter graph}\index{Coxeter graph} of $W$ is the graph
$D=D(W)$ with vertex set $I$, such that vertices $i,j\in I$ are
joined by an edge if and only if $m_{i,j}\geq3$. If $m_{i,j}\geq4$
then the edge $\{i,j\}$ is labelled by $m_{i,j}$.

By an \textit{automorphism of $D$} we mean a permutation of the
vertex set of $D$ that preserves adjacency and edge labels, that
is, a permutation $\s$ of~$I$ such that $m_{\s(i),\s(j)}=m_{i,j}$
for all $i,j\in I$. We write $\Aut(D)$ for the group of all
automorphisms of $D$.

An automorphism $\s$ of $D$ induces a group automorphism of $W$,
which we will also denote by $\s$, via the (well defined) action
\begin{align}
\label{sigmaaction}\s(w)=s_{\s(i_1)}\cdots s_{\s(i_n)}
\end{align}
whenever $s_{i_1}\cdots s_{i_n}$ is an expression for $w$. Note
that $\ell(\s(w))=\ell(w)$ for all $w\in W$.

Recall \cite[p.1]{ronan} that a set $\cc$ is a \textit{chamber
system over a set $I$} if each $i\in I$ determines a partition of
$\cc$, two elements in the same block of this partition being called
\textit{$i$-adjacent}. The elements of $\cc$ are called
\textit{chambers}, and we write $c\sim_i d$ to mean that the
chambers $c$ and $d$ are $i$-adjacent. By a \textit{gallery} of
type $i_1\cdots i_n\in I^*$ in $\cc$ we mean a finite sequence
$c_0,\ldots,c_n$ of chambers such that $c_{k-1}\sim_{i_k}c_k$ and
$c_{k-1}\neq c_k$ for $1\leq k\leq n$. If
$J\subseteq I$, we say that $d\in\cc$ is \textit{$J$-connected}
to $c\in\cc$ if $d$ can be joined to $c$ by a gallery
$c=c_0,\ldots,c_n=d$ of type $j_1\cdots j_n$ with each $j_k\in
J$. We call such a gallery a \textit{$J$-gallery}, and for
$c\in\cc$ we write $R_{J}(c)$ for the set of all chambers that
can be joined to $c$ by a $J$-gallery. We call $R_J(c)$ the
\textit{$J$-residue of $c$}. If $\mathcal{C}$ and $\mathcal{D}$
are chamber systems over a common index set $I$, we call a map
$\psi:\cc\to\mathcal{D}$ an \textit{isomorphism of chamber
systems} if $\psi$ is a bijection such that $c\sim_i d$ if and only if $\psi(c)\sim_i\psi(d)$.

To a Coxeter group $W$ over the index set $I$ we associate a
chamber system $\cc(W)$, called the \textit{Coxeter complex} of
$W$, by taking the elements $w\in W$ as chambers, and for each
$i\in I$ define $i$-adjacency by declaring $w\sim_i w$ and
$w\sim_i ws_i$.

For the present purpose it is most convenient to consider
\textit{buildings} as certain chamber systems. Thus we give the
definition of buildings from \cite{ronan}.

\begin{defn} \cite{ronan}\textbf{.}\label{def1} Let $M$ be the Coxeter matrix of a
Coxeter group $W$ over~$I$. Then $\scx$ is a \textit{building of
type $M$} if
\begin{enumerate}
\item[$\mathrm{(i)}$] $\scx$ is a chamber system over $I$ such that for each
$c\in\scx$ and $i\in I$, there is a chamber $d\neq c$ in $\scx$
such that $d\sim_i c$, and
\item[$\mathrm{(ii)}$] there exists a \textit{$W$-distance function}
$\delta:\scx\times\scx\to W$ such that if $f$ is a reduced word
then $\delta(c,d)=s_f$ if and only if $c$ and $d$ can be joined
by a gallery of type $f$.
\end{enumerate}
\end{defn}

We will always use the symbol $\scx$ to denote a building. It is
convenient to write $\cc=\cc(\scx)$ for the chamber set of
$\scx$, even though according to the above definition $\scx$ is
itself a set of chambers. We sometimes say that $\scx$ is a
building of type~$W$ if~$W$ is the Coxeter group of type $M$. A
building $\scx$ is said to be \textit{thick} if for each
$c\in\cc$ and $i\in I$ there exist at least two distinct chambers
$d\neq c$ such that $d\sim_i c$. The \textit{rank} of a building
of type $M$ is the cardinality of the index set $I$. We sometimes
call a building \textit{irreducible} if the associated Coxeter
group is irreducible (that is, has connected Coxeter graph).

\section{Regularity and Parameter Systems}\label{section2}

In this section we write $\scx$ for a building of type $M$, with
associated Coxeter group $W$ over index set $I$. We will assume
that $\scx$ is \textit{locally finite}, by which we mean
$|I|<\infty$ and $|\{b\in\cc\mid a\sim_i b\}|<\infty$ for all
$i\in I$ and $a\in\cc$.

For each $a\in\cc$ and $w\in W$, let
\begin{align}
\label{chamber sets}\cc_w(a)=\{b\in\cc\mid\delta(a,b)=w\}\,.
\end{align}
Observe that for each fixed $a\in\cc$, the family
$\{\cc_w(a)\}_{w\in W}$ forms a partition of $\cc$.

We say that $\scx$ is \textit{regular} if for each $s\in S$,
$|\cc_s(a)|$ is independent of $a\in\cc$. If $\scx$ is a regular
building we define $q_s=|\cc_s(a)|$ for each $s\in S$ (this is
independent of $a\in\cc$ by definition), and we call
$\{q_s\}_{s\in S}$ the \textit{parameter system of the building}.
Local finiteness implies that $q_s<\infty$ for all $s\in S$. We
often write $q_i$ in place of $q_{s_i}$ for $i\in I$.

The two main results of this section are
Proposition~\ref{well}(ii), where we give a method for finding
relationships that must hold between the parameters of buildings,
and Theorem~\ref{regular}, where we generalise
\cite[Proposition~3.4.2]{scharlau} and show that all thick
buildings with no rank 2 residues of type $\widetilde{A}_1$ are
regular.

\begin{prop}\label{well} Let $\scx$ be a locally finite regular building.
\begin{enumerate}
\item[$\mathrm{(i)}$] $|\cc_{w}(a)|=q_{i_1}q_{i_2}\cdots q_{i_n}$ whenever
$w=s_{i_1}\cdots s_{i_n}$ is a reduced expression, and
\item[$\mathrm{(ii)}$] $q_i=q_j$ whenever $m_{i,j}<\infty$ is odd.
\end{enumerate}
\end{prop}

\begin{proof} We first prove (i). The result is true when
$\ell(w)=1$ by regularity. We claim that whenever $s=s_i\in S$ and
$\ell(ws)=\ell(w)+1$,
\begin{align}\label{addition1}
\cc_{ws}(a)=\bigcup_{b\in\cc_w(a)}\cc_s(b)
\end{align}
where the union is disjoint, from which the result follows by
induction.

First suppose that $c\in\cc_{ws}(a)$ where $\ell(ws)=\ell(w)+1$.
Then there exists a minimal gallery $a=c_0,\ldots,c_k=c$ of type
$fi$ (where $w=s_f$ with $f\in I^*$ reduced) from $a$ to $c$, and
in particular $c\in\cc_{s}(c_{k-1})$ where $c_{k-1}\in\cc_w(a)$.
On the other hand, if $c\in\cc_s(b)$ for some $b\in\cc_w(a)$ then
$c\in\cc_{ws}(a)$ since $\ell(ws)=\ell(w)+1$, and so we have
equality in (\ref{addition1}). To see that the union is disjoint,
suppose that $b,b'\in\cc_w(a)$ and that
$\cc_s(b)\cap\cc_s(b')\neq\emptyset$. Then if $b'\neq b$ we have
$b'\in\cc_s(b)$, and thus $b'\in\cc_{ws}(a)$, a contradiction.

To prove (ii), suppose $m_{i,j}<\infty$ is odd. Since
$s_is_js_i\cdots s_i=s_js_is_j\cdots s_j$ ($m_{i,j}$ factors on each
side), by (i) we have $q_{i}q_jq_i\cdots q_i=q_jq_iq_j\cdots q_j$
($m_{i,j}$ factors on each side), and the result follows.
\end{proof}

\begin{cor}\label{conjno} Let $\scx$ be a locally finite regular building
of type $W$. If $s_j=ws_iw^{-1}$ for some $w\in W$ then $q_i=q_j$.
\end{cor}

\begin{proof} By \cite[IV, \ts1, No.3, Proposition~3]{bourbaki},
$s_j=ws_iw^{-1}$ for some $w\in W$ if and only if there exists a
sequence $s_{i_1},\ldots,s_{i_p}$ such that $i_1=i$, $i_p=j$, and
$m_{i_k,i_{k+1}}$ is finite and odd for each $1\leq k<p$. The
result now follows from Proposition~\ref{well}(ii).  \end{proof}

Proposition~\ref{well}(i) justifies the notation
$q_w=q_{i_1}\cdots q_{i_n}$ whenever $s_{i_1}\cdots s_{i_n}$ is a
reduced expression for $w$; it is independent of the particular
reduced expression chosen. Clearly we have $q_{w^{-1}}=q_w$ for
all $w\in W$.

\begin{exmp}\label{interesting} Using Proposition~\ref{well}(ii) it is now a simple exercise to
describe the relations between the parameters of any given
(locally finite) regular building. For example, in a building of
type
$\stackrel{4}{\bullet{\rule[0.078cm]{0.8cm}{0.01cm}\bullet}}\stackrel{}{{\rule[0.078cm]{0.8cm}{0.01cm}\bullet}}$
(with the nodes labelled $0,1$ and $2$ from left to right) we
must have $q_{1}=q_{2}$ since $m_{1,2}=3$ is odd. Note that we
cannot relate $q_{0}$ to $q_{1}$ since $m_{0,1}=4$ is even.
\end{exmp}

The following theorem seems to be well known (see
\cite[Proposition~3.4.2]{scharlau} for the case $|W|<\infty$), but
we have been unable to find a direct proof in the literature. For
the sake of completeness we will provide a proof here.

\begin{thm}\label{regular} Let $\scx$ be a thick building such that $m_{i,j}<\infty$ for each pair $i,j\in I$.
Then $\scx$ is regular.
\end{thm}

Before giving the proof of Theorem~\ref{regular} we make some
preliminary observations. First we note that the assumption that
$m_{i,j}<\infty$ in Theorem~\ref{regular} is essential, for
$\widetilde{A}_1$ buildings are not in general regular, as they
are just trees with no end vertices. Secondly we note that
Theorem~\ref{regular} shows that most `interesting' buildings are
regular, for examining the Coxeter graphs of the affine Coxeter
groups, for example, we see that $m_{i,j}=\infty$ only occurs in
$\widetilde{A}_1$ buildings. Thus regularity is not a very
restrictive hypothesis.

Recall that for $m\geq2$ or $m=\infty$ a \textit{generalised
$m$-gon} is a connected bipartite graph with diameter $m$ and
girth $2m$. By \cite[Proposition~3.2]{ronan}, a building of type
$\stackrel{m}{\bullet{\rule[0.078cm]{0.8cm}{0.01cm}\bullet}}$ is
a generalised $m$-gon, and vice versa (where the edge set of the
$m$-gon is taken to be the chamber set of the building, and vice
versa).

In a generalised $m$-gon we define the \textit{valency} of a
vertex $v$ to be the number of edges that contain $v$, and we call
the generalised $m$-gon \textit{thick} if every vertex has
valency at least 3. By \cite[Proposition~3.3]{ronan}, in a thick
generalised $m$-gon with $m<\infty$, vertices in the same
partition have the same valency.  In the statement of
\cite[Proposition~3.3]{ronan}, the assumption $m<\infty$ is
inadvertently omitted. The result is in fact false if $m=\infty$,
for a thick generalised $\infty$-gon is simply a tree in which
each vertex has valency at least~3.

\begin{proof}[Proof of Theorem~\ref{regular}] For each $a\in \cc$ and each $i\in I$, let $q_i(a)=|\cc_{s_i}(a)|$. By thickness,
we have $q_i(a)>1$. We will show that $q_i(a)=q_i(b)$ for all
$a,b\in\cc$ and for all $i\in I$.

Fix $a\in\cc$. By \cite[Theorem~3.5]{ronan} we know that for
$i,j\in I$, the residue $R_{\{i,j\}}(a)$ is a thick building of
type $M_{\{i,j\}}$ which is in turn a thick generalised
$m_{i,j}$-gon by \cite[Proposition~3.2]{ronan}. Thus, since
$m_{i,j}<\infty$ by assumption, \cite[Proposition~3.3]{ronan}
implies that
\begin{align}
\label{reg eqn}q_i(b)=q_i(a)&&\textrm{for all
$b\in\Res_{\{i,j\}}(a)$}\,.
\end{align}

Now, with $a$ fixed as before, let $b\in\cc$ be any other
chamber. Suppose firstly that $b\sim_k a$ for some $k\in I$. If
$k=i$, then $q_i(b)=q_i(a)$ since $\sim_i$ is an equivalence
relation. So suppose that $k\neq i$. Then
\begin{align*}
q_i(b)+1&=|\{c\in\cc:c\sim_i b\}|\\
&=|\{c\in\Res_{\{i,k\}}(b):c\sim_i b\}|\\
&=|\{c\in\Res_{\{i,k\}}(b):c\sim_i a\}|&&\textrm{by (\ref{reg eqn})}\\
&=|\{c\in\Res_{\{i,k\}}(a):c\sim_i a\}|&&\textrm{since
$\Res_{\{i,k\}}(b)=\Res_{\{i,k\}}(a)$}\\
&=|\{c\in\cc:c\sim_i a\}|=q_i(a)+1\,,
\end{align*}
and so $q_i(b)=q_i(a)$. Induction now shows that $q_i(a)$ is
independent of the particular~$a$, and so the building is regular.
 \end{proof}

\begin{rem} The description of parameter systems given in
this section by no means comes close to \textit{classifying} the
parameter systems of buildings. For example, it is an open
question as to whether thick $\widetilde{A}_2$ buildings exist
with parameters that are not prime powers. By the free
construction of certain buildings given in \cite{ronan2} this is
equivalent to the corresponding question concerning the
parameters of projective planes (generalised $3$-gons). See
\cite[Section 6.2]{batten} for a discussion of the known
parameter systems of generalised 4-gons.
\end{rem}

We conclude this section by recording a definition of later
reference.

\begin{defn}\label{defn:poincarep} Let
$\{q_s\}_{s\in S}$ be a set of indeterminates such that
$q_{s'}=q_s$ whenever $s'=wsw^{-1}$ for some $w\in W$. Then
\cite[IV, \ts1, No.5, Proposition~5]{bourbaki} implies that for
$w\in W$, the monomial $q_w=q_{s_{i_1}}\cdots q_{s_{i_n}}$ is
independent of the particular reduced decomposition
$w=s_{i_1}\cdots s_{i_n}$ of $w$. If $U$ is a finite subset of
$W$, the \textit{Poincar\'{e} polynomial} $U(q)$ of $U$ is
\begin{align*}
U(q)=\sum_{w\in U}q_w\,.
\end{align*}
Usually the set $\{q_s\}_{s\in S}$ will be the parameters of a
building (see Corollary~\ref{conjno}).
\end{defn}

\section{Chamber Set Operators and Chamber
Regularity}\label{section3}

The results of this section generalise the results in
\cite[Chapter 6]{garrett}, where it is assumed that there is a
group $G$ (of type preserving simplicial complex automorphisms)
acting \textit{strongly transitively} on $\scx$ (see
\cite[\ts5.2]{garrett}). As noted in the introduction, all
buildings admitting such a group action are necessarily regular,
whereas the converse is not true. Our proofs work for all locally
finite regular buildings, which, by Theorem~\ref{regular},
includes all thick buildings with no rank 2 residues of type
$\widetilde{A}_1$. It should be noted that our results also apply
to thin buildings (where $q_i=1$ for all $i\in I$), as well as to
regular buildings that are neither thick nor thin (that is,
buildings that have $q_i=1$ for some but not all $i\in I$). We
note that some of the results of this section are proved in
\cite{zieschang} in the context of \textit{association schemes}.

Let $\scx$ be a locally finite regular building. We say that
$\scx$ is \textit{chamber regular} if for all $w_1$ and $w_2$ in
$W$,
$$|\cc_{w_1}(a)\cap\cc_{w_2}(b)|=|\cc_{w_1}(c)\cap\cc_{w_2}(d)|\quad\textrm{whenever}\quad\delta(a,b)=\delta(c,d)\,,$$
where the sets $\cc_w(a)$ are as in (\ref{chamber sets}). In this
section we will prove that regularity implies chamber regularity
(Corollary \ref{cor1}), and we introduce an algebra $\scb$ of
chamber set averaging operators (Definition \ref{balg}) and show
that this algebra is isomorphic to a suitably parametrised Hecke
algebra (Proposition~\ref{usebalg}). Recall that for a regular
building we define $q_s=|\cc_s(a)|$, and we write $q_{s_i}=q_i$.

\begin{defn}\label{chamberops} For each $w\in W$, define an operator
$B_w$, acting on the space of all functions $f:\cc\to\bc$ as in
(\ref{short})
\end{defn}

Observe that $b\in \cc_w(a)$ if and only if $a\in
\cc_{w^{-1}}(b)$. If $\cc'\subseteq\cc$, write
$1_{\cc'}:\cc\to\{0,1\}$ for the characteristic function on
$\cc'$. Thus for $w_1,w_2\in W$ we have
\begin{align}
\nonumber(B_{w_1}B_{w_2}f)(a)&=\frac{1}{q_{w_1}}\sum_{b\in
\cc_{w_1}(a)}(B_{w_2}f)(b)\\
\nonumber&=\frac{1}{q_{w_1}q_{w_2}}\sum_{b\in
\cc_{w_1}(a)}\sum_{c\in
\cc_{w_2}(b)}f(c)\\
\label{one}&=\frac{1}{q_{w_1}q_{w_2}}\sum_{b\in\cc}\sum_{c\in\cc}1_{\cc_{w_1}(a)}(b)1_{\cc_{w_2}(b)}(c)f(c)\\
\nonumber&=\frac{1}{q_{w_1}q_{w_2}}\sum_{c\in\cc}\left(\sum_{b\in\cc}1_{\cc_{w_1^{\vphantom{-1}}}(a)}(b)1_{\cc_{w_2^{-1}}(c)}(b)\right)f(c)\\
\nonumber&=\frac{1}{q_{w_1}q_{w_2}}\sum_{c\in\cc}|\cc_{w_1^{\vphantom{-1}}}(a)\cap
\cc_{w_2^{-1}}(c)|\,f(c)\,.
\end{align}

We wish to explicitly compute the above when $w_2=s\in S$ (and so
$w_2^{-1}=w_2$). Thus we have the following lemmas.

\begin{lem}\label{firstcount} Let $w\in W$ and $s\in S$, and fix $a\in\cc$.
Then
\begin{align*}
\cc_w(a)\cap\cc_s(b)\neq\emptyset\Rightarrow\begin{cases}b\in\cc_{ws}(a)&\textrm{if
$\ell(ws)=\ell(w)+1$, and}\\
b\in\cc_w(a)\cup\cc_{ws}(a)&\textrm{if
$\ell(ws)=\ell(w)-1$.}\end{cases}
\end{align*}
\end{lem}

\begin{proof} Let $s=s_i$ where $i\in I$. Suppose first that
$\ell(ws)=\ell(w)+1$ and that $c\in \cc_w(a)\cap \cc_s(b)$. Let
$f$ be a reduced word in $I^*$ so that $s_f=w$, and so there
exists a gallery from $a$ to $c$ of type $f$. Since
$b\in\cc_{s}(c)$, there is a gallery of type $fi$ from $a$ to
$b$, which is a reduced word by hypothesis. It follows that $b\in
\cc_{ws}(a)$.

Suppose now that $\ell(ws)=\ell(w)-1$, and that $c\in \cc_w(a)\cap
\cc_s(b)$. Since $ws$ is not reduced, there exists a reduced word
$f'$ such that $f'i$ is a reduced word for~$w$. This shows that
there exist a minimal gallery $a=a_0,\ldots,a_{m}=c$ such that
$a_{m-1}\in\cc_s(c)$. Since $b\in\cc_{s}(c)$ too, it follows that
either $b=a_{m-1}$ or $b\in\cc_{s}(a_{m-1})$. In the former case
we have $b\in \cc_{ws}(a)$ and in the latter we have $b\in
\cc_{w}(a)$.  \end{proof}

We now perform counts that will be used to demonstrate chamber
regularity.
\begin{lem}\label{secondcount} Let $w\in W$ and $s\in S$. Fix $a,b\in\cc$. Then
\begin{align*}
|\cc_w(a)\cap \cc_s(b)|=\begin{cases}1&\textrm{if
$\ell(ws)=\ell(w)+1$ and $b\in \cc_{ws}(a)$,}\\
q_s&\textrm{if $\ell(ws)=\ell(w)-1$ and $b\in \cc_{ws}(a)$, and}\\
q_s-1&\textrm{if $\ell(ws)=\ell(w)-1$ and $b\in
\cc_w(a)$.}\end{cases}
\end{align*}
\end{lem}

\begin{proof} Suppose first that $\ell(ws)=\ell(w)+1$ and that
$b\in \cc_{ws}(a)$. Thus there is a minimal gallery
$a=a_0,\ldots,a_{m}=b$ such that $a_{m-1}\in\cc_s(b)$. There are
$q_s$ chambers $c$ in~$\cc_{s}(b)$. One of these chambers is
$a_{m-1}$, which lies in $\cc_{w}(a)$, and the remaining $q_s-1$
lie in $\cc_{ws}(a)$, so $a_{m-1}$ is the only element of
$\cc_{w}(a)\cap\cc_{s}(b)$. Thus $|\cc_w(a)\cap \cc_s(b)|=1$ as
claimed in this case.

Suppose now that $\ell(ws)=\ell(w)-1$ and that $b\in \cc_{ws}(a)$. Write
$s=s_i$, and
let $w=s_f$ where $f\in I^*$ is reduced. Since
$\ell(ws)=\ell(w)-1$, there exists a reduced word $f'$ such that
$f'i$ is a reduced word for $w$, and thus there exists a
minimal gallery of type $f'$ from $a$ to $b$. Thus each
$c\in\cc_s(b)$ can be joined to $a$ by a gallery of type $f'i\sim
f$, and hence $c\in\cc_w(a)$, verifying the count in this case.

Finally, suppose that $\ell(ws)=\ell(w)-1$ and $b\in \cc_w(a)$.
Then, as in the proof of Lemma~\ref{firstcount}, there exists a
minimal gallery $a=a_0,\ldots,a_{m}=b$ such that
$b\in\cc_s(a_{m-1})$. Exactly one of the $q_s$ chambers
$c\in\cc_s(b)$ equals $a_{m-1}$, and thus lies in $\cc_{ws}(a)$.
For the remaining $q_s-1$ chambers we have $c\in\cc_s(a_{m-1})$,
and thus $c\in\cc_w(a)$, completing the proof.  \end{proof}

\begin{thm}\label{recursive} Let $w\in W$ and $s\in S$. Then
\begin{align*}
B_wB_s=\begin{cases}B_{ws}&\textrm{when $\ell(ws)=\ell(w)+1$\,,}\\
\frac{1}{q_s}B_{ws}+\left(1-\frac{1}{q_s}\right)B_w&\textrm{when
$\ell(ws)=\ell(w)-1$\,.}\end{cases}
\end{align*}
\end{thm}

\begin{proof} Let us look at the case $\ell(ws)=\ell(w)-1$. The case $\ell(ws)=\ell(w)+1$ is similar.
By (\ref{one}) and Lemma~\ref{secondcount} we have
$$B_wB_s=\frac{q_{ws}}{q_w}B_{ws}+\left(1-\frac{1}{q_s}\right)B_w\,.$$
All that remains is to show that
$\frac{q_{ws}}{q_w}=\frac{1}{q_s}$. If $f$ is a reduced word with
$s_f=w$ and $s=s_i$, the hypothesis that $\ell(ws)=\ell(w)-1$
implies that there exists a reduced word $f'$ such that $f'i$ is
a reduced word for $w$. The result now follows.  \end{proof}

\begin{cor}\label{brec} $B_{w_1}B_{w_2}=B_{w_1w_2}$ whenever
$\ell(w_1w_2)=\ell(w_1)+\ell(w_2)$.
\end{cor}

\begin{cor}\label{numbers} Let $w_1,w_2\in W$. There exist numbers
$b_{w_1,w_2;w_3}\in\bq^+$ such that
$$B_{w_1}B_{w_2}=\sum_{w_3\in W}b_{w_1,w_2;w_3}B_{w_3}\qquad\textrm{and}\qquad\sum_{w_3\in W}b_{w_1,w_2;w_3}=1\,.$$
Moreover, $|\{w_3\in W\mid b_{w_1,w_2;w_3}\neq0\}|$ is finite for
all $w_1,w_2\in W$.
\end{cor}

\begin{proof} An induction on $\ell(w_2)$ shows existence of
the numbers $b_{w_1,w_2;w_3}\in\bq^+$ such that
$B_{w_1}B_{w_2}=\sum_{w_3}b_{w_1,w_2;w_3}B_{w_3}$, and shows that
only finitely many of the $b_{w_1,w_2;w_3}$'s are nonzero for
fixed $w_1$ and $w_2$. Evaluating both sides at the constant
function $1_{\cc}:\cc\to\{1\}$ shows that
$\sum_{w_3}b_{w_1,w_2;w_3}=1$.  \end{proof}

\begin{defn}\label{balg} Let $\scb$ be the linear span of the set $\{B_w\mid
w\in W\}$ over $\bc$. Corollary~\ref{numbers} shows that $\scb$ is
an associative algebra.
\end{defn}

\begin{prop}\label{vspbase} $\{B_w\mid w\in W\}$ is a vector
space basis of $\scb$, and $\scb$ is generated by $\{B_s\mid s\in
S\}$.
\end{prop}

\begin{proof} Suppose we have a relation $\sum_{k=1}^{n}b_kB_{w_k}=0$, and fix $a,b\in\cc$ with
$\delta(a,b)=w_j$ with $1\leq j\leq n$. Then writing
$\delta_b=1_{\{b\}}$ we have
$$0=\sum_{k=1}^{n}b_k(B_{w_k}\delta_b)(a)=\sum_{k=1}^nb_kq_{w_{k}}^{-1}\delta_{k,j}=b_{j}q_{w_j}^{-1},$$
and so $b_j=0$. From Corollary \ref{brec} we see that $\{B_s\mid
s\in S\}$ generates $\scb$.  \end{proof}

We refer to the numbers $b_{w_1,w_2;w_3}$ from Corollary
\ref{numbers} as the \textit{structure constants} of the algebra
$\scb$ (with respect to the natural basis $\{B_w\mid w\in W\}$).

\begin{prop}\label{cor1} Let $\scx$ be a regular building of type~$W$, and 
let $w_1,w_2,w_3\in W$. For any pair $a,b\in\cc$ with $b\in \cc_{w_3}(a)$ we have
$$|\cc_{w_1^{\vphantom{-1}}}(a)\cap
\cc_{w_2^{-1}}(b)|=\frac{q_{w_1}q_{w_2}}{q_{w_3}}b_{w_1,w_2;w_3}\,,$$
and so $\scx$ is chamber regular.
\end{prop}

\begin{proof} Using (\ref{one}) we compute
$(B_{w_1}B_{w_2}\delta_b)(a)=q_{w_1}^{-1}q_{w_2}^{-1}|\cc_{w_1^{\vphantom{-1}}}(a)\cap
\cc_{w_2^{-1}}(b)|$, whereas by Corollary \ref{numbers} we have
$(B_{w_1}B_{w_2}\delta_b)(a)=q_{w_3}^{-1}b_{w_1,w_2;w_3}$.
\end{proof}

Those readers familiar with Hecke algebras will notice immediately
from Theorem~\ref{recursive} the connection between $\scb$ and
Hecke algebras. For our purposes we define \textit{Hecke algebras} as follows (see \cite[Chapter 7]{h}). For each
$s\in S$, let $a_s$ and $b_s$ be complex numbers such that
$a_{s'}=a_s$ and $b_{s'}=b_s$ whenever $s'=wsw^{-1}$ for some
$w\in W$. The \textit{(generic) Hecke algebra} $\sch(a_s,b_s)$ is
the algebra over $\bc$ with presentation given by basis elements
$T_w$, $w\in W$, and relations
\begin{align}
\label{rr}T_wT_s=\begin{cases}T_{ws}&\textrm{when $\ell(ws)=\ell(w)+1$\,,}\\
a_sT_{ws}+b_sT_w&\textrm{when $\ell(ws)=\ell(w)-1$\,.}\end{cases}
\end{align}

\begin{thm}\label{usebalg} Suppose a building $\scx$ of type $W$ exists with parameters
$\{q_{s}\}_{s\in S}$. Then $\scb\cong\sch(q_s^{-1},1-q_s^{-1})$.
\end{thm}

\begin{proof} We note first that by Corollary~\ref{conjno}, the numbers $a_s=q_s^{-1}$ and
$b_s=1-q_s^{-1}$ satisfy the condition $a_{s'}=a_s$ and
$b_{s'}=b_s$ whenever $s'=wsw^{-1}$ for some $w\in W$.

Since $\{T_w\mid w\in W\}$ is a vector space basis of
$\sch(q_s^{-1},1-q_s^{-1})$ and $\{B_w\mid w\in W\}$ is a vector
space basis of $\scb$ (see Proposition~\ref{vspbase}) there exists
a unique vector space isomorphism
$\Phi:\sch(q_s^{-1},1-q_s^{-1})\to\scb$ such that $\Phi(T_w)=B_w$
for all $w\in W$. By (\ref{rr}) and Theorem~\ref{recursive} we
have $\Phi(T_wT_s)=\Phi(T_w)\Phi(T_s)$ for all $w\in W$ and $s\in
S$, and so $\Phi$ is an algebra homomorphism. It follows that
$\Phi$ is an algebra isomorphism.  \end{proof}

Recall that we write $D$ for the Coxeter graph of $W$.

\begin{defn}\label{autoq} Let $\scx$ be a locally finite regular building. Define
$$\Aut_{q}(D)=\{\s\in \Aut(D)\mid q_{\s(i)}=q_i\textrm{ for all
}i\in I\}.$$
\end{defn}

\begin{lem}\label{vertexuse} For all $w_1,w_2\in W$
and $\s\in\Aut_{q}(D)$ we have
$$|\cc_{\s(w_1)}(a')\cap\cc_{\s(w_2)}(b')|=|\cc_{w_1}(a)\cap\cc_{w_2}(b)|\,,$$
whenever $a,b,a',b'\in\cc$ are chambers with
$\delta(a',b')=\s(\delta(a,b))$.
\end{lem}

\begin{proof} We first show that, in the notation of Corollary
\ref{numbers},
\begin{align}
\label{eq}b_{w_1,w_2;w_3}=b_{\s(w_1),\s(w_2);\s(w_3)}
\end{align}
for all $w_1,w_2,w_3\in W$.

Theorem~\ref{recursive}, the definition of $\Aut_q(D)$ and the
fact that $\ell(\s(w))=\ell(w)$ for all $w\in W$ show that this
is true when $\ell(w_2)=1$, beginning an induction. Suppose
(\ref{eq}) holds whenever $\ell(w_2)<n$, and suppose
$w=s_{i_1}\cdots s_{i_{n-1}}s_{i_n}$ has length $n$. Write
$w'=s_{i_1}\cdots s_{i_{n-1}}$ and $s=s_{i_n}$. Observe that
$\s(w)=\s(w')\s(s)$ so that $B_{\s(w)}=B_{\s(w')}B_{\s(s)}$ by
Theorem~\ref{recursive}, and so
\begin{align*}
B_{\s(w_1)}B_{\s(w)}&=(B_{\s(w_1)}B_{\s(w')})B_{\s(s)}\\
&=\sum_{w_3\in W}b_{\s(w_1),\s(w');\s(w_3)}B_{\s(w_3)}B_{\s(s)}\\
&=\sum_{w_3\in W}\left(b_{w_1,w';w_3}\sum_{w_4\in W}b_{\s(w_3),\s
(s);\s(w_4)}B_{\s(w_4)}\right)\\
&=\sum_{w_4\in W}\left(\sum_{w_3\in
W}b_{w_1,w';w_3}b_{w_3,s;w_4}\right)B_{\s(w_4)}\,.
\end{align*}
Thus
\begin{align}
\label{fo}b_{\s(w_1),\s(w);\s(w_4)}=\sum_{w_3\in
W}b_{w_1,w';w_3}b_{w_3,s;w_4}\qquad\textrm{for all $w_4\in W$}.
\end{align}
The same calculation without the $\s$'s shows that this is also
$b_{w_1,w;w_4}$. This completes the induction step, and so
(\ref{eq}) holds for all $w_1,w_2$ and $w_3$ in $W$.

Thus for any chambers $a,b,a',b'$ with $\delta(a,b)=w_3$, and
$\delta(a',b')=\s(w_3)$ we have (using Proposition~\ref{cor1})
\begin{align*}
|\cc_{w_1}(a)\cap\cc_{w_2}(b)|&=\frac{q_{w_1^{\vphantom{-1}}}q_{w_2^{-1}}}{q_{w_3}}b_{w_1^{\vphantom{-1}},w_2^{-1};w_3^{\vphantom{-1}}}\\
&=\frac{q_{\s(w_1^{\vphantom{-1}})}q_{\s(w_2^{-1})}}{q_{\s(w_3)}}b_{\s(
w_1^{\vphantom{-1}}),\s(
w_2^{-1});\s(w_3^{\vphantom{-1}})}\\
&=|\cc_{\s(w_1)}(a')\cap\cc_{\s(w_2)}(b')|.\qedhere
\end{align*}
\end{proof}

\section{Preliminary Material}\label{prelim}

This section is preparation for our study of the vertex set
averaging operators associated to locally finite regular affine
buildings.

\subsection{Chamber Systems and Simplicial Complexes}
Recall that a \textit{simplicial complex} with vertex set $V$ is
a collection $X$ of finite subsets of $V$ (called
\textit{simplices}) such that for every $v\in V$, the singleton
$\{v\}$ is a simplex, and every subset of a simplex $x$ is a
simplex (called a \textit{face of $x$}). If $x$ is a simplex
which is not a proper subset of any other simplex, then we call
$x$ a \textit{maximal simplex}, or \textit{chamber} of $X$.

A \textit{labelled simplicial complex} with vertex set $V$ is a
simplicial complex equipped with a set $I$ of \textit{types}, and
a \textit{type map} $\tau:V\to I$ such that the restriction
$\tau|_{C}:C\to I$ of $\tau$ to any chamber $C$ is a bijection.

An \textit{isomorphism} of simplicial complexes is a bijection of
the vertex sets that maps simplices, and only simplices, to simplices. If both
simplicial complexes are labelled by the same set, then an
isomorphism which preserves types is said to be \textit{type
preserving}.

There is a well known method of producing labelled simplicial
complexes from chamber systems, and vice versa (see
\cite[\ts1.4]{donaldintro} for details). This allows us to
consider buildings and Coxeter complexes as certain labelled
simplicial complexes (with canonical labellings). The following
is an alternative (and of course equivalent) definition of
buildings from a simplicial complex approach.

\begin{defn}\label{def2}\cite{brown}\textbf{.} Let $W$ be a
Coxeter group of type $M$. A \textit{building of type $M$} is a
nonempty simplicial complex $\scx$ which contains a family of
subcomplexes called \textit{apartments} such that
\begin{enumerate}
\item[$\mathrm{(i)}$] each apartment is isomorphic to the (simplicial) Coxeter
complex of $W$,
\item[$\mathrm{(ii)}$] given any two chambers of $\scx$ there is an
apartment containing both, and
\item[$\mathrm{(iii)}$] given any
two apartments $\ca$ and $\ca'$ that contain a common chamber,
there exists an isomorphism $\psi:\ca\to\ca'$ fixing $\ca\cap\ca'$
pointwise.
\end{enumerate}
\end{defn}

We remark that Definition~\ref{def2}(iii) can be replaced with
the following \cite[p.76]{brown}.
\begin{itemize}
\item[(iii)$'$] If $\ca$ and $\ca'$ are apartments both containing
simplices $\rho$ and $\sigma$, then there is an isomorphism
$\psi:\ca\to\ca'$ fixing $\rho$ and $\s$ pointwise.
\end{itemize}

It is easy to see that $\scx$ is in fact a labellable simplicial
complex, and all the isomorphisms in the above definition may be
taken to be label preserving.

\subsection{Root Systems} For the purpose of fixing notation we will give a brief discussion
of root systems. A thorough reference to this well known material
is \cite{bourbaki}.

Let $E$ be an $n$-dimensional vector space over $\br$ with inner
product $\lan\cdot,\cdot\ran$, and for $\alpha\in
E\backslash\{0\}$ define
$\alpha^{\vee}=\frac{2\alpha}{\lan\alpha,\alpha\ran}$. Let $R$ be
an \textit{irreducible}, but not necessarily \textit{reduced},
\textit{root system} in $E$ (see \cite[VI, \ts1,
No.1-2]{bourbaki}).

The elements of $R$ are called \textit{roots}, and the
\textit{rank} of $R$ is $n$, the dimension of~$E$. A root system
that is not reduced is said to be \textit{non-reduced}. See
\cite[VI, \ts4, No.5--No.14]{bourbaki} for the classification of
irreducible root systems.

Let $B=\{\alpha_i\mid i\in I_0\}$ be a \textit{base} of $R$,
where $I_0=\{1,2,\ldots,n\}$. Thus $B$ is a subset of $R$ such
(i) a vector space basis of $E$, and (ii) each root in $R$ can be
written as a linear combination of elements of $B$ with integer
coefficients which are either all nonnegative or all nonpositive.
We say that $\alpha\in R$ is \textit{positive} (respectively
\textit{negative}) if the expression for $\alpha$ from (ii) has
only nonnegative (respectively nonpositive) coefficients. Let
$R^+$ (respectively $R^-$) be the set of all positive
(respectively negative) roots. Thus $R^-=-R^+$ and $R=R^{+}\cup
R^-$, where the union is disjoint.

Define the \textit{height (with respect to $B$)} of
$\alpha=\sum_{i\in I_0}k_i\alpha_i\in R$ by
$\mathrm{ht}(\alpha)=\sum_{i\in I_0 }k_{i}$. By \cite[VI, \ts 1
No.8, Proposition~25]{bourbaki} there exists a unique root
$\tilde{\alpha}\in R$ whose height is maximal, and defining
numbers $m_i$ by
\begin{align}
\label{highestroot}\tilde{\alpha}=\sum_{i\in I_0}m_i\alpha_i
\end{align}
we have $m_i\geq1$ for all $i\in I_0$. To complete the notation
we define $m_0=1$.

The \textit{dual} (or \textit{inverse}) of $R$ is
$\chR=\{\cha\mid\alpha\in R\}$. By \cite[VI, \ts1, No.1,
Proposition~2]{bourbaki} $\chR$ is an irreducible root system
which is reduced if and only if $R$ is.

We define a dual basis $\{\la_i\}_{i\in I_0}$ of $E$ by
$\lan\la_i,\alpha_j\ran=\delta_{i,j}$. Recall that the
\textit{coroot lattice} $Q$ of $R$ is the $\bz$-span of
$R$, and the \textit{coweight lattice} $P$ of $R$
is the $\bz$-span of $\{\la_i\}_{i\in I_0}$. Elements of $P$ are
called \textit{coweights} (of $R$), and it is clear that
$Q\subseteq P$. Note that in the literature $Q$ and $P$ are also
called the \textit{root} and \textit{weight} lattices of $R^{\vee}$.
We call a coweight $\la=\sum_{i\in I_0}a_i\la_i$ \textit{dominant}
if $a_i\geq0$ for all $i\in I_0$, and we write $P^+$ for the set
of all dominant coweights.

For each $n\geq1$ there is exactly one irreducible non-reduced
root system (up to isomorphism) of rank $n$, denoted by $BC_n$
\cite[VI, \ts4, No.14]{bourbaki}. We may take $E=\br^n$ with the
usual inner product, and let $\alpha_j=e_j-e_{j+1}$ for $1\leq
j<n$ and $\alpha_n=e_n$. Then $B=\{\alpha_j\}_{j=1}^{n}$, and
$$R^+=\{e_k,2e_k,e_i+e_j,e_i-e_j\mid 1\leq k\leq n,\,1\leq i<j\leq
n\}\,.$$ Notice that $\chR=R$, and one easily sees that $Q=P$.

\subsection{Hyperplane Arrangements and Reflection Groups}\label{43} Let $R$ be an irreducible root system, and for each $\alpha\in R$ and $k\in\bz$
let $H_{\alpha;k}=\{x\in E\mid\lan x,\alpha\ran=k\}$. Let $\ch$
denote the family of these (affine) \textit{hyperplanes}
$H_{\alpha;k}$, $\alpha\in R$, $k\in\bz$. We write $H_{\alpha}$
in place of $H_{\alpha;0}$, and denote by $\ch_0$ the family of
these hyperplanes $H_{\alpha}$, $\alpha\in R$.

Given $H_{\alpha;k}\in \ch$, the associated \textit{orthogonal
reflection} is the map $s_{\alpha;k}:E\to E$ given by
$s_{\alpha;k}(x)=x-(\lan x,\alpha\ran-k)\cha$ for all $x\in E$.
We write $s_{\alpha}$ in place of $s_{\alpha;0}$, and $s_i$ in
place of $s_{\alpha_i}$. The \textit{Weyl group of $R$}, denoted
$W_0(R)$, or simply $W_0$, is the subgroup of $\mathrm{GL}(E)$
generated by the reflections $s_{\alpha}$, $\alpha\in R$, and the
\textit{affine Weyl group of $R$}, denoted $W(R)$, or simply $W$,
is the subgroup of $\mathrm{Aff}(E)$ generated by the reflections
$s_{\alpha;k}$, $\alpha\in R$, $k\in\bz$. Here $\mathrm{Aff}(E)$
is the set of maps $x\mapsto Tx+v$, $T\in \mathrm{GL}(E)$, $v\in
E$. Writing $t_{v}$ for the translation $x\mapsto x+v$, we
consider $E$ as a subgroup of $\mathrm{Aff}(E)$ by identifying
$v$ and~$t_v$. We have $\mathrm{Aff}(E)=\mathrm{GL}(E)\ltimes E$, and $W\cong
W_0\ltimes Q$. Note that $W_0(R^{\vee})=W_0(R)$ \cite[VI,\ts1,No.1]{bourbaki}.

Let $s_{0}=s_{\tilde{\alpha};1}$, define $I=I_0\cup\{0\}$, and
let $S_0=\{s_i\mid i\in I_0\}$ and $S=\{s_i\mid i\in I\}$. The
group $W_0$ (respectively $W$) is a Coxeter group over $I_0$
(respectively $I$) generated by $S_0$ (respectively $S$).

We write $\S=\S(R)$ for the vector space $E$ equipped with the
sectors, chambers and vertices as defined below. The open
connected components of $E\backslash\bigcup_{H\in\ch}H$ are
called the \textit{chambers} of $\S$ (this terminology is
motivated by building theory, and differs from that used in
\cite{bourbaki} where there are \textit{chambers} and
\textit{alcoves}), and we write $\cc(\S)$ for the set of chambers
of $\S$. Since $R$ is irreducible, each $C\in\cc(\S)$ is an open
(geometric) simplex \cite[V, \ts3, No.9,
Proposition~8]{bourbaki}. Call the extreme points of the sets
$\overline{C}$, $C\in\cc(\S)$, \textit{vertices} of $\S$, and
write $V(\S)$ for the set of all vertices of $\S$.

In choice of $B$ gives a natural \textit{fundamental chamber}
\begin{align}
\label{fundamentalalcove}C_0=\{x\in E\mid\lan
x,\alpha_i\ran>0\textrm{ for all $i\in I_0$ and }\lan
x,\tilde{\alpha}\ran<1\},
\end{align}
where we use the notation of (\ref{highestroot}).

The \textit{fundamental sector} of $\S$ is
\begin{align}
\label{fundamentalsector}\cs_0=\{x\in E\mid\lan
x,\alpha_i\ran>0\textrm{ for all $i\in I_0$}\},
\end{align}
and the \textit{sectors} of $\S$ are the sets $\la+w\cs_0$, where
$\la\in P$ and $w\in W_0$. The sector $\cs=\la+w\cs_0$ is said to
have \textit{base vertex} $\la$ (we will see in
Section~\ref{laterefad} that $\la$ is indeed a vertex of $\S$).

The group $W_0$ acts simply transitively on the set of sectors
based at $0$, and $\overline{\cs}_0$ is a fundamental domain for the action of
$W_0$ on $E$.
Similarly, $W$ acts simply transitively on $\cc(\S)$, and $\overline{C}_0$ if a
fundamental domain for the action of $W$ on $E$ \cite[VI, \ts1-3]{bourbaki}.

The following fact follows easily from \cite[VI, \ts2, No.2,
Proposition~4(ii)]{bourbaki}.
\begin{lem}\label{rt} $W_0$ acts simply transitively on the set of $C\in\cc(\S)$ with $0\in\overline{C}$.
\end{lem}

\subsection{A Geometric Realisation of the Coxeter Complex}\label{45} The set $\cc(\S)$
from Section~\ref{43} forms a chamber system over $I$ if we
declare $wC_0\sim_i wC_0$ and $wC_0\sim_i ws_iC_0$ for each $w\in
W$ and each $i\in I$. The map $w\mapsto wC_0$ is an isomorphism
of the Coxeter complex $\cc(W)$ of Section~\ref{s1} onto this
chamber system, and so $\S$ may be regarded as a
\textit{geometric realisation} of $\cc(W)$.

The vertices of $\overline{C}_0$ are $\{0\}\cup\{\la_i/m_i\mid
i\in I_0\}$ (see \cite[VI, \ts2, No.2]{bourbaki}), and we declare
$\tau(0)=0$ and $\tau(\la_i/m_i)=i$ for $i\in I_0$. This extends
to a unique labelling $\tau:V(\S)\to I$ (see
\cite[Lemma~1.5]{donaldintro}), and the action of $W$ on $\S$ is
type preserving.\begin{comment}
 \cite[Theorem, p.58]{brown}.
\end{comment}

\subsection{Special and Good Vertices of $\S$}\label{laterefad} Following \cite[V, \ts3, No.10]{bourbaki},
a point $v\in E$ is said to be \textit{special} if for every
$H\in\ch$ there exists a hyperplane $H'\in\ch$ parallel to $H$
such that $v\in H'$. Note that in our set-up $0\in E$ is special.
Each special point is a vertex of $\S$ \cite[V, \ts3,
No.10]{bourbaki}, and thus we will call the special points
\textit{special vertices}. Note that in general not all vertices
are special (for example, in the $\widetilde{C}_2$ and
$\widetilde{G}_2$ complexes). When $R$ is reduced $P$ is the set
of special vertices of~$\S$ \cite[VI, \ts2, No.2,
Proposition~3]{bourbaki}. When $R$ is non-reduced then $P$ is a
proper subset of the special vertices of~$\S$ (see
Example~\ref{example4}).

To deal with the reduced and non-reduced cases simultaneously, we
define the \textit{good} vertices of $\S$ to be the elements of~$P$. On the first reading the reader is encouraged to think of~$P$
as the set of all special vertices, for this is true unless~$R$ is
of type $BC_n$. Note that, according to our definitions, every
sector of $\S$ is based at a good vertex of $\S$.

We write $I_{P}$ for the set of \textit{good types}. That is,
$I_P=\{\tau(\la)\mid \la\in P\}\subseteq I$. 

\begin{comment}The following gives a
simple method for computing $I_P$, and shows that $I_P=\{0\}$ if
$R$ is non-reduced.
\end{comment}
\begin{lem} Let the numbers $m_i$ be as in (\ref{highestroot}). Then $I_P=\{i\in I\mid
m_i=1\}$.
\end{lem}

\begin{proof} The vertices of $C_0$ are
$\{0\}\cup\{\la_i/m_i\mid i\in I_0\}$. The good vertices of $C_0$
are those in $P$, and thus have type $0$ or $i$ for some $i$ with
$m_i=1$.  \end{proof}

\subsection{Examples}\label{examples} \begin{comment}Let us briefly look at
two examples to highlight the difference between the reduced and
non-reduced cases.
\end{comment}

\begin{exmp}[$R=C_2$]\label{example3} Take $E=\br^2$,
$\alpha_1=e_1-e_2$ and $\alpha_2=2e_2$. Then
$B=\{\alpha_1,\alpha_2\}$ and
$R^+=\{\alpha_1,\alpha_2,\alpha_1+\alpha_2,2\alpha_1+\alpha_2\}$.

\begin{figure}[ht]
\begin{center}
\psset{xunit= 0.7 cm,yunit= 0.7 cm} \psset{origin={0,0}}

\vspace{3cm}

    \psline[linestyle=dashed,dash=3pt 2pt](0,-4.5)(0,4.5)
    \psline[linestyle=dashed,dash=3pt 2pt](-1,-4.5)(-1,4.5)
    \psline[linestyle=dashed,dash=3pt 2pt](-2,-4.5)(-2,4.5)
    \psline[linestyle=dashed,dash=3pt 2pt](-3,-4.5)(-3,4.5)
    \psline[linestyle=dashed,dash=3pt 2pt](-4,-4.5)(-4,4.5)
    \psline[linestyle=dashed,dash=3pt 2pt](1,-4.5)(1,4.5)
    \psline[linestyle=dashed,dash=3pt 2pt](2,-4.5)(2,4.5)
    \psline[linestyle=dashed,dash=3pt 2pt](3,-4.5)(3,4.5)
    \psline[linestyle=dashed,dash=3pt 2pt](4,-4.5)(4,4.5)
    \psline[linestyle=dashed,dash=3pt 2pt](-4.5,4)(4.5,4)
    \psline[linestyle=dashed,dash=3pt 2pt](-4.5,3)(4.5,3)
    \psline[linestyle=dashed,dash=3pt 2pt](-4.5,2)(4.5,2)
    \psline[linestyle=dashed,dash=3pt 2pt](-4.5,1)(4.5,1)
    \psline[linestyle=dashed,dash=3pt 2pt](-4.5,0)(4.5,0)
    \psline[linestyle=dashed,dash=3pt 2pt](-4.5,-1)(4.5,-1)
    \psline[linestyle=dashed,dash=3pt 2pt](-4.5,-2)(4.5,-2)
    \psline[linestyle=dashed,dash=3pt 2pt](-4.5,-3)(4.5,-3)
    \psline[linestyle=dashed,dash=3pt 2pt](-4.5,-4)(4.5,-4)
    \psline[linestyle=dotted,dotsep=2pt](-4.5,3.5)(-3.5,4.5)
    \psline[linestyle=dotted,dotsep=2pt](-4.5,1.5)(-1.5,4.5)
    \psline[linestyle=dotted,dotsep=2pt](-4.5,-0.5)(0.5,4.5)
    \psline[linestyle=dotted,dotsep=2pt](-4.5,-2.5)(2.5,4.5)
    \psline[linestyle=dotted,dotsep=2pt](-4.5,-4.5)(4.5,4.5)
    \psline[linestyle=dotted,dotsep=2pt](3.5,-4.5)(4.5,-3.5)
    \psline[linestyle=dotted,dotsep=2pt](1.5,-4.5)(4.5,-1.5)
    \psline[linestyle=dotted,dotsep=2pt](-0.5,-4.5)(4.5,0.5)
    \psline[linestyle=dotted,dotsep=2pt](-2.5,-4.5)(4.5,2.5)
    \psline[linestyle=dotted,dotsep=2pt](4.5,3.5)(3.5,4.5)
    \psline[linestyle=dotted,dotsep=2pt](4.5,1.5)(1.5,4.5)
    \psline[linestyle=dotted,dotsep=2pt](4.5,-0.5)(-0.5,4.5)
    \psline[linestyle=dotted,dotsep=2pt](4.5,-2.5)(-2.5,4.5)
    \psline[linestyle=dotted,dotsep=2pt](4.5,-4.5)(-4.5,4.5)
    \psline[linestyle=dotted,dotsep=2pt](-3.5,-4.5)(-4.5,-3.5)
    \psline[linestyle=dotted,dotsep=2pt](-1.5,-4.5)(-4.5,-1.5)
    \psline[linestyle=dotted,dotsep=2pt](0.5,-4.5)(-4.5,0.5)
    \psline[linestyle=dotted,dotsep=2pt](2.5,-4.5)(-4.5,2.5)
    \psline[linewidth=1pt]{<->}(-4,0)(4,0)
    \psline[linewidth=1pt]{<->}(0,4)(0,-4)
    \psline[linewidth=1pt]{<->}(-2,-2)(2,2)
    \psline[linewidth=1pt]{<->}(-2,2)(2,-2)
    \rput(0.70,0.25){$\scriptstyle C_0$}
    \rput(1.7,-1.4){$\scriptstyle \alpha_1$}
    \rput(0.3,3.3){$\scriptstyle \alpha_2$}

\vspace{3cm}

\end{center}

\caption{}\label{fig:1}

\end{figure}

The dotted lines in Figure~\ref{fig:1} are the hyperplanes
$\{H_{w\alpha_1;k}\mid w\in W_0,k\in\bz\}$, and the dashed lines
are the hyperplanes $\{H_{w\alpha_2;k}\mid w\in W_0,k\in\bz\}$.
In this example $\la_1=e_1$ and $\la_2=\frac{1}{2}(e_1+e_2)$, and
$\tau(0)=0$, $\tau(\frac{1}{2}e_1)=1$ and
$\tau(\frac{1}{2}(e_1+e_2))=2$. We have
$P=\{(x,y)\in(\frac{1}{2}\bz)^2\mid x+y\in\bz\}$, which coincides
with the set of all special vertices (as expected, since $R$ is
reduced here). Thus $I_P=\{0,2\}$.
\end{exmp}

\begin{exmp}[$R=BC_2$]\label{example4} Take $E=\br^2$,
$\alpha_1=e_1-e_2$ and $\alpha_2=e_2$. Then
$B=\{\alpha_1,\alpha_2\}$ and
$R^+=\{\alpha_1,\alpha_2,\alpha_1+\alpha_2,\alpha_1+2\alpha_2,2\alpha_2,2\alpha_1+2\alpha_2\}$.

\begin{figure}[ht]
\begin{center}
\psset{xunit= 0.7 cm,yunit= 0.7 cm} \psset{origin={0,0}}

\vspace{3cm}

    \psline(0,-4.5)(0,4.5)
    \psline[linestyle=dashed,dash=3pt 2pt](-1,-4.5)(-1,4.5)
    \psline(-2,-4.5)(-2,4.5)
    \psline[linestyle=dashed,dash=3pt 2pt](-3,-4.5)(-3,4.5)
    \psline(-4,-4.5)(-4,4.5)
    \psline[linestyle=dashed,dash=3pt 2pt](1,-4.5)(1,4.5)
    \psline(2,-4.5)(2,4.5)
    \psline[linestyle=dashed,dash=3pt 2pt](3,-4.5)(3,4.5)
    \psline(4,-4.5)(4,4.5)
    \psline(-4.5,4)(4.5,4)
    \psline[linestyle=dashed,dash=3pt 2pt](-4.5,3)(4.5,3)
    \psline(-4.5,2)(4.5,2)
    \psline[linestyle=dashed,dash=3pt 2pt](-4.5,1)(4.5,1)
    \psline(-4.5,0)(4.5,0)
    \psline[linestyle=dashed,dash=3pt 2pt](-4.5,-1)(4.5,-1)
    \psline(-4.5,-2)(4.5,-2)
    \psline[linestyle=dashed,dash=3pt 2pt](-4.5,-3)(4.5,-3)
    \psline(-4.5,-4)(4.5,-4)
    \psline[linestyle=dotted,dotsep=2pt](-4.5,3.5)(-3.5,4.5)
    \psline[linestyle=dotted,dotsep=2pt](-4.5,1.5)(-1.5,4.5)
    \psline[linestyle=dotted,dotsep=2pt](-4.5,-0.5)(0.5,4.5)
    \psline[linestyle=dotted,dotsep=2pt](-4.5,-2.5)(2.5,4.5)
    \psline[linestyle=dotted,dotsep=2pt](-4.5,-4.5)(4.5,4.5)
    \psline[linestyle=dotted,dotsep=2pt](3.5,-4.5)(4.5,-3.5)
    \psline[linestyle=dotted,dotsep=2pt](1.5,-4.5)(4.5,-1.5)
    \psline[linestyle=dotted,dotsep=2pt](-0.5,-4.5)(4.5,0.5)
    \psline[linestyle=dotted,dotsep=2pt](-2.5,-4.5)(4.5,2.5)
    \psline[linestyle=dotted,dotsep=2pt](4.5,3.5)(3.5,4.5)
    \psline[linestyle=dotted,dotsep=2pt](4.5,1.5)(1.5,4.5)
    \psline[linestyle=dotted,dotsep=2pt](4.5,-0.5)(-0.5,4.5)
    \psline[linestyle=dotted,dotsep=2pt](4.5,-2.5)(-2.5,4.5)
    \psline[linestyle=dotted,dotsep=2pt](4.5,-4.5)(-4.5,4.5)
    \psline[linestyle=dotted,dotsep=2pt](-3.5,-4.5)(-4.5,-3.5)
    \psline[linestyle=dotted,dotsep=2pt](-1.5,-4.5)(-4.5,-1.5)
    \psline[linestyle=dotted,dotsep=2pt](0.5,-4.5)(-4.5,0.5)
    \psline[linestyle=dotted,dotsep=2pt](2.5,-4.5)(-4.5,2.5)
    \psline[linewidth=1pt]{<->}(-4,0)(4,0)
    \psline[linewidth=1pt]{<->}(0,4)(0,-4)
    \psline[linewidth=1pt]{<->}(-2,-2)(2,2)
    \psline[linewidth=1pt]{<->}(-2,2)(2,-2)
    \rput(0.70,0.25){$\scriptstyle C_0$}
    \rput(1.7,-1.4){$\scriptstyle \alpha_1$}
    \rput(0.35,3.2){$\scriptstyle 2\alpha_2$}
    \rput(0.3,1.3){$\scriptstyle \alpha_2$}
    \psline[linewidth=1pt]{<->}(-2,0)(2,0)
    \psline[linewidth=1pt]{<->}(0,2)(0,-2)

\vspace{3cm}

\end{center}
\caption{}\label{fig:2}
\end{figure}

The dotted and solid lines in Figure~\ref{fig:2} represent the
hyperplanes in the sets $\{H_{w\alpha_1;k}\mid w\in
W_0,k\in\bz\}$ and $\{H_{w\alpha_2;k}\mid w\in W_0,k\in\bz\}$
respectively. The union of the dashed and solid lines represent
the hyperplanes in $\{H_{w(2\alpha_2);k}\mid w\in W_0,k\in\bz\}$.

In contrast to the previous example, here we have $\la_1=e_1$ and
$\la_2=e_1+e_2$. The set of special vertices and the vertex types
are as in Example~\ref{example3}, but here~$P=\bz^2$ (and so
$I_P=\{0\}$).
\end{exmp}

\subsection{The Extended Affine Weyl Group}\label{tilde{W}} The \textit{extended affine
Weyl group of~$R$}, denoted $\tilde{W}(R)$ or simply $\tilde{W}$,
is $\tilde{W}=W_0\ltimes P$. In general $\tilde{W}$ is larger
than~$W$. In fact, $\tilde{W}/W\cong P/Q$ \cite[VI, \ts2,
No.3]{bourbaki}. We note that while $W(C_n)= W(BC_n)$,
$\tilde{W}(C_n)$ is not isomorphic to $\tilde{W}(BC_n)$.

In particular, notice that for each $\la\in P$, the translation $t_{\la}:E\to
E$, $t_{\la}(x)=x+\la$, is in $\tilde{W}$.

The group $\tilde{W}$ permutes the chambers of $\S$, but in
general does not act simply transitively. Recall \cite[\ts
2.2]{m} that for $w\in \tilde{W}$, the \textit{length} of $w$ is
defined by
\begin{align*}
\ell(w)=|\{H\in \ch\mid H\textrm{ separates $C_0$ and
$w^{-1}C_0$}\}|\,.
\end{align*}
When $w\in W$, this definition agrees with the definition of
$\ell(w)$ given previously for Coxeter groups.

The subgroup $G=\{g\in\tilde{W}\mid\ell(g)=0\}$ will play an
important role; it is the stabiliser of $C_0$ in~$\tilde{W}$. We have
$\tilde{W}\cong W\rtimes G$ \cite[VI,
\ts2, No.3]{bourbaki}, and furthermore, $G\cong P/Q$, and so $G$
is a finite abelian group. Let $w_0$ and $w_{0\la}$ denote the
longest elements of $W_0$ and $W_{0\la}$ respectively, where
for $\la\in P$,
\begin{align}
\label{earlier?}W_{0\la}=\{w\in W_0\mid w\la=\la\}.
\end{align}
Recall the definition of the numbers $m_i$ (with $m_0=1$) from
(\ref{highestroot}). Then
\begin{align}
\label{G}G=\{g_i\mid m_i=1\}
\end{align}
where $g_0=1$ and $g_i=t_{\la_i}w_{0\la_i}w_0$ for $i\in
I_P\backslash\{0\}$ (see \cite[VI, \ts2, No.3]{bourbaki} in the
reduced case and note that $G=\{1\}$ in the non-reduced case since
$G\cong P/Q$).

\subsection{Automorphisms of $\S$ and $D$}\label{cp} An
\textit{automorphism} of $\S$ is a bijection $\psi$ of~$E$ that
maps chambers, and only chambers, to chambers with the property that $C\sim_i D$
if and only if $\psi(C)\sim_{i'}\psi(D)$ for some $i'\in I$
(depending on $C,D$ and $i$). Let $\Aut(\S)$ denote the
automorphism group of~$\S$. Clearly $W_0$, $W$ and $\tilde{W}$
can be considered as subgroups of $\Aut(\S)$, and we have $W_0\leq
W\leq\tilde{W}\leq\Aut(\S)$. Note that in some cases $\tilde{W}$
is a proper subgroup of $\Aut(\S)$. For example, if $R$ is of type
$A_2$, then the map $a_1\la_1+a_2\la_2\mapsto a_1\la_2+a_2\la_1$
is in $\Aut(\S)$ but is not in $\tilde{W}$.

Write $D$ for the Coxeter graph of $W$ (see Section~\ref{s1}).
Recall the definition of the type map $\tau:V(\S)\to I$ from
Section~\ref{45}.

\begin{prop}\label{induction} Let $\psi\in\Aut(\S)$. Then there exists $\s\in\Aut(D)$ such that
$(\tau\circ\psi)(v)=(\s\circ\tau)(v)$ for all $v\in V(\S)$. If
$C\sim_i D$, then $\psi(C)\sim_{\s(i)}\psi(D)$.
\end{prop}

\begin{proof} The result follows from \cite[p.64--65]{brown}.  \end{proof}

For each $g_i\in G$ (see (\ref{G})), let $\s_i\in\Aut(D)$ be the
automorphism induced as in Proposition~\ref{induction}. We call
the automorphisms $\s_i\in\Aut(D)$ \textit{type rotating} (for in
the $\widetilde{A}_n$ case they are the permutations $k\mapsto
k+i\mod n+1$), and we write $\Auttr(D)$ for the group of all type
rotating automorphisms of $D$. Thus
\begin{align}
\label{auttrd}\Auttr(D)=\{\s_i\mid i\in I_P\}\,.
\end{align}
Note that since $g_0=1$, $\s_0=\mathrm{id}$.

Let $D_0$ be the Coxeter graph of $W_0$. We have \cite[VI, \ts4,
No.3]{bourbaki}
\begin{align}
\label{gammanormal}\Aut(D)=\Aut(D_0)\ltimes\Auttr(D)\,.
\end{align}

The group $\tilde{W}$ has a presentation with generators $s_i$, $i\in I$,
and $g_j$, $j\in I_P$, and relations (see \cite[(1.20)]{ram})
\begin{equation}
    \label{pres}
    \begin{aligned}
     (s_is_j)^{m_{i,j}}&=1&&\textrm{for all $i,j\in I$, and}\\
     g_j^{\vphantom{-1}}s_ig_j^{-1}&=s_{\s_j(i)}&&\textrm{for all
$i\in I$ and $j\in I_P$\,.}
    \end{aligned}
\end{equation}

\begin{prop}\label{lete} Let $i\in I_P$ and $\s\in\Auttr(D)$.
\begin{enumerate}
\item[(i)] $\s_i(0)=i$.
\item[(ii)] If $\s(i)=i$, then $\s=\s_0=\mathrm{id}$.
\item[(iii)] $\Auttr(D)$ acts simply transitively on the good
types of $D$.
\end{enumerate}
\end{prop}

\begin{proof} (i) follows from the formula $g_i=t_{\la_i}w_{0\la_i}w_0$ ($i\in I_0$) given in Section~\ref{tilde{W}}.
By~(i) we have $(\s_i^{-1}\circ\s\circ\s_i)(0)=0$, and so
$\s_i^{-1}\circ\s\circ\s_i=\s_0=\mathrm{id}$. Thus (ii) holds,
and (iii) is now clear.
\end{proof}

\begin{prop}\label{prop:addlate} Let $\psi\in\Aut(\S)$.
\begin{enumerate}
\item[(i)] The image under $\psi$ of a gallery in $\S$ is again a gallery in~$\S$.
\item[(ii)] A gallery in $\S$ is minimal if and only if its image
under $\psi$ is minimal.
\item[(iii)] There exists a unique $\s\in\Aut(D)$ so that $\psi$
maps galleries of type $f$ to galleries of type $\s(f)$. If
$\psi=w\in \tilde{W}$ then $\s\in\Auttr(D)$. If $w=w'g_i$, where
$w'\in W$, then $\s=\s_i$.
\item[(iv)] If $\psi\in\tilde{W}$ maps $\la\in P$ to $\mu\in P$,
then the induced automorphism from (iii) is
$\s=\s_m\circ\s_l^{-1}$, where $l=\tau(\la)$ and $m=\tau(\mu)$.
\end{enumerate}
\end{prop}

\begin{proof} (i) and (ii) are obvious.

(iii) The first statement follows easily from
Proposition~\ref{induction}, and the remaining statements follow
from the definition of $\Auttr(D)$.

(iv) Since $\s(l)=m$, we have $(\s\circ\s_l)(0)=\s_m(0)$, and so
$\s=\s_m\circ\s_l^{-1}$ by Proposition~\ref{lete}.  \end{proof}

\begin{prop}\label{negmap} $x\mapsto-x$ is an automorphism of $\S$.
\end{prop}
\begin{proof} The map $x\mapsto-x$ maps $\ch$ to itself and is
continuous, and so maps chambers to chambers. If $C\sim_i D$ and
$C\neq D$ then there is only one $H\in\ch$ separating $C$ and
$D$, and then $-H$ is the only hyperplane in $\ch$ separating
$-C$ and $-D$, and so $-C\sim_{i'}-D$ for some $i'\in I$.
\end{proof}

\begin{defn}\label{defn:stardefn} Let $\s_*\in\Aut(D)$ be the automorphism of $D$ induced by
the automorphism $x\mapsto-x$ of $\S$ (see
Proposition~\ref{negmap}). Furthermore, for $\la\in P$ let
$\la^*=w_0(-\la)$, where $w_0$ is the longest element of $W_0$.
Finally, for $l\in I_P$ let $l^*=\tau(\la^*)$, where $\la\in P$
is any vertex with $\tau(\la)=l$.
\end{defn}

We need to check that the definition of $l^*$ is unambiguous. If
$\tau(\la)=\tau(\mu)$, then $\la=w\mu$ for some $w\in W$. Since
$W=W_0\ltimes Q$ we have $w=w't_{\gamma}$ for some $w'\in W_0$
and $\gamma\in Q$, and so
$-\la=-w'(\gamma+\mu)=w't_{-\gamma}(-\mu)=w''(-\mu)$ for some
$w''\in W$. Thus $\tau(-\la)=\tau(-\mu)$, and so
$\tau(\la^*)=\tau(\mu^*)$.

Note that in general $\s_*$ is not an element of $\Auttr(D)$. In
the $BC_n$ case, $\s_*$ is the identity, for the map $x\mapsto-x$
fixes the good type $0$, implying that $\s_*=\mathrm{id}$ by
direct consideration of the Coxeter graph.

\begin{prop}\label{prop:3} If $\la\in P^+$, then $\la^*\in P^+$.
\end{prop}
\begin{proof}
Observe that $w_0(-\cs_0)=\cs_0$ since $-\cs_0$ is a sector that
lies on the opposite side of every wall to $\cs_0$. Thus
$w_0(-\la)\in P^+$.
\end{proof}

\subsection{Special Group Elements and Technical Results}\label{wla}
For $i\in I$, let $W_i=W_{I\backslash\{i\}}$ (this extends our
notation for $W_0$). Given $\la\in P^+$, define $t_{\la}'$ to be
the unique element of $W$ such that $t_{\la}=t_{\la}'g$ for some
$g\in G$, and, using \cite[VI, \ts1, Exercise~3]{bourbaki}, let
$w_{\la}$ be the unique minimum length representative of the
double coset $W_0t_{\la}'W_{l}$, where $l=\tau(\la)$. Fix a
reduced word $f_{\la}\in I^*$ such that $s_{f_{\la}}=w_{\la}$.

\begin{prop}\label{coll} Let $\la\in P^+$ and $i\in I_{P}$. Suppose
that $\tau(\la)=l$, and write $j=\s_i(l)$. Then $g_j=g_ig_l$ and
$t_{\la}=t_{\la}'g_l$.
\end{prop}

\begin{proof} We see that $g_j=g_ig_l$ since the image of $0$ under both
functions is the same. Temporarily write
$t_{\la}=t_{\la}'g_{\la}$, and so
$g_{\la}={t_{\la}'}^{-1}t_{\la}$. Observe that $g_{\la}(0)=v_k$
for some $k\in I_P$ (here $v_k$ is the type $k$ vertex of $C_0$).
But $({t_{\la}'}^{-1}t_{\la})(0)={t_{\la}'}^{-1}(\la)=v_l$, since
$t_{\la}'$ is type preserving. Thus $v_k=v_l$, so $k=l$, and so
$g_{\la}=g_l$.  \end{proof}

Recall that $\s\in\Aut(D)$ induces an automorphism (which we also
denote by $\s$) of $W$ as in (\ref{sigmaaction}). From
(\ref{pres}) we have the following.

\begin{lem}\label{end} Let $\la\in P$ and $l=\tau(\la)$.
Then $g_lW_0g_l^{-1}=W_l=\s_l(W_0)$, and so $W_l$ is the
stabiliser of the type $l$ vertex $v_l$ of $C_0$.
\end{lem}

\begin{prop}\label{closest} Let $\la\in P^+$. Then
\begin{enumerate}
\item[$\mathrm{(i)}$] $w_{\la}=t_{\la}w_{0\la}w_0g_l^{-1}=t_{\la}'\s_l(w_{0\la}w_0)$, where
$l=\tau(\la)$, and $w_{0\la}$ and $w_0$ are the longest elements
of $W_{0\la}$ and $W_0$ respectively.
\item[$\mathrm{(ii)}$] $\la\in w_{\la}\overline{C}_0$, and $w_{\la}C_0$ is
the unique chamber nearest $C_0$ with this property,
\item[$\mathrm{(iii)}$] $w_{\la}C_0\subseteq \cs_0$.
\end{enumerate}
\end{prop}

\begin{proof} (i) By Proposition~\ref{coll} and Lemma~\ref{end} we
have $W_0t_{\la}W_0=W_0t_{\la}'g_lW_0=W_0t_{\la}'W_lg_l$, and so
the double coset $W_0t_{\la}W_0$ has unique minimal length
representative $m_{\la}=w_{\la}g_l$. By \cite[(2.4.5)]{m} (see
also \cite[(2.16)]{ram}) we have $m_{\la}=t_{\la}w_{0\la}w_0$,
proving the first equality in (i). Then
$$
w_{\la}=m_{\la}g_l^{-1}=t_{\la}w_{0\la}w_0g_l^{-1}=t_{\la}'g_lw_{0\la}w_0g_l^{-1}=t_{\la}'\s_l(w_{0\la}w_0).
$$

(ii) With $m_{\la}$ as above we have
$m_{\la}(0)=(t_{\la}w_{0\la}w_0)(0)=\la$, so $\la\in
m_{\la}\overline{C}_0$. Now $w_{\la}=m_{\la}g_l^{-1}$, and since
$g_l^{-1}\in G$ fixes $C_0$ we have $\la\in
w_{\la}\overline{C}_0$.

To see that $w_{\la}$ is the unique chamber nearest $C_0$ that
contains $\la$ in its closure, notice that by Lemma~\ref{end} the
stabiliser of $\la$ in $W$ is $t_{\la}'W_lt_{\la}'^{-1}$, which
acts simply transitively on the set of chambers containing $\la$
in their closure. So if $wC_0$ is a chamber containing $\la$ in
its closure, then
$wC_0=(t_{\la}'w_lt_{\la}'^{-1})t_{\la}'(C_0)=t_{\la}'w_lC_0$ for
some $w_l\in W_l$. Thus $w=t_{\la}'w_l\in t_{\la}'W_l\subset
W_0t_{\la}'W_l$, and so $\ell(w_{\la})\leq \ell(w)$. The
uniqueness follows from \cite[Theorem~2.9]{ronan}.

We now prove (iii). The result is clear if $\la=0$, so let $\la\in
P^+\backslash\{0\}$. If $\la\in\cs_0$ then $\cs_0\cap
w_{\la}C_0\neq\emptyset$, and so $w_{\la}C_0\subseteq \cs_0$ since
$w_{\la}C_0$ is connected and contained in
$E\backslash\bigcup_{H\in \ch_0}H$.

Now suppose that $\la\in\overline{\cs}_0\backslash \cs_0$, so
$\la\in H_{\alpha}$ for some $\alpha\in B$. Let
$C_0,C_1,\ldots,C_m=w_{\la}C_0$ be the gallery of type $f_{\la}$
from $C_0$ to $w_{\la}C_0$. If $w_{\la}C_0\nsubseteq \cs_0$ then
this gallery crosses the wall $H_{\alpha}$, so let $C_{k}$ be the
first chamber on the opposite side of $H_{\alpha}$ to~$C_0$. The
sequence
$C_0,\ldots,C_{k-1},s_{\alpha}(C_k),\ldots,s_{\alpha}(w_{\la}C_0)$
joins $0$ to $\la$ as $s_{\alpha}(\la)=\la$. Since
$C_{k-1}=s_{\alpha}(C_k)$, we can construct a gallery joining $0$ to $\la$ of length strictly less
than $m$, a contradiction.  \end{proof}

Each coset $wW_{0\la}$, $w\in W_0$, has a unique minimal length
representative. To see this, notice that by Lemma~\ref{inW2},
$W_{0\la}$ is the subgroup of $W_0$ generated by $S_{0\la}=\{s\in
S_0\mid s\la=\la\}$, and apply \cite[IV, \ts1,
Exercise~3]{bourbaki}. We write $W_0^{\la}$ for the set of minimal
length representatives of elements of $W_{0}/W_{0\la}$. The
following proposition records some simple facts.

\begin{prop}\label{early2} Let $\la\in P^+$ and write $l=\tau(\la)$. Then
\begin{enumerate}

\item[$\mathrm{(i)}$] $t_{\la}'=w_{\la}w_l$ for some $w_l\in W_l$, and
$\ell(t_{\la}')=\ell(w_{\la})+\ell(w_l)$.

\item[$\mathrm{(ii)}$] Each $w\in W_0$ can be written uniquely as $w=uv$ with $u\in
W_0^{\la}$ and $v\in W_{0\la}$, and moreover
$\ell(w)=\ell(u)+\ell(v)$.

\item[$\mathrm{(iii)}$] For $v\in W_{0\la}$, $vw_{\la}=w_{\la}w_l\s_l(v)w_l^{-1}$
where $w_l\in W_{l}$ is as in (i). Moreover
$$\ell(vw_{\la})=\ell(v)+\ell(w_{\la})=\ell(w_{\la})+\ell(w_l\s_l(v)w_l^{-1}).$$

\item[$\mathrm{(iv)}$] Each $w\in W_0w_{\la}W_l$ can be written uniquely as
$w=uw_{\la}w'$ for some $u\in W_0^{\la}$ and $w'\in W_l$, and
moreover $\ell(w)=\ell(u)+\ell(w_{\la})+\ell(w')$.
\end{enumerate}
\end{prop}

\begin{proof} (i) follows from the proof of Proposition~\ref{closest} and \cite[VI,
\ts1, Exercise~3]{bourbaki}.

(ii) is immediate from the definition of $W_0^{\la}$, and
\cite[VI, \ts1, Exercise~3]{bourbaki}.

(iii) Observe first that $vt_{\la}=t_{\la}v$ in the extended
affine Weyl group, for $vt_{\la}v^{-1}=t_{v\la}$ for all $v\in
W_0$, and $t_{v\la}=t_{\la}$ if $v\in W_{0\la}$. Since
$t_{\la}=t_{\la}'g_{l}$ (see Proposition~\ref{coll}) we have
$$vt_{\la}'=vt_{\la}g_l^{-1}=t_{\la}vg_l^{-1}=t_{\la}'(g_l^{\vphantom{-1}}vg_l^{-1})=t_{\la}'\s_l(v),$$
and so from (i), $vw_{\la}=w_{\la}w_l\s_l(v)w_l^{-1}$. By
\cite[IV, \ts1, Exercise~3]{bourbaki} we have
$\ell(vw_{\la})=\ell(v)+\ell(w_{\la})$; in fact,
$\ell(ww_{\la})=\ell(w)+\ell(w_{\la})$ for all $w\in W_0$.
Observe now that $ws_{\alpha}w^{-1}=s_{w\alpha}$ for $w\in W_0$,
and it follows that
$\ell(w_l^{\vphantom{-1}}\s_l(v)w_{l}^{-1})=\ell(v)$.

(iv) By \cite[IV, \ts1, Exercise~3]{bourbaki} each $w\in
W_0w_{\la}W_l$ can be written as $w=w_1w_{\la}w_2$ for some
$w_1\in W_0$ and $w_2\in W_l$ with
$\ell(w)=\ell(w_1)+\ell(w_{\la})+\ell(w_2)$. Write $w_1=uv$ where
$u\in W_0^{\la}$ and $v\in W_{0\la}$ as in (ii). Then by (iii)
$$w_1w_{\la}w_2=uvw_{\la}w_2=uw_{\la}(w_l^{\vphantom{-1}}\s_l(v)w_l^{-1}w_2),$$
and so each $w\in W_0w_{\la}W_l$ can be written as $w=uw_{\la}w'$
for some $u\in W_0^{\la}$ and $w'\in W_l$ with
$\ell(w)=\ell(u)+\ell(w_{\la})+\ell(w')$. Suppose that we have
two such expressions $w=u_1w_{\la}w'_1=u_2w_{\la}w'_2$ where
$u_1,u_2\in W_0^{\la}$ and $w_1',w_2'\in W_l$. Write $v_l$ for
the type $l$ vertex of $C_0$. Then
$(u_1w_{\la}w_l')(v_l)=(u_1w_{\la})(v_l)=u_1\la$, and similarly
$(u_2w_{\la}w_2')(v_l)=u_2\la$. Thus $u_1^{-1}u_2\in W_{0\la}$,
and so $u_1W_{0\la}=u_2W_{0\la}$, forcing $u_1=u_2$. This clearly
implies that $w_1'=w_2'$ too, completing the proof.  \end{proof}

Recall the definitions of $\s_*$, $\la^*$ and $l^*$ from
Definition~\ref{defn:stardefn}.

\begin{prop}\label{little} Let $\la\in P^+$ (so $\la^*\in P^+$ too), and write $\tau(\la)=l$.
\begin{enumerate}
\item[$\mathrm{(i)}$] $\s_*^2=\mathrm{id}$ and $\s_*(0)=0$.
\item[$\mathrm{(ii)}$] $\s_*(w_{\la})=w_{\la^*}$ and
$\s_*(l)=l^*$.
\item[$\mathrm{(iii)}$] $\s_*\circ\s_i\circ\s_*^{-1}=\s_{i^*}$ for all $i\in
I_P$.
\item[$\mathrm{(iv)}$] $w_{\la^*}=\s_l^{-1}(w_{\la}^{-1})$.
\end{enumerate}
\end{prop}

\begin{proof} (i) is clear, since $-(-x)=x$ for all
$x\in E$.

(ii) Let $\psi$ be the automorphism of $\S$ given by
$\psi(x)=w_0(-x)$ for all $x\in E$. Then the automorphism of $D$
induced by $\psi$ is $\s_*$ (see Proposition~\ref{induction}).
Let $C_0,\ldots,C_m=w_{\la}C_0$ be the gallery of type $f_{\la}$
in $\S$ starting at~$C_0$, and so $\psi(C_0),\ldots,\psi(C_m)$ is
a minimal gallery of type $\s_*(f_{\la})$ (see
Proposition~\ref{prop:addlate}). Observe that $\psi(C_0)=C_0$ and
$\la^*\in\psi(\overline{C}_m)$. The gallery
$\psi(C_0),\ldots,\psi(C_m)$ from $C_0$ to $\la^*$ cannot be
replaced by any shorter gallery joining $C_0$ and $\la^*$, for if
so, by applying $\psi^{-1}$ we could obtain a gallery from $C_0$
to $\la$ of length $<\ell(w_{\la})$. Thus $\psi(C_m)=C_{\la^*}$
by Proposition~\ref{closest}, and so $\s_*(f_{\la})\sim
f_{\la^*}$. Therefore $\s_*(w_{\la})=w_{\la^*}$, and so
$\s_*(l)=l^*$.

(iii) Since $\Auttr(D)$ is normal in $\Aut(D)$ (see
(\ref{gammanormal})) we know that
$\s_*\circ\s_i\circ\s_*^{-1}=\s_k$ for some $k\in I_P$. By (i) and
(ii) we have $(\s_*\circ\s_i\circ\s_*^{-1})(0)=i^*$ and the
result follows.

(iv) Let $C_0,\ldots,C_m$ be the gallery from (ii) and write
$f_{\la}=i_1\cdots i_m$. Then $C_m,\ldots,C_0$ is a gallery of
type $\mathrm{rev}(f_{\la})=i_m\cdots i_1$ joining $\la$ to $0$.
Let $\psi=w_0\circ w_{0\la}^{-1}\circ t_{-\la}:\S\to\S$ where
$w_{0\la}$ is the longest element of $W_{0\la}$. By
Proposition~\ref{closest}(i) we have
$$
\psi(C_m)=(w_0\circ w_{0\la}^{-1}\circ t_{-\la}\circ
w_{\la})(C_0)=C_0.
$$
Thus by Proposition~\ref{prop:addlate}
$C_0=\psi(C_m),\ldots,\psi(C_0)$ is a gallery of type
$\s_l^{-1}(\mathrm{rev}(f_{\la}))$ joining $0$ to $\la^*$ (since
$\la^*\in\psi(\overline{C}_0)$). Since no shorter gallery joining
$0$ to $\la^*$ exists (for if so apply $\psi^{-1}$ to obtain a
contradiction) it follows that
\begin{align*}
w_{\la^*}&=\s_l^{-1}(s_{\mathrm{rev}(f_{\la})})=\s_l^{-1}(s_{f_{\la}}^{-1})=\s_l^{-1}(w_{\la}^{-1})\,.\qedhere
\end{align*}
\end{proof}

\subsection{Affine Buildings}\label{affinebuildings} A building $\scx$ is called \textit{affine} if the associated
Coxeter group $W$ is an affine Weyl group. To study the algebra
$\sca$ of the next section, it is convenient to associate a root
system $R$ to each irreducible locally finite regular affine
building. If $\scx$ is of type $W$, we wish to choose $R$ so that
(among other things) (i) the affine Weyl group of $R$ is
isomorphic to $W$, and (ii) $q_{\s(i)}=q_i$ for all $i\in I$ and
$\s\in\Auttr(D)$ (note that $\Auttr(D)$ depends on the choice of
$R$, see (\ref{auttrd})).

It turns out (as should be expected) that the choice of $R$ is in
most cases straight forward; for example, if $\scx$ is of type
$\widetilde{F}_4$ then choose $R$ to be a root system of
type~$F_4$ (and call $\scx$ an \textit{affine building of type
$F_4$}). The regular buildings of types $\widetilde{A}_1$ and
$\widetilde{C}_n$ $(n\geq2)$ are the only exceptions to this
rule, and in these cases the non-reduced root systems $BC_n$
$(n\geq1$) play an important role. Let us briefly describe why.

Using Proposition~\ref{well}(ii) we see that the parameters of a
regular $\widetilde{C}_n$ ($n\geq2$) building must be as follows:
\begin{figure}[ht]
 \begin{center}
 \psset{xunit= 0.9 cm,yunit= 0.9 cm}
 \psset{origin={0,0}}
\vspace{0.3cm}

\rput(-3,0.3){$q_0$}\pscircle*(-3,0){2pt}
\rput(-2,0.3){$q_1$}\pscircle*(-2,0){2pt} \psline(-3,0)(-2,0)
\rput(-2.5,0.3){$4$} \rput(-1,0.3){$q_1$}\pscircle*(-1,0){2pt}
\psline(-2,0)(-1,0) \pscircle*(0.4,0){1.5pt}
\pscircle*(-0.4,0){1.5pt} \pscircle*(0,0){1.5pt}
\rput(1,0.3){$q_1$}\pscircle*(1,0){2pt}
\rput(2,0.3){$q_1$}\pscircle*(2,0){2pt} \psline(1,0)(3,0)
\rput(3,0.3){$q_n$}\pscircle*(3,0){2pt} \rput(2.5,0.3){$4$}
\end{center}
\end{figure}

\noindent If we choose $R$ to be a $C_n$ root system then the
automorphism $\s_n\in\Auttr(D)$ interchanges the left most and
right most nodes, and so condition (ii) is not satisfied (unless
$q_0=q_n$). If, however, we take $R$ to be a $BC_n$ root system,
then $\Auttr(D)=\{\mathrm{id}\}$, and so both conditions (i) and
(ii) are satisfied.

Thus, in order to facilitate the statements of later results, we
rename regular $\widetilde{C}_n$ $(n\geq2)$ buildings, and call
them \textit{affine buildings of type $BC_n$} (or $\bct_n$ ($n\geq2$) buildings). We
reserve the name `$\widetilde{C}_n$ building' for the special
case when $q_0=q_n$ in the above parameter system. For a similar
reason we rename regular $\widetilde{A}_1$ buildings (which are
\textit{semi-homogeneous trees}) and call them \textit{affine
buildings of type $BC_1$} (or $\bct_1$ buildings), and reserve the
name `$\widetilde{A}_1$ building' for homogeneous trees. With
these conventions we make the following definitions.

\begin{defn}\label{goodbuilding} Let $\scx$ be an affine building of type~$R$ with
vertex set~$V$, and let~$\S=\S(R)$. Let $V_{\mathrm{sp}}(\S)$
denote the set of all special vertices of~$\S$ (see
Section~\ref{laterefad}), and let
$I_{\mathrm{sp}}=\{\tau(\la)\mid \la\in V_{\mathrm{sp}}(\S)\}$.
\begin{enumerate}
\item[(i)] A vertex $x\in V$ is said to be \textit{special}\index{vertex!of $\scx$!special} if
$\tau(x)\in I_{\mathrm{sp}}$. We write $V_{\mathrm{sp}}$ for the
set of all special vertices of~$\scx$.
\item[(ii)] A vertex $x\in V$ is said to be \textit{good}\index{vertex!of $\scx$!good} if
$\tau(x)\in I_P$, where $I_P$ is as in Section~\ref{laterefad}.
We write $V_P$ for the set of all good vertices of~$\scx$.
\end{enumerate}
\end{defn}

Clearly $V_P\subset V_{\mathrm{sp}}$. In fact if $R$ is reduced,
then by the comments made in Section~\ref{laterefad},
$V_P=V_{\mathrm{sp}}$. If $R$ is non-reduced (so $R$ is of type
$BC_n$ for some $n\geq1$), then $V_P$ is the set of all type~$0$
vertices of $\scx$,
whereas $V_{\mathrm{sp}}$ is the set of all type~$0$ and type~$n$
vertices of~$\scx$.

\begin{prop} A vertex $x\in V$ is good if and only if there exists an apartment~$\ca$ containing~$x$ and a type
preserving isomorphism~$\psi:\ca\to\S$ such that $\psi(x)\in P$.
\end{prop}

\begin{proof} Let $x\in V_P$, and choose any apartment~$\sca$
containing~$x$. Let $\psi:\ca\to\S$ be a type preserving
isomorphism (from the building axioms). Then $\psi(x)$ is a
vertex in~$\S$ with type~$\tau(x)\in I_P$, and so~$\psi(x)\in P$.
The converse is obvious.
\end{proof}

\begin{rem}\label{IDRG} We note that \textit{infinite distance
regular graphs} are just $\bct_1$ buildings in very thin
disguise. To see the connection, given any $p,q\geq1$, construct
a $\bct_1$ building (that is, a semi-homogeneous tree) with
parameters $q_0=p$ and $q_1=q$. Construct a new graph
$\Gamma_{p,q}$ with vertex set $V_P$ and vertices $x,y\in V_P$
connected by an edge if and only if $d(x,y)=2$. It is simple to
see that $\Gamma_{p,q}$ is the (graph) free product
$\mathbb{K}_q*\cdots*\mathbb{K}_q$ ($p$ copies) where
$\mathbb{K}_q$ is the complete graph on $q$ letters. By the
classification (\cite{ivanov}, \cite{macpherson}) $\Gamma_{p,q}$
is infinite distance regular, and all infinite distance regular
graphs occur in this way.
\end{rem}

Recall the definition of $\Aut_q(D)$ from (\ref{autoq}).

\begin{thm}\label{pictures} The diagrams in the Appendix characterise the parameter
systems of the locally finite regular affine buildings. In each
case $\Auttr(D)\cup\{\s_*\}\subseteq\Aut_q(D)$.
\end{thm}

\begin{proof} These parameter systems are found case by case using Proposition~\ref{well}(ii)
and the classification of the irreducible affine Coxeter graphs.
Note that $\Auttr(D)\cup\{\s_*\}=\{\mathrm{id}\}$ if $\scx$ is a
$\bct_n$ building. Thus the final result follows by considering
each Coxeter graph.  \end{proof}

\section{Vertex Set Operators and Vertex
Regularity}\label{section5}

Let $\scx$ be a locally finite regular affine building of type~$R$ (see
Section~\ref{affinebuildings}). Recall
(Definition \ref{goodbuilding}) that we write $V_P$ for the set of all good
vertices of~$\scx$.

For
each $\la\in P^+$ we will define an averaging operator $A_{\la}$
acting on the space of all functions $f:V_P\to\bc$, and we will
introduce an algebra $\sca$ of these operators. The operators
$A_{\la}$ were defined in \cite[II, \ts1.1.2, Excercise 3]{serre}
for homogeneous trees, \cite{cm} and \cite{MZ} for
$\widetilde{A}_2$ buildings, and \cite{C2} for $\widetilde{A}_n$
buildings. Our definition gives the generalisation of the
operators $A_{\la}$ and the algebra $\sca$ to all (irreducible)
root systems.

\subsection{Initial Observations}

Recall the definition of type preserving isomorphisms of
simplicial complexes.

\begin{defn}\label{typepreserving} Let $\ca_1$ and $\ca_2$ be apartments of
$\scx$.
\begin{enumerate}
\item[$\mathrm{(i)}$] An isomorphism $\psi:\ca_1\to\ca_2$ is called
\textit{type-rotating} if it is of the form
$\psi=\psi_2^{-1}\circ w\circ\psi_1$ where $\psi_1:\ca_1\to\S$
and $\psi_2:\ca_2\to\S$ are type preserving isomorphisms, and
$w\in \tilde{W}$.
\item[$\mathrm{(ii)}$] We have an analogous definition to (i) for isomorphisms
$\psi:\ca_1\to\S$ by omitting $\psi_2$.
\end{enumerate}
\end{defn}

\begin{prop}\label{useful} Let $\ca,\ca'$ be any apartments and suppose that
$\psi:\ca\to\ca'$ is an isomorphism. Then
\begin{enumerate}
\item[$\mathrm{(i)}$] The image under $\psi$ of a gallery in $\ca$ is a gallery
in~$\ca'$.
\item[$\mathrm{(ii)}$] A gallery in $\ca$ is minimal if and only if its image under
$\psi$ is minimal in~$\ca'$.
\item[$\mathrm{(iii)}$] There exists a unique
$\s\in\Aut(D)$ so that $\psi$ maps galleries of type $f$ in $\ca$
to galleries of type $\s(f)$ in $\ca'$. If $\psi$ is type
rotating, then $\s\in\Auttr(D)$, and
$(\tau\circ\psi)(x)=(\s\circ\tau)(x)$ for all vertices $x$ of
$\ca$.
\item[$\mathrm{(iv)}$] If $\psi$ is type rotating and maps a type $i\in I_P$ vertex in $\ca$ to a type $j\in I_P$
vertex in $\ca'$, then the induced automorphism from~(iii) is
$\s=\s_j\circ\s_i^{-1}$.
\end{enumerate}
\end{prop}

\begin{proof} This follows from Proposition~\ref{prop:addlate} and the definition of type rotating isomorphisms. \end{proof}

\begin{lem}\label{quicker} Suppose $x\in V_P$ is contained in the apartments $\ca$ and
$\ca'$ of $\scx$, and suppose that $\psi:\ca\to\S$ and
$\psi':\ca'\to\S$ are type rotating isomorphisms such that
$\psi(x)=0=\psi'(x)$. Let $\psi'':\ca\to\ca'$ be a type
preserving isomorphism mapping $x$ to $x$ (the existence of which
is guaranteed by Definition \ref{def2}). Then
$\phi=\psi'\circ\psi''\circ\psi^{-1}$ is in $W_0$.
\end{lem}

\begin{proof} Observe that $\phi:\S\to\S$ has $\phi(0)=0$. Since
$\psi$ and $\psi'$ are type rotating isomorphisms we have
$\psi=w\circ\psi_1$ and $\psi'=w'\circ\psi_1'$ for some
$w,w'\in\tilde{W}$ and $\psi_1:\ca\to\S$, $\psi_1':\ca'\to\S$
type preserving isomorphisms. Therefore
$$\phi=w'\circ\psi_1'\circ\psi''\circ\psi_1^{-1}\circ w^{-1}=w'\circ\phi'\circ w^{-1}\,,\quad\textrm{say}\,.$$
Now $\phi'=\psi_1'\circ\psi''\circ\psi_1^{-1}:\S\to\S$ is a type
preserving automorphism, as it is a composition of type
preserving isomorphisms. By \cite[Lemma~2.2]{ronan} we have
$\phi'=v$ for some $v\in W$, and hence $\phi=w'\circ v\circ
w^{-1}\in\tilde{W}$. Since $\phi(0)=0$ and $\tilde{W}=W_0\ltimes
P$ we in fact have $\phi\in W_0$, completing the proof.
\end{proof}

\subsection{The Sets $V_{\la}(x)$} The following definition gives the analogue of the partition
$\{\cc_w(a)\}_{w\in W}$ used for the chamber set of $\scx$. Let
us first record the following lemma from \cite[p.24]{brown} (or
\cite[\ts10.3, Lemma~B]{h2}). Recall the definition of the
fundamental sector $\cs_0$ from (\ref{fundamentalsector}).

\begin{lem}\label{inW2}
Let $w\in W_0$ and $\la\in E$. If
$\la'=w\la\in\overline{\cs}_0\cap w\overline{\cs}_0$ then
$\la'=\la$, and $w\in\lan \{s_i\mid s_i\la=\la\}\ran$.
\end{lem}

\begin{defn}\label{vertices in sigma} Given $x\in V_P$ and $\la\in P^+$, we define $V_{\la}(x)$ to be the set
of all $y\in V_P$ such that there exists an apartment $\ca$
containing $x$ and $y$ and a type rotating isomorphism
$\psi:\ca\to\S$ such that $\psi(x)=0$ and $\psi(y)=\la$.
\end{defn}

\goodbreak\begin{prop}\label{to} Let $V_{\la}(x)$ be as in
Definition~\ref{vertices in sigma}.
\begin{enumerate}
\item[(i)] Given $x,y\in V_P$, there exists some $\la\in P^+$ such
that $y\in V_{\la}(x)$.
\item[(ii)] If $y\in V_{\la}(x)\cap V_{\la'}(x)$ then $\la=\la'$.
\item[(iii)] Let $y\in V_{\la}(x)$. If $\ca$ is any apartment containing $x$ and $y$,
then there exists a type-rotating isomorphism $\psi:\ca\to\S$
such that $\psi(x)=0$ and $\psi(y)=\la$.
\end{enumerate}
\end{prop}

\begin{proof} First we prove (i). By Definition \ref{def2} there exists an apartment
$\ca$ containing $x$ and $y$ and a type preserving isomorphism
$\psi_1:\ca\to\S$. Let $\mu=\psi_1(x)$ and $\nu=\psi_1(y)$, so
$\mu,\nu\in P$. There exists a $w\in W_0$ such that
$w(\nu-\mu)\in\overline{\cs}_0\cap P$ \cite[p.55, exercise
14]{h2}, and so the isomorphism $\psi=w\circ t_{-\mu}\circ\psi_1$
satisfies $\psi(x)=0$ and $\psi(y)=w(\nu-\mu)\in P^+$,
proving~(i).

We now prove (ii). Suppose that there are apartments $\ca$ and
$\ca'$ containing $x$ and $y$, and type-rotating isomorphisms
$\psi:\ca\to\S$ and $\psi':\ca'\to\S$ such that
$\psi(x)=\psi'(x)=0$ and $\psi(y)=\la\in P^+$ and
$\psi'(y)=\la'\in P^+$. We claim that $\la=\la'$.

By Definition \ref{def2}(iii)$'$ there exists a type preserving
isomorphism $\psi'':\ca\to\ca'$ which fixes $x$ and $y$. Then
$\phi=\psi'\circ\psi''\circ\psi^{-1}:\S\to\S$ is a type-rotating
automorphism of $\S$ that fixes $0$ and maps $\la$ to $\la'$. By
Lemma~\ref{quicker} we have $\phi=w$ for some $w\in W_0$, and so
we have $\la'=w\la\in\overline{\cs}_0\cap w\overline{\cs}_0$.
Thus by Lemma~\ref{inW2} we have $\la'=\la$.

Note first that (iii) is not immediate from the definition of
$V_{\la}(x)$. To prove (iii), by the definition of $V_{\la}(x)$
there exists an apartment $\ca'$ containing $x$ and $y$, and a
type-rotating isomorphism $\psi':\ca'\to\S$ such that
$\psi'(x)=0$ and $\psi'(y)=\la$. Then by
Definition~\ref{def2}(iii)$'$ there is a type preserving
isomorphism $\phi:\ca\to\ca'$ fixing $x$ and $y$. Then
$\psi=\psi'\circ\phi:\ca\to\S$ has the required properties.
\end{proof}

\begin{rem} Note that the assumption that $\psi$ is
type-rotating in Definition~\ref{vertices in sigma} is essential
for Proposition~\ref{to}(ii) to hold. To see this we only need to
look at an apartment of an $\widetilde{A}_2$ building. The map
$a_1\la_1+a_2\la_2\mapsto a_1\la_2+a_2\la_1$ is an automorphism which maps $\la_1$ to
$\la_2$. Thus if we omitted the hypothesis that $\psi$ is
type-rotating in Definition \ref{vertices in sigma}, part (ii) of
Proposition~\ref{to} would be false.
\end{rem}

\begin{prop}\label{ease} If $y\in V_{\la}(x)$, then $x\in V_{\la^*}(y)$ where
$\la^*$ is as in Definition~\ref{defn:stardefn}.
\end{prop}

\begin{proof}
If $\psi:\ca\to\S$ is a type rotating isomorphism mapping $x$ to
$0$ and $y$ to $\la$, then $w_0\circ t_{-\la}\circ\psi:\ca\to\S$
is a type rotating isomorphism mapping $y$ to $0$ and $x$ to
$\la^*=w_0(-\la)\in P^+$ (see Proposition~\ref{prop:3}).
\end{proof}

\begin{lem} Let $x\in V_P$ and $\la\in P^+$. If $y,y'\in V_{\la}(x)$
then $\tau(y)=\tau(y')$.
\end{lem}

\begin{proof} Let $\ca$ be an apartment containing $x$ and $y$,
and $\ca'$ be an apartment containing $x$ and $y'$. Let
$\psi:\ca\to\S$ and $\psi':\ca'\to\S$ be type rotating
isomorphisms with $\psi(x)=\psi'(x)=0$ and
$\psi(y)=\psi'(y')=\la$. Thus
$\chi=\psi'^{-1}\circ\psi:\ca\to\ca'$ is a type preserving
automorphism since $\chi(x)=x$ (see Proposition~\ref{lete}). Since
$\chi(y)=y'$ we have $\tau(y)=\tau(y')$.  \end{proof}

In light of the above lemma we define $\tau(V_{\la}(x))=\tau(y)$
for any $y\in V_{\la}(x)$. \vspace{0.1cm}

Clearly the sets $V_{\la}(x)$ are considerably more complicated
objects than the sets $\cc_w(a)$. The following theorem provides
an important connection between the sets $V_{\la}(x)$ and
$\cc_w(a)$ that will be relied on heavily in subsequent work.
Given a chamber $c\in\cc$ and an index $i\in I$, we define
$\pi_i(c)$ to be the type $i$ vertex of $c$. For the following
theorem the reader is reminded of the definition of $w_{\la}\in
W$ and $f_{\la}\in I^*$ from Section \ref{wla}.

\begin{thm}\label{key} Let $x\in V_P$ and $\la\in P^+$.
Suppose $\tau(x)=i$ and $\tau(V_{\la}(x))=j$, and let $a\in\cc$ be
any chamber with $\pi_i(a)=x$. Then
$$\{b\in\cc:\pi_{j}(b)\in V_{\la}(x)\}= \bigcup_{w\in W_i\s_i(w_{\la})W_j}\cc_{w}(a)\,,$$
where the union is disjoint.
\end{thm}

\begin{proof} Suppose first that $y=\pi_{j}(b)\in V_{\la}(x)$. Let
$a=c_0,c_1,\ldots,c_n=b$ be a minimal gallery from $a$ to $b$ of
type $f$, say. By \cite[Theorem~3.8]{ronan}, all the $c_{k}$ lie
in some apartment, $\ca$, say. Let $\psi:\ca\to\S$ be a type
rotating isomorphism such that $\psi(x)=0$ and $\psi(y)=\la$. Then
$\psi(c_0),\psi(c_1),\ldots,\psi(c_n)$ is a minimal gallery of
type $\s_i^{-1}(f)$ by Proposition~\ref{useful}.

Recall the definition of the fundamental chamber $C_0$ from
(\ref{fundamentalalcove}). Since $0$ is a vertex of $\psi(c_0)$,
we can construct a gallery from $\psi(c_0)$ to
$C_0$ of type $e_1$, say, where $s_{e_1}\in W_0$. Similarly there
is a gallery from $w_{\la}C_0$ to $\psi(c_n)$ of
type $e_2$, where $s_{e_2}\in W_{\s_i^{-1}(j)}$. Thus we have a
gallery
$$\psi(c_0)\xrightarrow{e_1} C_0\xrightarrow{f_{\la}} w_{\la}C_0\xrightarrow{e_2}\psi(c_n)$$
of type $e_1f_{\la}e_2$. Since $\S$ is a
Coxeter complex, galleries (reduced or not) from
one chamber to another of types $f_1$ and $f_2$, say, satisfy
$s_{f_1}=s_{f_2}$ \cite[p.12]{ronan}, so
$s_{\s_i^{-1}(f)}=s_{e_1f_{\la}e_2}$. Thus
$$\delta(a,b)=s_f=\s_i(s_{\s_i^{-1}(f)})=\s_i(s_{e_1f_{\la}e_2})=s_{e_1'}s_{\s_i(f_{\la})}s_{e_2'}$$
where $e_1'\in W_i$ and $e_2'\in W_{j}$. Thus $b\in\cc_w(a)$ for
some $w\in W_i\s_i(w_{\la})W_j$.

Now suppose that $b\in\cc_w(a)$ for some $w\in
W_i\s_i(w_{\la})W_j$. Let $y=\pi_j(b)$. By \cite[p.35,
Exercise~1]{ronan}, there exists a gallery of type
$e_1'\s_i(f_{\la})e_2'$ from $a$ to $b$ where $e_1'\in W_i$ and
$e_2'\in W_{j}$. Let $c_{k},c_{k+1},\ldots,c_{l}$ be the
subgallery of type $\s_i(f_{\la})$. Note that $\pi_i(c_{k})=x$ and
$\pi_{j}(c_{l})=y$. Observe that $\s_i(f_{\la})$ is reduced since
$\s_i\in\Aut(D)$, and so all of the chambers $c_{m}$, $k\leq m\leq
l$, lie in an apartment $\ca$, say. Let $\psi:\ca\to\S$ be a type
rotating isomorphism such that $\psi(x)=0$. Thus
$\psi(c_k),\ldots,\psi(c_{l})$ is a gallery of
type $f_{\la}$ in $\S$ (Proposition~\ref{useful}). Since $W_0$
acts transitively on the chambers $C\in\cc(\S)$ with $0\in
\overline{C}$ (Lemma~\ref{rt}) there exists $w\in W_0$ such that
$w(\psi(c_k))=C_0$. Then $\psi'=w\circ\psi:\ca\to\S$ is a type
rotating isomorphism that takes the gallery
$c_k,\ldots,c_{l}$ in $\ca$ of type $\s_i(f_{\la})$ to a
gallery $C_0,\ldots ,\psi'(c_{l})$ of type
$f_{\la}$. But in a Coxeter complex there is only one
gallery of each type. So $\psi'(c_{l})$ must be
$w_{\la}(C_0)$, and by considering types $\psi'(y)=\la$, and so
$y\in V_{\la}(x)$.  \end{proof}

For $x\in V$ we write $\mathrm{st}(x)$ for the set of all
chambers that have $x$ as a vertex. Recall the definition of
Poincar\'{e} polynomials from Definition~\ref{defn:poincarep}.

\begin{lem}\label{W_0(q)} Let $x\in V_P$. Then
$|\mathrm{st}(x)|=W_0(q)$. In particular, this value is
independent of the particular $x\in V_P$.
\end{lem}

\begin{proof} Suppose $\tau(x)=i\in I_P$ and let $c_0$ be any chamber that has $x$ as a vertex.
Then
\begin{align*}
\mathrm{st}(x)=\{c\in\cc\mid\delta(c_0,c)\in W_i\}=\bigcup_{w\in
W_{i}}\cc_w(c_0)
\end{align*}
where the union is disjoint, and so $|\mathrm{st}(x)|=\sum_{w\in
W_i}q_w$. Theorem~\ref{pictures} now shows that
\begin{align*}
|\st(x)|&=\sum_{w\in W_0}q_{\s_i(w)}=\sum_{w\in
W_0}q_w=W_0(q).\qedhere
\end{align*}
\end{proof}

Note that if the hypothesis `let $x\in V_P$' in Lemma~\ref{W_0(q)}
is replaced by the hypothesis `let $x$ be a special vertex', then
in the non-reduced case it is no longer true in general that
$|\mathrm{st}(x)|=W_0(q)$.

\subsection{The Cardinalities $|V_{\la}(x)|$} In this subsection
we will find a closed form for $|V_{\la}(x)|$. We need to return
to the operators $B_w$ introduced in Section \ref{section3}.

For each $i\in I$ define an element $\ooo_i\in\scb$ by
\begin{align}
\label{ooo}\ooo_i=\frac{1}{W_i(q)}\sum_{w\in W_i}q_wB_w\,.
\end{align}

\begin{lem}\label{early1} Let $i\in I$. Then $\mathds{1}_iB_w=B_w\mathds{1}_i=\mathds{1}_i$ for all
$w\in W_i$, and $\ooo_i^2=\ooo_i$.
\end{lem}

\begin{proof} Suppose $s$ is a generator of $W_i$ and set
$W_i^{\pm}=\{w\in W_i\mid\ell(ws)=\ell(w)\pm1\}$. Then
\begin{align*}
W_i(q)\ooo_i B_s&=\sum_{w\in W_i^+}q_wB_{ws}+\sum_{w'\in
W_i^-}q_w\left(\frac{1}{q_s}B_{ws}+\left(1-\frac{1}{q_s}\right)B_w\right)\\
&=\sum_{w\in W_i^-}\frac{q_w}{q_s}B_{w}+\sum_{w'\in
W_i^-}q_w\left(\frac{1}{q_s}B_{ws}+\left(1-\frac{1}{q_s}\right)B_w\right)\\
&=\sum_{w\in W_i^-}\left(\frac{q_w}{q_s}B_{ws}+q_w B_w\right)\\
&=\sum_{w\in W_i^+}q_wB_w+\sum_{w\in W_i^-}q_wB_w=W_i(q)\ooo_i\,.
\end{align*}
A similar calculation works for $B_s\ooo_i$ too. It follows that
$\ooo_i B_w=B_w\ooo_i=\ooo_i$ for all $w\in W_i$ and so
$\ooo_i^2=\ooo_i$.  \end{proof}

Recall the definition of $W_{0\la}$ from (\ref{earlier?}).

\begin{thm}\label{early3} Let $\la\in P^+$ and write $l=\tau(\la)$. Then
$$\sum_{w\in
W_0w_{\la}W_l}q_wB_w=\frac{W_0^2(q)}{W_{0\la}(q)}q_{w_{\la}}\ooo_0B_{w_{\la}}\ooo_l.$$
\end{thm}

\begin{proof} Recall from Corollary \ref{brec} That
$B_{w_1}B_{w_2}=B_{w_1w_2}$ whenever
$\ell(w_1w_2)=\ell(w_1)+\ell(w_2)$. Then by
Proposition~\ref{early2}(ii), Proposition~\ref{early2}(iii),
Lemma~\ref{early1} and Proposition~\ref{early2}(iv) (in that
order)
\begin{align*}
\ooo_0B_{w_{\la}}\ooo_l&=\frac{1}{W_0(q)}\sum_{u\in
W_0^{\la}}\sum_{v\in
W_{0\la}}q_uq_vB_uB_vB_{w_{\la}}\ooo_l\\
&=\frac{1}{W_0(q)}\sum_{u\in W_0^{\la}}\sum_{v\in
W_{0\la}}q_uq_vB_uB_{w_{\la}}B_{w_l^{\vphantom{-1}}\s_l(v)w_l^{-1}}\ooo_l\\
&=\frac{1}{W_0(q)}\sum_{u\in W_0^{\la}}\sum_{v\in
W_{0\la}}q_uq_vB_uB_{w_{\la}}\ooo_l\\
&=\frac{W_{0\la}(q)}{W_0(q)W_l(q)}q_{w_{\la}}^{-1}\sum_{w\in
W_0w_{\la}W_l}q_wB_w,
\end{align*}
and the result follows, since
$$W_l(q)=\sum_{w\in W_l}q_{w}=\sum_{w\in W_0}q_{\s_l(w)}=W_0(q)$$ by Proposition~\ref{pictures}.
 \end{proof}

\begin{lem}\label{lit} Let $\la\in P^+$, $x\in V_P$, and $y\in V_{\la}(x)$.
Write $\tau(x)=i$, $\tau(y)=j$ and $\tau(\la)=l$. Then
$\s_i^{-1}(j)=l$, and so $\s_j=\s_i\circ\s_l$.
\end{lem}

\begin{proof} Since $y\in V_{\la}(x)$, there exists an apartment $\ca$ containing $x$ and $y$ and
a type rotating isomorphism $\psi:\ca\to\S$ such that $\psi(x)=0$
and $\psi(y)=\la$. Since $\psi(x)=0$, the $\s$ from
Proposition~\ref{useful}(iii) maps $i$ to $0$ and so is
$\s_i^{-1}$. Thus $\la=\psi(y)$ has type $\s(j)=\s_i^{-1}(j)$ and
so $l=\s_i^{-1}(j)$. Thus $\s_j(0)=(\s_i\circ\s_l)(0)$, and so
$\s_j=\s_i\circ\s_l$.
\end{proof}

\begin{thm}\label{N} Let $x\in V_P$ and $\la\in P^+$ with $\tau(\la)=l\in I_P$. Then
$$|V_{\la}(x)|=\frac{1}{W_0(q)}\sum_{w\in
W_0w_{\la}W_{l}}q_w=\frac{W_0(q)}{W_{0\la}(q)}q_{w_{\la}}=|V_{\la^*}(x)|\,.$$
\end{thm}

\begin{proof} Suppose $\tau(x)=i\in I_P$ and
$\tau(y)=j\in I_P$ for all $y\in V_{\la}(x)$. Let
$\cc_{\la}(x)=\{c\in\cc\mid\pi_{j}(c)\in V_{\la}(x)\}$ and
construct a map $\psi:\cc_{\la}(x)\to V_{\la}(x)$ by $c\mapsto
\pi_{j}(c)$ for all $c\in\cc_{\la}(x)$. Clearly $\psi$ is
surjective.

Observe that for each $y\in V_{\la}(x)$ the set
$\{c\in\cc_{\la}(x)\mid \psi(c)=y\}$ has $|\mathrm{st}(y)|$
distinct elements, and so by Lemma~\ref{W_0(q)} we see that
$\psi:\cc_{\la}(x)\to V_{\la}(x)$ is a $W_0(q)$-to-one
surjection. Let $c_0\in\cc$ be any chamber that has $x$ as a
vertex. Then by the above and Theorem~\ref{key} we have
\begin{align*}
|V_{\la}(x)|&=\frac{|\cc_{\la}(x)|}{W_0(q)}=\frac{1}{W_0(q)}\sum_{w\in
W_i\s_{i}(w_{\la})W_{j}}|\cc_{w}(c_0)|=\frac{1}{W_0(q)}\sum_{w\in
W_i\s_{i}(w_{\la})W_{j}}q_{w}\,.
\end{align*}
Since $\s_i^{-1}(j)=l$ (Lemma~\ref{lit}) we have
$W_i\s_i(w_{\la})W_{j}=\s_i(W_0w_{\la}W_{l})$, and so by
Theorem~\ref{pictures}
$$|V_{\la}(x)|=\frac{1}{W_0(q)}\sum_{w\in
W_0w_{\la}W_{l}}q_{\s_i(w)}=\frac{1}{W_0(q)}\sum_{w\in
W_0w_{\la}W_{l}}q_w\,.$$

Let $1_{\cc}:\cc\to\{1\}$ be the constant function. Then
$(B_w1_{\cc})(c)=1$ for all $c\in\cc$, and so we compute
$(\ooo_l1_{\cc})(c)=1$ for all $c\in\cc$. Thus by
Theorem~\ref{early3}
$$\sum_{w\in
W_0w_{\la}W_l}q_w=\frac{W_0^2(q)}{W_{0\la}(q)}q_{w_{\la}}.$$

Now, by Proposition~\ref{little} and Theorem~\ref{pictures} we
have
\begin{align*}
|V_{\la^*}(x)|&=\frac{1}{W_0(q)}\sum_{w\in\s_*(W_0w_{\la}W_{l})}q_w=\frac{1}{W_0(q)}\sum_{w\in
W_0w_{\la}W_{l}}q_w=|V_{\la}(x)|.\qedhere
\end{align*}
\end{proof}

\begin{defn}
For $\la\in P^+$ we define $N_{\la}=|V_{\la}(x)|$, which is
independent of $x\in V_P$ by Theorem~\ref{N}. 
\end{defn}
By the above we have
$N_{\la}=N_{\la^*}$.

\subsection{The Operators $A_{\la}$ and the Algebra $\sca$} We now
define the \textit{vertex set averaging operators} on $\scx$.

\begin{defn} For each $\la\in P^+$, define an operator
$A_{\la}$, acting on the space of all functions $f:V_P\to\bc$ as
in (\ref{short2}).
\end{defn}

\begin{lem}\label{indep} The operators $A_{\la}$ are linearly independent.
\end{lem}

\begin{proof} Suppose we have a relation
$\sum_{\la\in P^+} a_{\la}A_{\la}=0\,,$ and fix $x,y\in V_P$ with
$y\in V_{\mu}(x)$. Then writing $\delta_y$ for the function
taking the value 1 at $y$ and 0 elsewhere,
$$0=\sum_{\la\in P^+}a_{\la}(A_{\la}\delta_y)(x)=\sum_{\la\in P^+}a_{\la}N_{\la}^{-1}\delta_{\la,\mu}=a_{\mu}N_{\mu}^{-1}\,,$$
and so $a_{\mu}=0$.  \end{proof}

Following the same technique used in (\ref{one}) for the chamber
set averaging operators, we have
\begin{align}
\label{vertexop}(A_{\la}A_{\mu}f)(x)=\frac{1}{N_{\la}N_{\mu}}\sum_{y\in
V_P}|V_{\la}(x)\cap V_{\mu^*}(y)|\,f(y)\qquad\textrm{for all $x\in
V_P$}\,.
\end{align}
Our immediate goal now is to understand the cardinalities
$|V_{\la}(x)\cap V_{\mu^*}(y)|$.

\begin{defn}\label{svr} We say that $\scx$ is \textit{vertex regular} if, for
all $\la,\mu,\nu\in P^+$,
$$|V_{\la}(x)\cap V_{\mu^*}(y)|=|V_{\la}(x')\cap
V_{\mu^*}(y')|\qquad\textrm{whenever $y\in V_{\nu}(x)$ and $y'\in
V_{\nu}(x')$}\,,$$ and \textit{strongly vertex regular} if for
all $\la,\mu,\nu\in P^+$
$$|V_{\la}(x)\cap V_{\mu^*}(y)|=|V_{\la^*}(x')\cap
V_{\mu}(y')|\qquad\textrm{whenever $y\in V_{\nu}(x)$ and $y'\in
V_{\nu^*}(x')$}\,.$$
\end{defn}

Strong vertex regularity implies vertex regularity. To see this,
suppose we are given $x,y,x',y'\in V_P$ with $y\in V_{\nu}(x)$ and
$y'\in V_{\nu}(x')$, and choose any pair $x'',y''\in V_P$ with
$y''\in V_{\nu^*}(x'')$. Then if strong vertex regularity holds,
we have
$$|V_{\la}(x)\cap V_{\mu^*}(y)|=|V_{\la^*}(x'')\cap
V_{\mu}(y'')|=|V_{\la}(x')\cap V_{\mu^*}(y')|\,.$$

\begin{lem}\label{track} Let $y\in V_{\nu}(x)$ and  suppose that $z\in V_{\la}(x)\cap
V_{\mu^*}(y)$. Write $\tau(x)=i$, $\tau(y)=j$, $\tau(z)=k$,
$\tau(\la)=l$, $\tau(\mu)=m$, and $\tau(\nu)=n$.
\begin{enumerate}
\item[(i)] $\s_i^{-1}(k)=l$, $\s_k^{-1}(j)=m$ and
$\s_i^{-1}(j)=n$. Thus $\s_i^{-1}\circ\s_k=\s_l$,
$\s_k^{-1}\circ\s_j=\s_m$ and $\s_i^{-1}\circ\s_j=\s_n$.
\item[(ii)] $\s_n=\s_l\circ\s_m$.
\end{enumerate}
\end{lem}

\begin{proof} (i) follows immediately from Lemma~\ref{lit}. To prove (ii), we have
\begin{align*}
\s_l\circ\s_m&=\s_i^{-1}\circ\s_k\circ\s_k^{-1}\circ\s_j=\s_i^{-1}\circ\s_j=\s_n.\qedhere
\end{align*}
\end{proof}

Recall the definition of the automorphism $\s_*\in\Aut(D)$ from
Section \ref{cp}.

\begin{thm}\label{strongvertexregularity} $\scx$ is strongly vertex regular.
\end{thm}

\begin{proof} Let $x,y\in V_P$ with $y\in V_{\nu}(x)$ and suppose
that $z\in V_{\la}(x)\cap V_{\mu^*}(y)$. Let $\tau(x)=i$,
$\tau(y)=j$ and $\tau(z)=k$. With the notation used in the proof
of Theorem~\ref{N}, define a map
$\psi:\cc_{\la}(x)\cap\cc_{\mu^*}(y)\to V_{\la}(x)\cap
V_{\mu^*}(y)$ by the rule $\psi(c)=\pi_k(c)$. As in the proof of
Theorem~\ref{N} we see that this is a $W_{0}(q)$-to-one
surjection, and thus by Theorem~\ref{key}
$$|V_{\la}(x)\cap
V_{\mu^*}(y)|=\frac{1}{W_0(q)}\sum_{\substack{w_1\in
W_i\s_i(w_{\la})W_k\\ w_2\in
W_j\s_j(w_{\mu^*})W_k}}|\cc_{w_1}(a)\cap\cc_{w_2}(b)|$$ where $a$
and $b$ are any chambers with $\pi_i(a)=x$ and $\pi_j(b)=y$.
Notice that this implies that $\delta(a,b)\in
W_i\s_i(w_{\nu})W_j$, by Theorem~\ref{key}.

Writing $\tau(\la)=l$ and $\tau(\nu)=n$, Lemma~\ref{track}(i)
implies that
$$W_i\s_i(w_{\la})W_k=\s_i(W_0w_{\la}\s_i^{-1}(W_k))=\s_i(W_0w_{\la}W_{\s_i^{-1}(k)})=\s_i(W_0w_{\la}W_{l})\,,$$
$$W_j\s_{j}(w_{\mu^*})W_k=\s_i(W_{\s_i^{-1}(j)}(\s_i^{-1}\circ\s_j)(w_{\mu^*})W_{\s_i^{-1}(k)})=\s_i(W_n\s_n(w_{\mu^*})W_{l})$$
and similarly $W_i\s_i(w_{\nu})W_j=\s_i(W_0w_{\nu}W_n)$. Applying
Lemma~\ref{vertexuse} (with $\s=\s_i$) we therefore have
\begin{align}
\label{vreg}|V_{\la}(x)\cap
V_{\mu^*}(y)|=\frac{1}{W_0(q)}\sum_{\substack{w_1\in
W_0w_{\la}W_{l}\\ w_2\in
W_n\s_n(w_{\mu^*})W_{l}}}|\cc_{w_1}(a')\cap\cc_{w_2}(b')|
\end{align}
where $a',b'$ are any chambers with $\delta(a',b')\in
W_0w_{\nu}W_n$.

Vertex regularity follows from (\ref{vreg}), for the value of
$|V_{\la}(x)\cap V_{\mu^*}(y)|$ is seen to only depend on
$\la,\mu$ and $\nu$. To see that strong vertex regularity holds,
we use Proposition~\ref{little} to see that
$$W_0w_{\la}W_{l}=\s_*(W_{\s_*^{-1}(0)}\s_*^{-1}(w_{\la})W_{\s_*^{-1}(l)})=\s_*(W_0w_{\la^*}W_{l^*})\,,$$
$$W_n\s_n(w_{\mu^*})W_{l}=\s_*(W_{n^*}(\s_*^{-1}\circ\s_n\circ\s_*)(w_{\mu})W_{l^*})=\s_*(W_{n^*}\s_{n^*}(w_{\mu})W_{l^*})\,,$$
and similarly $W_0w_{\nu}W_n=\s_*(W_0w_{\nu^*}W_{n^*})$. A
further application of Lemma~\ref{vertexuse} (with $\s=\s_*$)
implies that
\begin{align*}
|V_{\la}(x)\cap
V_{\mu^*}(y)|&=\frac{1}{W_0(q)}\sum_{\substack{w_1\in
W_0w_{\la^*}W_{l^*}\\ w_2\in
W_{n^*}\s_{n^*}(w_{\mu})W_{l^*}}}|\cc_{w_1}(a'')\cap\cc_{w_2}(b'')|
\end{align*}
where $a'',b''$ are any chambers with $\delta(a'',b'')\in W_0
w_{\nu^*}W_{n^*}$. Thus by comparison with (\ref{vreg}) we have
$$|V_{\la}(x)\cap V_{\mu^*}(y)|=|V_{\la^*}(x')\cap V_{\mu}(y')|,$$
where $x',y'\in V_P$ are any vertices with $y'\in V_{\nu^*}(x')$;
that is, strong vertex regularity holds.  \end{proof}

\begin{cor}\label{numbers2} There exist numbers $a_{\la,\mu;\nu}\in\bq^+$ such
that for $\la,\mu\in P^+$,
$$A_{\la}A_{\mu}=\sum_{\nu\in
P^+}a_{\la,\mu;\nu}A_{\nu}\qquad\textrm{and}\qquad\sum_{\nu\in
P^+}a_{\la,\mu;\nu}=1\,.$$ Moreover, $|\{\nu\in P^+\mid
a_{\la,\mu;\nu}\neq0\}|$ is finite for all $\la,\mu\in P^+$.
\end{cor}

\begin{proof} Let $v\in V_{\nu}(u)$ and set
\begin{align}
\label{finalchange}a_{\la,\mu;\nu}=\frac{N_{\nu}}{N_{\la}N_{\mu}}|V_{\la}(u)\cap
V_{\mu^*}(v)|\,,
\end{align}
which is independent of the particular pair $u,v$ by vertex
regularity. The numbers $a_{\la,\mu;\nu}$ are clearly nonnegative
and rational, and from (\ref{vertexop}) we have
\begin{align*}
(A_{\la}A_{\mu}f)(x)&=\sum_{\nu\in P^+}\left(\sum_{y\in
V_{\nu}(x)}\frac{|V_{\la}(x)\cap
V_{\mu^*}(y)|}{N_{\la}N_{\mu}}f(y)\right)\\
&=\sum_{\nu\in
P^+}a_{\la,\mu;\nu}\left(\frac{1}{N_{\nu}}\sum_{y\in
V_{\nu}(x)}f(y)\right)\\
&=\sum_{\nu\in P^+} a_{\la,\mu;\nu}(A_{\nu}f)(x)\,.
\end{align*}
When $f=1_{V_P}:V_P\to\{1\}$ we see that $\sum a_{\la,\mu;\nu}=1$.

We now show that only finitely many of the $a_{\la,\mu;\nu}$'s
are nonzero for each fixed pair $\la,\mu\in P^+$. Fix $x\in V_P$
and observe that $a_{\la,\mu;\nu}\neq0$ if and only if
$V_{\la}(x)\cap V_{\mu^*}(y)\neq \emptyset$ for each $y\in
V_{\nu}(x)$. Applying $(N_{\la}A_{\la})(N_{\mu}A_{\mu})$ to the
constant function $1_{V_P}:V_P\to\{1\}$, we obtain
$$
\sum_{y\in V_P}|V_{\la}(x)\cap V_{\mu^*}(y)|=N_{\la}N_{\mu},
$$
and hence $V_{\la}(x)\cap V_{\mu^*}(y)\neq\emptyset$ for only
finitely many $y\in V_P$. \end{proof}

\begin{defn}\label{algebraa} Let $\sca$ be the linear span of $\{A_{\la}\mid\la\in
P^+\}$ over $\bc$. The previous corollary shows that $\sca$ is an
associative algebra.
\end{defn}

We refer to the numbers $a_{\la,\mu;\nu}$ in Corollary
\ref{numbers2} as the \textit{structure constants} of the algebra
$\sca$.

\begin{thm}\label{maina} The algebra $\sca$ is commutative.
\end{thm}

\begin{proof} We need to show that
$a_{\la,\mu;\nu}=a_{\mu,\la,\nu}$ for all $\la,\mu,\nu\in P^+$.
Fixing any pair $u,v$ in $V_P$ with $v\in V_{\nu}(u)$, strong
vertex regularity implies that
\begin{align*}
a_{\la,\mu;\nu}&=\frac{N_{\nu}}{N_{\la}N_{\mu}}|V_{\la}(u)\cap
V_{\mu^*}(v)|=\frac{N_{\nu}}{N_{\la}N_{\mu}}|V_{\la^*}(v)\cap
V_{\mu}(u)|=a_{\mu,\la;\nu}
\end{align*}
completing the proof.  \end{proof}

We note that a similar calculation using Theorem~\ref{N}
(specifically the fact that $N_{\la}=N_{\la^*}$) shows that
$a_{\la,\mu;\nu}=a_{\la^*,\mu^*;\nu^*}$ for all $\la,\mu,\nu\in
P^+$.

\begin{rem}\label{AS} Let $X$ be a set and let $K$ be a partition of
$X\times X$ such that $\emptyset\notin K$ and $\{(x,x)\mid x\in
X\}\in K$. For $k\in K$, define $k^*=\{(y,x)\mid (x,y)\in k\}$,
and for each $x\in X$ and $k\in K$ define $k(x)=\{y\in
X\mid(x,y)\in k\}$. Recall \cite{zieschang} that an
\textit{association scheme} is a pair ($X,K$) as above such that
(i) $k\in K$ implies that $k^*\in K$, and (ii) for each $k,l,m\in
K$ there exists a cardinal number $e_{k,l;m}$ such that
$$(x,y)\in m\quad\textrm{implies that}\quad|k(x)\cap
l^*(y)|=e_{k,l;m}\,.$$

Let $X=V_P$, and for each $\la\in P^+$ let $\la'=\{(x,y)\mid y\in
V_{\la}(x)\}$. The set $L=\{\la'\mid\la\in P^+\}$ forms a
partition of $V_P\times V_P$, and $\la'(x)=V_{\la}(x)$ for $x\in
V_P$.

By vertex regularity it follows that the pair $(V_P,L)$ forms an
association scheme, and the cardinal numbers $e_{\la',\mu';\nu'}$
are simply $N_{\la}N_{\mu}N_{\nu}^{-1}a_{\la,\mu;\nu}$. By strong
vertex regularity this association scheme also satisfies the
condition $e_{\la',\mu';\nu'}=e_{\mu',\la';\nu'}$ for all
$\la,\mu,\nu\in P^+$ (see \cite[p.1, footnote]{zieschang}).

Note that the algebra $\sca$ is essentially the
\textit{Bose-Mesner algebra} of the association scheme $(V_P,L)$
(see \cite[Chapter 2]{rosemary}). With reference to
Remark~\ref{IDRG}, the above construction generalises the familiar
construction of association schemes from infinite distance
regular graphs (see \cite[\ts1.4.4]{rosemary} for the case of
\textit{finite} distance regular graphs).
\end{rem}

Recall the definition of the numbers $b_{w_1,w_2;w_3}$ given in
Corollary \ref{numbers}.

\begin{prop}\label{motive} Let $\tau(\la)=l$ and $\tau(\nu)=n$. Suppose that $y\in
V_{\nu}(x)$ and $V_{\la}(x)\cap V_{\mu^*}(y)\neq\emptyset$. Then
\begin{align*}
a_{\la,\mu;\nu}&=\frac{W_{0\la}(q)W_{0\mu}(q)}{W_{0\nu}(q)W_0^2(q)q_{w_{\la}}q_{w_{\mu}}}\sum_{\substack{w_1\in
W_0w_{\la}W_{l}\\ w_2\in
W_{l}\s_{l}(w_{\mu})W_n}}q_{w_1}q_{w_2}b_{w_1,w_2;w_{\nu}}
\end{align*}
\end{prop}

\begin{proof} By Lemma~\ref{track}(ii) we have
$\s_n=\s_l\circ\s_m$. Thus by Proposition~\ref{little}(iv) we have
$W_n\s_n(w_{\mu^*})W_{l}=(W_{l}\s_{l}(w_{\mu})W_n)^{-1}$, and so
by (\ref{vreg}) we see that
\begin{align}
\label{sss}|V_{\la}(x)\cap
V_{\mu^*}(y)|=\frac{1}{W_0(q)}\sum_{\substack{w_1\in
W_0w_{\la}W_{l}\\ w_2\in
W_{l}\s_{l}(w_{\mu})W_n}}|\cc_{w_1^{\vphantom{-1}}}(a)\cap\cc_{w_2^{-1}}(b)|
\end{align}
whenever $\delta(a,b)\in W_0w_{\nu}W_n$.

By Proposition~\ref{cor1} (and the proof thereof) we have
$$|\cc_{w_1^{\vphantom{-1}}}(a)\cap\cc_{w_2^{-1}}(b)|=q_{w_1}q_{w_2}(B_{w_1}B_{w_2}\delta_b)(a)\,,$$
and the result now follows from (\ref{sss}) by using
Theorem~\ref{N} and the definitions of $a_{\la,\mu;\nu}$ and
$b_{w_1,w_2;w_3}$, by choosing $b\in~\cc_{w_{\nu}}(a)$.
\end{proof}

\section{Affine Hecke Algebras and Macdonald Spherical
Functions}\label{section6}

In Section~\ref{6.3} we make an important connection between the algebra~$\sca$
and affine Hecke algebras. In particular, in Theorem~\ref{i} we show that $\sca$
is isomorphic to $Z(\tilde{\sch})$, the centre of an appropriately parametrised
affine Hecke algebra~$\tilde{\sch}$.

In Sections~\ref{6.1} and~\ref{6.2} we give an outline of some known results
regarding affine Hecke algebras. The main references for this material are \cite{m} and \cite{ram}.
Note that in~\cite{ram} there is only one parameter~$q$, although the results
there go through without any serious difficulty in the more general case of
multiple parameters $\{q_s\}_{s\in S}$. Note also that in \cite{ram} $Q=Q(R)$ and $P=P(R)$, whereas for us
$Q=Q(\chR)$ and $P=P(\chR)$.

\subsection{Affine Hecke Algebras}\label{6.1}

Let $\{q_{s}\}_{s\in S}$ be a set of positive real numbers with
$q_{s_i}=q_{s_j}$ whenever $s_i$ and $s_j$ are conjugate in
$\tilde{W}$. The \textit{affine Hecke algebra} $\tilde{\sch}$
with parameters $\{q_s\}_{s\in S}$ is the algebra over $\bk$ with
presentation given by the generators $T_w$, $w\in\tilde{W}$, and
relations
\begin{align}
\label{f1}T_{w_1}T_{w_2}&=T_{w_1w_2}&&\textrm{if
$\ell(w_1w_2)=\ell(w_1)+\ell(w_2)$}\,,\\
\label{f2}T_wT_s&=\frac{1}{q_s}T_{ws}+\left(1-\frac{1}{q_s}\right)T_{w}&&\textrm{if
$\ell(ws)<\ell(w)$ and $s\in S$\,.}
\end{align}

By~(\ref{f1}), $T_1T_w=T_wT_1=T_w$ for all $w\in\tilde{W}$, and
hence $T_1=I$ since $\{T_w\}_{w\in \tilde{W}}$
generates~$\tilde{\sch}$. Then~(\ref{f2}) implies that each $T_s$,
$s\in S$, is invertible, and from~(\ref{f1}) we see that each
$T_g$, $g\in G$, is invertible, with inverse $T_{g^{-1}}$ (recall
the definition of $G$ from Section~\ref{tilde{W}}). Since each
$w\in\tilde{W}$ can be written as $w=w'g$ for $w'\in W$ and $g\in
G$ it follows that each $T_w$, $w\in\tilde{W}$, is invertible.

\begin{rem} (i) In \cite{m} the numbers $\{q_{s}\}_{s\in S}$ are taken as
positive real variables. Our choice to fix the numbers
$\{q_s\}_{s\in S}$ does not change the algebraic structure of
$\tilde{\sch}$ in any serious way (for our purposes, at least).

(ii) The condition that $q_{s_i}=q_{s_j}$ whenever
$s_i=ws_jw^{-1}$ for some $w\in \tilde{W}$ is equivalent to the
condition that $q_{s_i}=q_{s_j}$ whenever $s_i=us_{\s(j)}u^{-1}$
for some $\s\in\Auttr(D)$ and $u\in W$. This condition is quite
restrictive, and it is easy to see that we obtain the parameter
systems given in the Appendix. Thus connections with our earlier
results on the algebra $\sca$ will become apparent
when we take the numbers $\{q_s\}_{s\in S}$ to be the parameters
of a locally finite regular affine building.
\end{rem}

\begin{defn}\label{defnum} (i) We write $q_w=q_{s_{i_1}}\cdots q_{s_{i_m}}$ if $s_{i_1}\cdots s_{i_m}$ is a reduced expression for
$w\in W$. This is easily seen to be independent of the particular
reduced expression (see \cite[IV, \ts1, No.5,
Proposition~5]{bourbaki}). Each $\tilde{w}\in \tilde{W}$ can be
written uniquely as $\tilde{w}=wg$ for $w\in W$ and $g\in G$, and
we define $q_{\tilde{w}}=q_w$. In particular $q_g=1$ for all $g\in
G$. Furthermore, if $s=s_i$ we write $q_s=q_i$.

(ii) To conveniently state later results we make the following
definitions. Let $R_1=\{\alpha\in R\mid 2\alpha\notin R\}$,
$R_2=\{\alpha\in R\mid \frac{1}{2}{\alpha}\notin R\}$ and
$R_3=R_1\cap R_2$ (so $R_1=R_2=R_3=R$ if $R$ is reduced). For
$\alpha\in R_2$, write $q_{\alpha}=q_i$ if $\alpha\in
W_0\alpha_i$ (note that if $\alpha\in W_0\alpha_i$ then
necessarily $\alpha\in R_2$). It follows easily from
Corollary~\ref{conjno} that this definition is unambiguous.

Note that $R$ is the disjoint union of $R_3$, $R_1\backslash R_3$
and $R_2\backslash R_3$, and define set of numbers
$\{\tau_{\alpha}\}_{\alpha\in R}$ by
\begin{align*}
\tau_{\alpha}=\begin{cases} q_{\alpha}&\textrm{if $\alpha\in
R_3$}\\ q_0&\textrm{if $\alpha\in R_1\backslash R_3$}\\
q_{\alpha}q_{0}^{-1}&\textrm{if $\alpha\in R_2\backslash R_3$},
\end{cases}
\end{align*}
where $q_0=q_{s_0}$ (with $s_0=s_{\tilde{\alpha};1}$ and
$\tilde{\alpha}$ is as in (\ref{highestroot})). It is convenient
to also define $\tau_{\alpha}=1$ if $\alpha\notin R$. The reader
only interested in the reduced case can simply read
$\tau_{\alpha}$ as $q_{\alpha}$. Note that
$\tau_{w\alpha}=\tau_{\alpha}$ for all $\alpha\in R$ and $w\in
W_0$.
\end{defn}

\begin{rem}\label{converttof} We have chosen a slight distortion of the usual
definition of the algebra $\tilde{\sch}$. This choice has been
made so as to make the connection between the algebras $\sca$ and
$\tilde{\sch}$ more transparent, as the reader will shortly see.
To allow the reader to convert between our notation and that in
\cite{m}, we provide the following instructions. With reference
to our presentation for $\tilde{\sch}$ given above, let
$\tau_i=\sqrt{q_i}$ and $T_w'=\sqrt{q_w}\,T_w$ (these $\tau$'s
are unrelated to those in Definition~\ref{defnum}(ii)). Our
presentation then transforms into that given in \cite[4.1.2]{m}
(with the $T$'s there replaced by $T'$'s). This transformation
also makes it clear why the $\sqrt{q_w}\,$'s appear in the
following discussion.
\end{rem}

If $\la\in P^+$ let $x^{\la}=\sqrt{q_{t_{\la}}}\,T_{t_{\la}}$,
and if $\la=\mu-\nu$ with $\mu,\nu\in P^+$ let
$x^{\la}=x^{\mu}(x^{\nu})^{-1}$. This is well defined by
\cite[p.40]{m}, and for all $\la,\mu\in P$ we have
$x^{\la}x^{\mu}=x^{\la+\mu}=x^{\mu}x^{\la}$.

We write $\bk[P]$ for the $\bk$-span of $\{x^{\la}\mid \la\in
P\}$. The group $W_0$ acts on $\bk[P]$ by linearly extending the
action $wx^{\la}=x^{w\la}$. We write $\bk[P]^{W_{0}}$ for the set
of elements of $\bk[P]$ that are invariant under the action of
$W_0$. By Corollary~\ref{cor:center}, the centre $Z(\tilde{\sch})$
of $\tilde{\sch}$ is $\bk[P]^{W_0}$.

Let $\sch$ be the subalgebra of $\tilde{\sch}$ generated by
$\{T_{s}\mid s\in S\}$. The following relates the algebra $\sch$
to the algebra $\scb$ of chamber set averaging operators on an
irreducible affine building.

\begin{prop}\label{connectb} Suppose a building $\scx$ of type $R$ exists with parameters
$\{q_{s}\}_{s\in S}$. Then $\sch\cong\scb$.
\end{prop}

\begin{proof} This follows in the same way as Theorem~\ref{usebalg}.  \end{proof}

We make the following parallel definition to (\ref{ooo}). Recall
the definition of Poincar\'{e} polynomials from
Definition~\ref{defn:poincarep}. For each $i\in I$, let
\begin{align}
\label{oo}\mathds{1}_i=\frac{1}{W_i(q)}\sum_{w\in W_i}q_wT_w,
\end{align}
where $W_i=W_{I\backslash\{i\}}$ (as before). Thus $\ooo_i$ is an
element of $\sch$. As a word of warning, we have used the same
notation as in (\ref{ooo}) where we defined the analogous element
in $\scb$. There should be no confusion caused by this decision.

The following lemma follows in exactly the same way as
Lemma~\ref{early1}.

\begin{lem}\label{magic} $\mathds{1}_iT_w=T_w\mathds{1}_i=\mathds{1}_i$ for all
$w\in W_i$ and $i\in I$. Furthermore $\ooo_i^2=\ooo_i$.
\end{lem}

\subsection{The Macdonald Spherical Functions}\label{6.2}

The following relations are of fundamental significance.

\goodbreak\begin{thm}\label{bb} Let $\la\in P$ and $i\in I_0$.
\begin{itemize}
\item[(i)] If $(R,i)\neq(BC_n,n)$ for any $n\geq1$, then
$$x^{\la}T_{s_i}-T_{s_i}x^{s_i\la}=(1-q_i^{-1})\frac{x^{\la}-x^{s_i\la}}{1-x^{-\alpha_i^{\vee}}}.$$
\item[(ii)] If $R=BC_n$ for some $n\geq1$ and $i=n$, then
$$x^{\la}T_{s_n}-T_{s_n}x^{s_n\la}=\left[1-q_n^{-1}+q_n^{-1/2}\big(q_0^{1/2}-q_0^{-1/2}\big)x^{-(2\alpha_n)^{\vee}}\right]\frac{x^{\la}-x^{s_n\la}}{1-x^{-2(2\alpha_n)^{\vee}}}.$$
\end{itemize}
\end{thm}

\begin{proof} This follows from \cite[(4.2.4)]{m} (see
Remark~\ref{converttof}), taking into account \cite[(1.4.3) and
(2.1.6)]{m} in case(ii).
\end{proof}

We note that the fractions appearing in Theorem~\ref{bb} are in
fact finite linear combinations of the $x^{\mu}$'s
\cite[(4.2.5)]{m}. We refer to the relations in Theorem~\ref{bb}
as the \textit{Bernstein relations},
for they are a crucial ingredient in the so-called
\textit{Bernstein presentation} of the Hecke algebra.

\begin{cor}\label{cor:center} The centre $Z(\tilde{\sch})$ of $\tilde{\sch}$ is
$\bc[P]^{W_0}$.
\end{cor}

\begin{proof} This well known fact can be proved using the
Bernstein relations, exactly as in \cite[(4.2.10)]{m}.
\end{proof}

For each $\la\in P^+$, define an element
$P_{\la}(x)\in\bk[P]^{W_0}$ by
\begin{align}
\label{function}P_{\la}(x)&=\frac{q_{t_{\la}}^{-1/2}}{W_0(q)}\sum_{w\in
W_0}w\left(x^{\la}\prod_{\alpha\in R
^+}\frac{\tau_{\alpha\vphantom{/2}}^{\vphantom{-1}}\tau_{\alpha/2}^{1/2}\,x^{\cha}-1}{\tau_{\alpha/2}^{1/2}\,x^{\cha}-1}\right).
\end{align}
We call the elements $P_{\la}(x)$ the \textit{Macdonald spherical
functions} of $\tilde{\sch}$.

\begin{rem}\label{macnorm} (i) We have chosen a slightly different normalisation of the
Macdonald spherical function from that in \cite{m}. Our formula
uses the normalisation of \cite[Theorem~4.1.2]{macsph}.

(ii) Notice that the formula simplifies in the reduced case
(namely, $\tau_{\alpha/2}=1$).

(iii) It is not immediately clear that $P_{\la}(x)$ as defined in
(\ref{function}) is in $\bk[P]^{W_0}$, although this is a
consequence of \cite[VI, \ts3, No.3, Proposition~2]{bourbaki}.
\end{rem}

The proof of Theorem~\ref{sig} below follows
\cite[Theorem~2.9]{ram} very closely. 

\begin{thm}\label{sig}\cite[Theorem~2.9]{ram}\textbf{.} For $\la\in P^+$ we have
$q_{t_{\la}}^{1/2}P_{\la}(x)\oo=\oo x^{\la}\oo$.
\end{thm}

\begin{proof} By the Satake isomorphism (see \cite[Theorem~2.4]{ram} and
\cite[5.2]{knop} for example) there exists some
$P_{\la}'(x)\in\bc[P]^{W_0}$ such that $P'_{\la}(x)\oo=\oo
x^{\la}\oo$. If $i\in I_0$ and $(R,i)\neq(BC_n,n)$, then by Theorem~\ref{bb}(i)
(and using Lemma~\ref{magic}) we have
\begin{align}\begin{aligned}\label{eee1}
(1+q_iT_{s_i})x^{\la}\ooo_0&=x^{\la}\ooo_0+q_ix^{s_i\la}T_{s_i}\ooo_0+(q_i-1)\frac{x^{\la}-x^{s_i\la}}{1-x^{-\alpha_i^{\vee}}}\ooo_0\\
&=\frac{q_ix^{\la}-x^{\la-\alpha_i^{\vee}}-q_ix^{s_i\la-\alpha_i^{\vee}}+x^{s_i\la}}{1-x^{-\alpha_i^{\vee}}}\ooo_0\\
&=\bigg(\frac{q_ix^{\alpha_i^{\vee}}-1}{x^{\alpha_i^{\vee}}-1}x^{\la}+\frac{q_ix^{-\alpha_i^{\vee}}-1}{x^{-\alpha_i^{\vee}}-1}x^{s_i\la}\bigg)\ooo_0\\
&=(1+s_i)\frac{q_ix^{\alpha_i^{\vee}}-1}{x^{\alpha_i^{\vee}}-1}x^{\la}\ooo_0.
\end{aligned}\end{align} A similar calculation, using Theorem~\ref{bb}(ii),
shows that if $i\in I_0$ and $(R,i)=(BC_n,n)$, then
\begin{align}\label{eee2}
(1+q_nT_{s_n})x^{\la}\ooo_0=(1+s_n)\frac{\big(\sqrt{q_0q_n}\,x^{(2\alpha_n)^{\vee}}-1\big)\big(\sqrt{q_n/q_0}\,x^{(2\alpha_n)^{\vee}}+1\big)}{x^{2(2\alpha_n)^{\vee}}-1}x^{\la}\ooo_0.
\end{align}

It will be convenient to write (\ref{eee1}) and (\ref{eee2}) as
one equation, as follows. In the reduced case, let
$\beta_i=\alpha_i$ for all $i\in I_0$, and in the non-reduced
case (so $R=BC_n$ for some $n\geq1$) let $\beta_i=\alpha_i$ for
$1\leq i\leq n-1$ and let $\beta_n=2\alpha_n$. For $\alpha\in R$
and $i\in I_0$, write
$$a_i(x^{\alpha^{\vee}})=\frac{(\tau_{\b_i\vphantom{/2}}^{\vphantom{-1}}\tau_{\b_i/2}^{1/2}\,x^{\alpha^{\vee}}-1)(\tau_{\b_i/2}^{1/2}\,x^{\alpha^{\vee}}+1)}{x^{2\alpha^{\vee}}-1},$$
and so in all cases
\begin{align}
\label{getstarted}(1+q_{i}T_{s_i})x^{\la}\oo=(1+s_i)a_i(x^{\beta_i^{\vee}})x^{\la}\oo\,.
\end{align}
By induction we see that (writing $T_i$ for $T_{s_i}$)
\begin{align}
\label{starmac}\left[\prod_{k=1}^m(1+q_{i_k}T_{i_k})\right]x^{\la}\oo=\left[\prod_{k=1}^{m}(1+s_{i_k})a_{i_k}(x^{\beta_{i_k}^{\vee}})\right]x^{\la}\oo,
\end{align}
where we write $\prod_{k=1}^{m}x_k$ for the ordered product
$x_1\cdots x_m$. Therefore $\oo x^{\la}\oo$ can be written as
$fx^{\la}\oo$, where $f$ is independent of $\la$ and is a finite
linear combination of terms of the form
$$(1+s_{i_1})a_{i_1}(x^{\beta_{i_1}^{\vee}})\cdots(1+s_{i_m})a_{i_m}(x^{\beta_{i_m}^{\vee}})$$
where $i_1,\ldots,i_m\in I_0$.

Thus we have
$$P_{\la}'(x)=\sum_{w\in W_0}w\left(b_{w}(x)x^{\la}\right)$$
where each $b_w(x)$ is a linear combination of products of terms
$a_i(x^{\beta_i^{\vee}})$ and is independent of $\la\in P^+$. It
is easily seen that this expression is unique, and since
$P_{\la}'(x)\in\bc[P]^{W_0}$ it follows that $b_{w}(x)=b_{w'}(x)$
for all $w,w'\in W_0$, and we write $b(x)$ for this common value.
Thus
$$P_{\la}'(x)=\sum_{w\in W_0}w\left(b(x)x^{\la}\right)=\sum_{w\in
W_0}w\left(x^{w_0\la}w_0b(x)\right)$$
where $w_0$ is the longest element of $W_0$.

We now compute the coefficient of $x^{w_{0}\la}$ in the above
expression. Since this coefficient is independent of $\la\in P^+$
we may assume that $\lan\la,\alpha_i\ran>0$ for all $i\in I_0$
and so $w\la\neq w_0\la$ for all $w\in W_0\backslash\{w_0\}$.

If $w_0=s_{i_1}\cdots s_{i_m}$ is a reduced expression, then
\begin{align*}
\oo&=\frac{1}{W_0(q)}\bigg((1+q_{i_1}T_{i_1})\cdots(1+q_{i_m}T_{i_m})\\&\quad+\textrm{terms
$(1+q_{j_1}T_{j_1})\cdots(1+q_{j_l}T_{j_l})$ with $j_k\in I_0$ and
$l<m$}\bigg)\,.
\end{align*}
Thus, by (\ref{starmac})
\begin{align*}
\oo
x^{\la}\oo&=\frac{1}{W_0(q)}\Bigg[\left(\prod_{k=1}^{m}s_{i_{k}}a_{i_k}(x^{\beta_{i_k}^{\vee}})\right)x^{\la}\oo\\
&\quad+\textrm{terms
}\left(\prod_{k=1}^{l}s_{j_{k}}a_{j_k}(x^{\beta_{j_k}^{\vee}})\right)x^{\la}\oo\textrm{
with $j_k\in I_0$ and $l<m$}\Bigg].
\end{align*}

Thus the coefficient of $x^{w_{0}\la}$ is
\begin{align*}
w_0b(x)&=\frac{1}{W_0(q)}s_{i_1}a_{i_1}(x^{\beta_{i_1}^{\vee}})\cdots
s_{i_m}a_{i_m}(x^{\beta_{i_m}^{\vee}})\\
&=\frac{1}{W_0(q)}w_0\prod_{\beta\in
R_1^+}\frac{(\tau_{\beta\vphantom{/2}}^{\vphantom{-1}}\tau_{\beta/2}^{1/2}\,x^{\beta^{\vee}}-1)(\tau_{\beta/2}^{1/2}\,x^{\beta^{\vee}}+1)}{x^{2\beta^{\vee}}-1}\,,
\end{align*}
where we have used the fact that
$$\{\chb_{i_m},s_{i_m}\chb_{i_{m-1}},\cdots,s_{i_m}s_{i_{m-1}}\cdots
s_{i_2}\chb_{i_1}\}=(R_1^+)^{\vee}\,,$$ (see \cite[(2.2.9)]{m})
and the fact that $\tau_{w\alpha}=\tau_{\alpha}$ for all $w\in
W_0$ and $\alpha\in R$. The result now follows after an
elementary manipulation.  \end{proof}

Since $x^{\la}=q_{t_{\la}}^{1/2}T_{t_{\la}}$ by definition, we have the following.

\begin{cor}\label{sat} For $\la\in P^+$ we have
$$\oo T_{t_{\la}}\oo=P_{\la}(x)\oo\,.$$
\end{cor}

We write $Q^+$ for the $\bn$-span of $\{\cha\mid\alpha\in R^+\}$.
Define a partial order $\preceq$ on $P$ by $ \mu\preceq\la$ if
and only if $\la-\mu\in Q^+$.

\begin{thm}\label{macrec} $\{P_{\la}(x)\mid\la\in\ P^+\}$ is a basis
of $\bc[P]^{W_0}$. Furthermore, the Macdonald spherical functions
satisfy
$${P}_{\la}(x){P}_{\mu}(x)=\sum_{\nu\preceq\la+\mu}c_{\la,\mu;\nu}{P}_{\nu}(x)$$
for some numbers $c_{\la,\mu;\nu}$, with $c_{\la,\mu;\la+\mu}>0$.
\end{thm}

\begin{proof} This is a simple application of the triangularity
condition of the Macdonald spherical functions, see
\cite[\ts10]{loeth}.  \end{proof}

\subsection{Connecting $\sca$ and $Z(\tilde{\sch})$}\label{6.3}
We can now see how to relate the vertex set
averaging operators $A_{\la}$ from Section \ref{section5} to the
algebra elements $P_{\la}(x)$. Let us recall (and make) some
definitions. For $\la,\mu,\nu\in P^+$ and $w_1,w_2,w_3\in W$,
define numbers $a_{\la,\mu;\nu}$, $b_{w_1,w_2;w_3}$,
$c_{\la,\mu;\nu}$ and $d_{w_1,w_2;w_3}$ by
\begin{align*}
A_{\la}A_{\mu}&=\sum_{\nu\in P^+}a_{\la,\mu;\nu}A_{\nu}&&B_{w_1}B_{w_2}=\sum_{w_3\in W}b_{w_1,w_2;w_3}B_{w_3}\\
P_{\la}(x)P_{\mu}(x)&=\sum_{\nu\in
P^+}c_{\la,\mu;\nu}P_{\nu}(x)&&\,\,T_{w_1}T_{w_2}=\sum_{w_3\in
W}d_{w_1,w_2;w_3}T_{w_3}\,.
\end{align*}
Thus the numbers are the structure constants of the algebras
$\sca$, $\scb$, $\bk[P]^{W_0}$ and $\sch$ with respect to the
bases $\{A_{\la}\mid\la\in P^+\}$, $\{B_w\mid w\in W\}$,
$\{P_{\la}(x)\mid\la\in P^+\}$ and $\{T_w\mid w\in W\}$
respectively.

Note that by Proposition~\ref{connectb} we have
$b_{w_1,w_2;w_3}=d_{w_1,w_2;w_3}$ whenever a building with
parameter system $\{q_{s}\}_{s\in S}$ exists. We stress that
$d_{w_1,w_2;w_3}$ is a more general object, for it makes sense
for a much more general set of $q_s$'s.

Recall the definition of $w_{\la}$ from Section~\ref{wla}, and
recall the definition of $W_{0\la}$ from (\ref{earlier?}). We
give the following lemma linking double cosets in $W$ with double
cosets in $\tilde{W}$.

\begin{lem}\label{form} Let $\la\in P^+$ and $i\in I_{P}$. Suppose
that $\tau(\la)=l$, and write $j=\s_i(l)$ (so
$\s_j=\s_i\circ\s_{l}$). Then
$$W_i\s_i(t_{\la}')W_{j}=g_iW_0t_{\la}W_0g_{j}^{-1}\,,$$
where the elements $g_i$ are defined in (\ref{G}).
\end{lem}

\begin{proof} By Proposition~\ref{coll}, $g_j=g_ig_{l}$ and
$t_{\la}=t_{\la}'g_{l}$, and by (\ref{pres}),
$\s_k(w)=g_kwg_{k}^{-1}$ for all $w\in W$ and $k\in I_P$. Thus
\begin{align*}
W_i\s_i(t_{\la}')W_j&=(g_iW_0g_i^{-1})(g_it_{\la}g_{l}^{-1}g_{i}^{-1})(g_jW_0g_j^{-1})=g_iW_0t_{\la}W_0g_j^{-1}\,.\qedhere
\end{align*}
 \end{proof}

\begin{lem}\label{importantalso}\cite[Lemma~2.7]{ram}\textbf{.} Let $\la\in P^+$. Then
$$\sum_{w\in
W_0t_{\la}W_0}q_wT_w=\frac{W_0^2(q)}{W_{0\la}(q)}q_{w_{\la}}\oo
T_{t_{\la}}\oo\,.$$
\end{lem}

\begin{proof} This can be deduced from Theorem~\ref{early3}, or
see the proof in \cite{ram}.  \end{proof}

The following important Theorem will be used (along with
Proposition~\ref{motive}) to prove that $\sca\cong
Z(\tilde{\sch})$.

\begin{thm}\label{bigone} Let $\la,\mu,\nu\in P^+$ and write $\tau(\la)=l$, $\tau(\mu)=m$ and $\tau(\nu)=n$.
Then if $c_{\la,\mu;\nu}\neq 0$ we have
\begin{align*}
c_{\la,\mu;\nu}&=\frac{W_{0\la}(q)W_{0\mu}(q)}{W_{0\nu}(q)W_0^2(q)q_{w_{\la}}q_{w_{\mu}}}\sum_{\substack{w_1\in
W_0w_{\la}W_{l}\\ w_2\in
W_{l}\s_{l}(w_{\mu})W_n}}q_{w_1}q_{w_2}d_{w_1,w_2;w_{\nu}}.
\end{align*}
\end{thm}

\begin{proof} To abbreviate notation we write
$P_{\la}=P_{\la}(x)$. First observe that by Theorem~\ref{macrec} we
have $c_{\la,\mu;\nu}=0$ unless $\nu\preceq\la+\mu$. In
particular we have $c_{\la,\mu;\nu}=0$ when
$\tau(\nu)\neq\tau(\la+\mu)$. It follows that
$\s_n=\s_{l}\circ\s_{m}$, and so $g_n=g_{l}g_m$ (see
Proposition~\ref{coll}). We will use this fact later.

By Corollary \ref{sat} and Lemma~\ref{importantalso}, for any
$\la\in P^+$ we have
$$P_{\la}\oo=\oo
T_{t_{\la}}\oo=\frac{W_{0\la}(q)}{W_0^2(q)q_{w_{\la}}}\sum_{w\in
W_0t_{\la}W_0}q_wT_w\,,$$ and so if $i\in I_P$, $\tau(\la)=l$ and
$j=\s_i(l)$ we have (see Lemma~\ref{form})
\begin{align}
\label{astar}T_{g_i}P_{\la}\oo
T_{g_j^{-1}}=\frac{W_{0\la}(q)}{W_0^2(q)q_{w_{\la}}}\sum_{w\in
W_i\s_i(t_{\la}')W_j}q_wT_w\,.
\end{align}
We can replace the $t_{\la}'$ by $w_{\la}$ in the above because
$W_i\s_i(t_{\la}')W_j=W_i\s_i(w_{\la})W_j$ by
Proposition~\ref{early2}(i) and the fact that $\s_i(W_l)=W_j$.

Using the fact that $g_n=g_{l}g_m$ if $c_{\la,\mu;\nu}\neq 0$ we
have, by (\ref{astar})
\begin{align*}
P_{\la}\oo P_{\mu}\oo T_{g_n^{-1}}&=(P_{\la}\oo
T_{g_{l}^{-1}})(T_{g_{l}}P_{\mu}\oo
T_{g_n^{-1}})\\
&=\frac{W_{0\la}(q)W_{0\mu}(q)}{W_0^4(q)q_{w_{\la}}q_{w_{\mu}}}\sum_{\substack{w_1\in
W_0w_{\la}W_{l}\\ w_2\in W_{l}
\s_{l}(w_{\mu})W_{n}}}q_{w_1}q_{w_2}T_{w_1}T_{w_2}\\
&=\frac{W_{0\la}(q)W_{0\mu}(q)}{W_0^4(q)q_{w_{\la}}q_{w_{\mu}}}\sum_{w_3\in
W} \Bigg(\sum_{\substack{w_1\in W_0w_{\la}W_{l}\\ w_2\in W_{l}
\s_{l}(w_{\mu})W_{n}}}q_{w_1}q_{w_2}d_{w_1,w_2;w_3}T_{w_3}\Bigg)\,.
\end{align*}
So the coefficient of $T_{w_{\nu}}$ in the expansion of
$P_{\la}\oo P_{\mu}\oo T_{g_n^{-1}}$ in terms of the $T_w$'s is
\begin{align}
\label{countoneway}\frac{W_{0\la}(q)W_{0\mu}(q)}{W_0^4(q)q_{w_{\la}}q_{w_{\mu}}}\sum_{\substack{w_1\in
W_0w_{\la}W_{l}\\ w_2\in W_{l}
\s_{l}(w_{\mu})W_{n}}}q_{w_1}q_{w_2}d_{w_1,w_2;w_{\nu}}\,.
\end{align}

On the other hand, by Theorem~\ref{macrec} we have
\begin{align*}
P_{\la}\oo P_{\mu}\oo
T_{g_n^{-1}}&=\sum_{\eta\preceq\la+\mu}c_{\la,\mu;\eta}P_{\eta}\oo
T_{g_n^{-1}}\\
&=\sum_{\eta\preceq\la+\mu}\Bigg(\frac{W_{0\eta}(q)}{W_0^2(q)q_{w_{\eta}}}c_{\la,\mu;\eta}\sum_{w\in
W_0w_{\eta}W_n}q_wT_w\Bigg)\,.
\end{align*}
Since the double cosets $W_0w_{\eta}W_n$ are disjoint over
$\{\eta\in P^+\mid\eta\preceq\la+\mu\}$ we see that the
coefficient of $T_{w_{\nu}}$ is
\begin{align}
\label{andanother}\frac{W_{0\nu}(q)}{W_0^2(q)}c_{\la,\mu;\nu}\,.
\end{align}
The theorem now follows by equating (\ref{countoneway}) and
(\ref{andanother}).  \end{proof}

\begin{cor}\label{eqq} Suppose that an irreducible affine building exists with parameter
system $\{q_s\}_{s\in S}$. Then for all $\la,\mu,\nu\in P^+$ we
have $a_{\la,\mu;\nu}=c_{\la,\mu;\nu}$.
\end{cor}
\begin{proof} This follows from Theorem~\ref{bigone} and Propositions~\ref{motive} and~\ref{connectb}.
\end{proof}

\begin{thm}\label{i} Suppose that an irreducible affine building exists with parameters $\{q_s\}_{s\in S}$.
Then the map $P_{\la}(x)\mapsto A_{\la}$ determines an algebra isomorphism,
and so $\sca\cong Z(\tilde{\sch})$.
\end{thm}

\begin{proof} Since $\{P_{\la}(x)\mid\la\in P^+\}$ is a basis of $\bc[P]^{W_0}$ and
$\{A_{\la}\mid\la\in P^+\}$ is a basis of~$\sca$, there exists a
unique vector space isomorphism $\Phi:Z(\tilde{\sch})\to\sca$
with $\Phi(P_{\la})=A_{\la}$ for all $\la\in P^+$. Since
$a_{\la,\mu;\nu}=c_{\la,\mu;\nu}$ by Corollary \ref{eqq}, we see
that $\Phi$ is an algebra isomorphism.  \end{proof}

\begin{thm}\label{generators} The algebra $Z(\tilde{\sch})$ is generated by $\{P_{\la_i}(x)\mid i\in I_0\}$, and
so $\sca$ is generated by $\{A_{\la_i}\mid i\in I_0\}$.
\end{thm}

\begin{proof} First we define a less restrictive partial order on $P^+$ than
$\preceq$. For $\la,\mu\in P ^+$ we define $\mu<\la$ if and only
if $\la-\mu$ is an $\br^+$--linear combination of $(\chR)^+$ and
$\la\neq\mu$. Clearly if $\mu\prec\la$ then $\mu<\la$. Observe
also that $\la_i> 0$ for all $i\in I_0$ (see exercises 7 and 8 on
p.72 of \cite{h2}). Thus if $\la=\la'+\la_i$ for some $\la'\in
P^+$ and $i\in I_0$, we have $\la-\la'=\la_i>0$ and so $\la'<\la$.

Let $\mathcal{P}(\la)$ be the statement that $P_{\la}$ is a
polynomial in $P_{\la_1},\ldots,P_{\la_n}$ (and $P_0=1$). Suppose
that $\mathcal{P}(\la)$ fails for some $\la\in P^+$. Since
$\{\mu\in P^+\mid\mu\leq\la\}$ is finite (by the proof of
\cite[Lemma~13.2B]{h2}) we can pick $\la\in P^+$ minimal with
respect to $\leq$ such that $\mathcal{P}(\la)$ fails. There is an
$i$ such that $\la-\la_i=\la'$ is in $P^+$. Then $\la'<\la$ and
$P_{\la}=cP_{\la'}P_{\la_i}+$ a linear combination of $P_{\mu}$'s
where $\mu<\la$, $\mu\neq\la$. Then $\mathcal{P}(\la')$ holds, as
does $\mathcal{P}(\mu)$ for all these $\mu$'s. So
$\mathcal{P}(\la)$ holds, a contradiction.  \end{proof}

\section{A Positivity Result and Hypergroups}\label{section8}

Here we show that the structure constants $c_{\la,\mu;\nu}$ of the
algebra $\bc[P]^{W_0}$ are, up to positive normalisation factors,
polynomials with nonnegative integer coefficients in the variables
$\{q_{s}-1\mid s\in S\}$. This result has independently been
obtained by Schwer in~\cite{schwer}, where a formula for
$c_{\la,\mu;\nu}$ is given (in the case $q_s=q$ for all $s\in S$).

Thus if $q_s\geq1$ for all $s\in S$,
then $c_{\la,\mu;\nu}\geq0$ for all $\la,\mu,\nu\in P^+$. This
result was proved for root systems of type $A_n$ by Miller Malley
in \cite{mm}, where the numbers $c_{\la,\mu;\nu}$ are Hall
polynomials (up to positive normalisation factors). Note that it
is clear from (\ref{finalchange}) and Corollary \ref{eqq} that
$c_{\la,\mu;\nu}\geq0$ when there exists a building with
parameters $\{q_s\}_{s\in S}$.

In a recent series of papers \cite{rapoport}, \cite{haines},
\cite{tupan} and \cite{schwer} the numbers $a_{\la,\mu}$
appearing in $P_{\la}(x)=\sum_{\mu}a_{\la,\mu}m_{\mu}$ are
studied. Here $m_{\mu}$ is the monomial symmetric function
$\sum_{\gamma\in W_0\mu}x^{\gamma}$, where $W_0\mu$ is the orbit
$\{w\mu\mid w\in W_0\}$. We will provide a connection with the
results we prove here and the numbers $a_{\la,\mu}$ in
\cite[Theorem~6.2]{p2}. In particular, for $\la\in P^+$, let
$\Pi_{\la}\subset P$ denote the saturated set (see
\cite[VI,\ts1,Exercise~23]{bourbaki}) with highest coweight~$\la$.
If $\mu\notin\Pi_{\la}$ then $a_{\la,\mu}=0$, and for all
$\mu\in\Pi_{\la}$,
$$a_{\la,\mu}=\sqrt{\frac{N_{\nu-\mu}}{N_{\nu}}}c_{\la,\mu;\nu},$$
where $\nu$ is any dominant coweight with each
$\lan\nu,\alpha_i\ran$ `sufficiently large'.

The results of this section show how to construct a (commutative)
polynomial hypergroup, in the sense of \cite{bloom} (see also
\cite{hyper} where the $A_2$ case is discussed).

For each $w_1,w_2,w_3\in W$, let
$d'_{w_1,w_2;w_3}=q_{w_1}^{\vphantom{-1}}q_{w_2}^{\vphantom{-1}}q_{w_3}^{-1}d_{w_1,w_2;w_3}$.

\begin{lem}\label{neat} For all $w_1,w_2,w_3\in W$, $d_{w_1,w_2;w_3}'$ is a polynomial with nonnegative
integer coefficients in the variables $q_{s}-1$, $s\in S$.
\end{lem}

\begin{proof} We prove the result by induction on $\ell(w_2)$.
When $\ell(w_2)=1$, so $w_2=s$ for some $s\in S$, we have
\begin{align*}
d'_{w_1,s;w_3}=\begin{cases}1&\textrm{if $\ell(w_1s)=\ell(w_1)+1$
and $w_3=w_1s$\,,}\\
q_s&\textrm{if $\ell(w_1s)=\ell(w_1)-1$ and $w_3=w_1s$\,,}\\
q_s-1&\textrm{if $\ell(w_1s)=\ell(w_1)-1$ and $w_3=w_1$\,,}\\
0&\textrm{otherwise\,,}
\end{cases}
\end{align*}
proving the result in this case.

Suppose that $n\geq2$ and that the result is true for
$\ell(w_2)<n$. Then if $\ell(w_2)=n$, write $w_2=ws$ with
$\ell(w)=n-1$. Thus
\begin{align*}
T_{w_1}T_{w_2}&=(T_{w_1}T_w)T_s=\sum_{w'\in
W}d_{w_1,w;w'}T_{w'}T_s=\sum_{w_3\in W}\left(\sum_{w'\in
W}d_{w_1,w;w'}d_{w',s;w_3}\right)T_{w_3}\,,
\end{align*}
which implies that
$$d'_{w_1,w_2;w_3}=\sum_{w'\in W}d'_{w_1,w;w'}d'_{w',s;w_3}\,.$$
The result follows since $\ell(w)<n$ and $\ell(s)=1$.
\end{proof}

For each $\la,\mu,\nu\in P^+$, let
\begin{align}
\label{lastref}c'_{\la,\mu;\nu}=\frac{W_0(q)W_{0\nu}(q)}{W_{0\la}(q)W_{0\mu}(q)}\frac{q_{w_{\la}}q_{w_{\mu}}}{q_{w_{\nu}}}c_{\la,\mu;\nu}.
\end{align}

\begin{thm}\label{hyp} For all $\la,\mu,\nu\in P^+$, the structure constants $c_{\la,\mu;\nu}'$ are polynomials with nonnegative
integer coefficients in the variables $q_{s}-1$, $s\in S$.
\end{thm}

\begin{proof} We will use the same notation as in Theorem~\ref{bigone}, so
let $\tau(\la)=l$, $\tau(\mu)=m$ and $\tau(\nu)=n$. By
Theorem~\ref{bigone} we have
$$c_{\la,\mu;\nu}'=\frac{1}{W_0(q)}\sum_{\substack{w_1\in W_0w_{\la}W_l\\ w_2\in
W_l\s_l(w_{\mu})W_n}}d_{w_1,w_2;w_{\nu}}',$$ and so it
immediately follows from Lemma~\ref{neat} that
$W_0(q)c_{\la,\mu;\nu}'$ is a polynomial in the variables
$q_s-1$, $s\in S$, with nonnegative integer coefficients. The
result stated in the theorem is stronger than this, and so we
need to sharpen the methods used in the proof of
Theorem~\ref{bigone}.

We make the following observations. See Proposition~\ref{early2}
for proofs of similar facts (we use the notations of
Proposition~\ref{early2} here). Firstly, each $w_1\in
W_0w_{\la}W_l$ can be written uniquely as $w_1=u_1w_{\la}w_l$ for
some $u_1\in W_0^{\la}$ and $w_l\in W_l$, and moreover
$\ell(w_1)=\ell(u_1)+\ell(w_{\la})+\ell(w_l)$. Similarly, each
$w_2\in W_l\s_l(w_{\mu})W_n$ can be written uniquely as
$w_2=w_l'\s_l(w_{\mu})u_2$ for some $u_2\in W_n^{\mu}$ and
$w_l'\in W_l$, and moreover
$\ell(w_2)=\ell(w_l')+\ell(\s_l(w_{\mu}))+\ell(u_2)$.

Secondly, each $w\in W_0w_{\la}$ can be written uniquely as
$w=uw_{\la}$ for some $u\in W_0^{\la}$, and moreover
$\ell(w)=\ell(u)+\ell(w_{\la})$. Similarly, each $w'\in
\s_l(w_{\mu})W_n$ can be written uniquely as $w'=\s_l(w_{\mu})u'$
for some $u'\in W_n^{\mu}$, and moreover
$\ell(w')=\ell(\s_l(w_{\mu}))+\ell(u')$.

Using these facts, along with the facts that $\ooo_l^2=\ooo_l$ and
$W_l(q)=W_0(q)$, we have (compare with the proof of
Theorem~\ref{bigone})
\begin{align*}
P_{\la}&\oo P_{\mu}\oo
T_{g_n^{-1}}=\frac{W_{0\la}(q)W_{0\mu}(q)}{W_0^4(q)q_{w_{\la}}q_{w_{\mu}}}\sum_{\substack{
w_1\in W_0w_{\la}W_l\\ w_2\in
W_l\s_l(w_{\mu})W_n}}q_{w_1}q_{w_2}T_{w_1}T_{w_2}\\
&=\frac{W_{0\la}(q)W_{0\mu}(q)W_l^2(q)}{W_0^4(q)q_{w_{\la}}q_{w_{\mu}}}\left(\sum_{u_1\in
W_{0}^{\la}}q_{u_1w_{\la}}T_{u_1w_{\la}}\right)\ooo_l^2\left(\sum_{u_2\in
W_n^{\mu}}q_{\s_l(w_{\mu})u_2}T_{\s_l(w_{\mu})u_2}\right)\\
&=\frac{W_{0\la}(q)W_{0\mu}(q)}{W_0^2(q)q_{w_{\la}}q_{w_{\mu}}}\left(\sum_{w\in
W_0w_{\la}}q_{w}T_{w}\right)\ooo_l\left(\sum_{w'\in \s_l(w_{\mu})W_n}q_{w'}T_{w'}\right)\\
&=\frac{W_{0\la}(q)W_{0\mu}(q)}{W_0^3(q)q_{w_{\la}}q_{w_{\mu}}}\sum_{\substack{w_1\in
W_0w_{\la},\,w_2\in W_l\\
w_3\in\s_l(w_{\mu})W_n}}q_{w_1}q_{w_2}q_{w_3}T_{w_1}T_{w_2}T_{w_3}\,.
\end{align*}

It is simple to see that
\begin{align*}
\sum_{\substack{w_1\in
W_0w_{\la},\,w_2\in W_l\\
w_3\in\s_l(w_{\mu})W_n}}q_{w_1}q_{w_2}q_{w_3}T_{w_1}T_{w_2}T_{w_3}=\sum_{w\in
W}d_w(\la,\mu)q_wT_w
\end{align*}
where $d_w(\la,\mu)$ is a linear combination of products of
$d_{w_1,w_2;w_3}'$'s with nonnegative integer coefficients, and so
\begin{align*}
P_{\la}\oo P_{\mu}\oo
T_{g_n^{-1}}&=\frac{W_{0\la}(q)W_{0\mu}(q)}{W_0^3(q)q_{w_{\la}}q_{w_{\mu}}}\sum_{w\in
W}d_{w}(\la,\mu)q_wT_w.
\end{align*}
So the coefficient of $T_{w_{\nu}}$ when $P_{\la}\oo P_{\mu}\oo
T_{g_n^{-1}}$ is expanded in terms of the $T_w$'s is
\begin{align}\label{secondeval}
\frac{W_{0\la}(q)W_{0\mu}(q)}{W_0^3(q)}\frac{q_{w_{\nu}}}{q_{w_{\la}}q_{w_{\mu}}}d_{w_{\nu}}(\la,\mu).
\end{align}
Comparing (\ref{secondeval}) with (\ref{andanother}) we see that
$c_{\la,\mu;\nu}'=d_{w_{\nu}}(\la,\mu)$, and so the result
follows from Lemma~\ref{neat} and the fact that
$d_{w_{\nu}}(\la,\mu)$ is a linear combination of products of
$d_{w_1,w_2;w_3}'$'s with nonnegative integer coefficients.
\end{proof}

\newpage

\section*{Appendix: Parameter Systems of Regular Affine Buildings}

For an $\widetilde{X}_n$ building there $n+1$ vertices in the
Coxeter graph. The special vertices are marked with an $s$. If
all of the parameters are equal we write $q_i=q$.

\begin{figure}[ht]
 \begin{center}
 \psset{xunit= 0.85 cm,yunit= 0.85 cm}
 \psset{origin={0,0}}
\vspace{5cm}

\rput(-4,5.5){$\widetilde{A}_1$:}\pscircle*(-2,5.5){2pt}
\pscircle*(-1,5.5){2pt}\rput(-1.5,5.7){$\infty$}
\psline(-2,5.5)(-1,5.5)\rput(-2,5.75){$q$}\rput(-1,5.75){$q$}
\rput(-2,5.3){$s$}\rput(-1,5.3){$s$} \rput(1.5,5.5){$\bct_1$:}
\pscircle*(3,5.5){2pt}\pscircle*(4,5.5){2pt}\rput(3.5,5.7){$\infty$}
\psline(3,5.5)(4,5.5)\rput(3,5.75){$q_0$}\rput(4,5.75){$q_1$}
\rput(3,5.3){$s$}\rput(4,5.3){$s$}
\rput(-4,3.5){$\widetilde{A}_n(n\geq2)$:}
\rput(-2,3.8){$q$}\rput(-2,3.3){$s$}\pscircle*(-2,3.5){2pt}
\rput(-1,3.8){$q$}\rput(-1,3.3){$s$}\pscircle*(-1,3.5){2pt}
\rput(1,3.8){$q$}\rput(1,3.3){$s$}\pscircle*(1,3.5){2pt}
\rput(2,3.8){$q$}\rput(2,3.3){$s$}\pscircle*(2,3.5){2pt}
\rput(0,4.8){$q$}\rput(0,4.3){$s$}\pscircle*(0,4.5){2pt}
\pscircle*(0,3.5){1.5pt} \pscircle*(0.4,3.5){1.5pt}
\pscircle*(-0.4,3.5){1.5pt} \psline(0,4.5)(-2,3.5)
\psline(0,4.5)(2,3.5) \psline(-2,3.5)(-1,3.5)
\psline(1,3.5)(2,3.5) \rput(-4,2){$\widetilde{B}_n(n\geq3)$:}
\rput(-1,2.3){$q_0$}\pscircle*(-1,2){2pt}
\rput(-2.2,2.5){$s$}\rput(-2,2.8){$q_0$}\pscircle*(-2,2.5){2pt}
\rput(-2.2,1.5){$s$}\rput(-2,1.8){$q_0$}\pscircle*(-2,1.5){2pt}
\psline(-1,2)(-2,2.5) \psline(-1,2)(-2,1.5)
\rput(0,2.3){$q_0$}\pscircle*(0,2){2pt} \psline(-1,2)(0,2)
\pscircle*(0.6,2){1.5pt} \pscircle*(1.4,2){1.5pt}
\pscircle*(1,2){1.5pt} \rput(2,2.3){$q_0$}\pscircle*(2,2){2pt}
\rput(3,2.3){$q_0$}\pscircle*(3,2){2pt} \psline(2,2)(4,2)
\rput(4,2.3){$q_n$}\pscircle*(4,2){2pt} \rput(3.5,2.3){$4$}
\rput(-4,0.5){$\widetilde{C}_n(n\geq2)$:}
\rput(-2.2,0.5){$s$}\rput(-2,0.8){$q_0$}\pscircle*(-2,0.5){2pt}
\rput(-1,0.8){$q_1$}\pscircle*(-1,0.5){2pt}
\psline(-2,0.5)(-1,0.5) \rput(-1.5,0.8){$4$}
\rput(0,0.8){$q_1$}\pscircle*(0,0.5){2pt} \psline(-1,0.5)(0,0.5)
\pscircle*(1.4,0.5){1.5pt} \pscircle*(0.6,0.5){1.5pt}
\pscircle*(1,0.5){1.5pt} \rput(2,0.8){$q_1$}\pscircle*(2,0.5){2pt}
\rput(3,0.8){$q_1$}\pscircle*(3,0.5){2pt} \psline(2,0.5)(4,0.5)
\rput(4.2,0.5){$s$}\rput(4,0.8){$q_0$}\pscircle*(4,0.5){2pt}
\rput(3.5,0.8){$4$}
\rput(-4,-0.5){$\bct_n(n\geq2)$:}\rput(-2.2,-0.5){$s$}\rput(-2,-0.2){$q_0$}\pscircle*(-2,-0.5){2pt}
\rput(-1,-0.2){$q_1$}\pscircle*(-1,-0.5){2pt}
\psline(-2,-0.5)(-1,-0.5) \rput(-1.5,-0.2){$4$}
\rput(0,-0.2){$q_1$}\pscircle*(0,-0.5){2pt}
\psline(-1,-0.5)(0,-0.5) \pscircle*(1.4,-0.5){1.5pt}
\pscircle*(0.6,-0.5){1.5pt} \pscircle*(1,-0.5){1.5pt}
\rput(2,-0.2){$q_1$}\pscircle*(2,-0.5){2pt}
\rput(3,-0.2){$q_1$}\pscircle*(3,-0.5){2pt}
\psline(2,-0.5)(4,-0.5)
\rput(4.2,-0.5){$s$}\rput(4,-0.2){$q_n$}\pscircle*(4,-0.5){2pt}
\rput(3.5,-0.2){$4$} \rput(-4,-2){$\widetilde{D}_n(n\geq4)$:}
\rput(-1,-1.7){$q$}\pscircle*(-1,-2){2pt}
\rput(-2.2,-1.5){$s$}\rput(-2,-1.2){$q$}\pscircle*(-2,-1.5){2pt}
\rput(-2.2,-2.5){$s$}\rput(-2,-2.2){$q$}\pscircle*(-2,-2.5){2pt}
\psline(-1,-2)(-2,-1.5) \psline(-1,-2)(-2,-2.5)
\rput(0,-1.7){$q$}\pscircle*(0,-2){2pt} \psline(-1,-2)(0,-2)
\pscircle*(0.6,-2){1.5pt} \pscircle*(1.4,-2){1.5pt}
\pscircle*(1,-2){1.5pt} \rput(2,-1.7){$q$}\pscircle*(2,-2){2pt}
\rput(3,-1.7){$q$}\pscircle*(3,-2){2pt} \psline(2,-2)(3,-2)
\rput(4.2,-1.5){$s$}\rput(4,-1.2){$q$}\pscircle*(4,-1.5){2pt}
\rput(4.2,-2.5){$s$}\rput(4,-2.2){$q$}\pscircle*(4,-2.5){2pt}
\psline(3,-2)(4,-1.5) \psline(3,-2)(4,-2.5)
\rput(-4,-5){$\widetilde{E}_6$:}
\rput(-2.2,-5){$s$}\rput(-2,-4.7){$q$}\pscircle*(-2,-5){2pt}
\rput(-1,-4.7){$q$}\pscircle*(-1,-5){2pt}
\rput(-0.2,-4.7){$q$}\pscircle*(0,-5){2pt}
\rput(-0.2,-4){$q$}\pscircle*(0,-4){2pt}
\rput(0,-2.8){$s$}\rput(-0.2,-3){$q$}\pscircle*(0,-3){2pt}
\rput(1,-4.7){$q$}\pscircle*(1,-5){2pt}
\rput(2.2,-5){$s$}\rput(2,-4.7){$q$}\pscircle*(2,-5){2pt}
\psline(-2,-5)(2,-5) \psline(0,-5)(0,-3)
\rput(-4,-6.5){$\widetilde{E}_7$:}
\rput(-2.2,-6.5){$s$}\rput(-2,-6.2){$q$}\pscircle*(-2,-6.5){2pt}
\rput(-1,-6.2){$q$}\pscircle*(-1,-6.5){2pt}
\rput(0,-6.2){$q$}\pscircle*(0,-6.5){2pt}
\rput(0.8,-6.2){$q$}\pscircle*(1,-6.5){2pt}
\rput(0.8,-5.5){$q$}\pscircle*(1,-5.5){2pt}
\rput(2,-6.2){$q$}\pscircle*(2,-6.5){2pt}
\rput(3,-6.2){$q$}\pscircle*(3,-6.5){2pt}
\rput(4.2,-6.5){$s$}\rput(4,-6.2){$q$}\pscircle*(4,-6.5){2pt}
\psline(-2,-6.5)(4,-6.5) \psline(1,-6.5)(1,-5.5)
\rput(-4,-8){$\widetilde{E}_8$:}
\rput(-2,-7.7){$q$}\pscircle*(-2,-8){2pt}
\rput(-1,-7.7){$q$}\pscircle*(-1,-8){2pt}
\rput(-0.2,-7.7){$q$}\pscircle*(0,-8){2pt}
\rput(-0.2,-7){$q$}\pscircle*(0,-7){2pt}
\rput(1,-7.7){$q$}\pscircle*(1,-8){2pt}
\rput(2,-7.7){$q$}\pscircle*(2,-8){2pt}
\rput(3,-7.7){$q$}\pscircle*(3,-8){2pt}
\rput(4,-7.7){$q$}\pscircle*(4,-8){2pt}
\rput(5.2,-8){$s$}\rput(5,-7.7){$q$}\pscircle*(5,-8){2pt}
\psline(-2,-8)(5,-8) \psline(0,-8)(0,-7)
\rput(-4,-9.5){$\widetilde{F}_4$:}
\rput(-2.2,-9.5){$s$}\rput(-2,-9.2){$q_0$}\pscircle*(-2,-9.5){2pt}
\rput(-1,-9.2){$q_0$}\pscircle*(-1,-9.5){2pt}
\rput(0,-9.2){$q_0$}\pscircle*(0,-9.5){2pt}
\rput(1,-9.2){$q_4$}\pscircle*(1,-9.5){2pt}
\rput(2,-9.2){$q_4$}\pscircle*(2,-9.5){2pt}
\psline(-2,-9.5)(2,-9.5) \rput(0.5,-9.2){$4$}
\rput(-4,-11){$\widetilde{G}_2$:}
\rput(-2.2,-11){$s$}\rput(-2,-10.7){$q_0$}\pscircle*(-2,-11){2pt}
\rput(-1,-10.7){$q_0$}\pscircle*(-1,-11){2pt}
\rput(0,-10.7){$q_1$}\pscircle*(0,-11){2pt} \psline(-2,-11)(0,-11)
\rput(-0.5,-10.7){$6$}

\end{center}
\end{figure}

\vspace{9.5cm}

\bibliography{PhD.bib}
\bibliographystyle{plain}

\end{document}